\documentclass{amsart}

\usepackage{enumerate}
\usepackage{amsmath,xspace,amssymb,mathrsfs}
\input xy
\xyoption{all}
\xyoption{2cell}
\UseAllTwocells
\CompileMatrices

\title{A 2-categories companion}
\author{Stephen Lack}
\address{School of Computing and Mathematics\\
University of Western Sydney}
\email{s.lack@uws.edu.au}
\thanks{The support of the Australian Research Council and
DETYA is gratefully acknowledged.}

\renewcommand{\phi}{\varphi}
\renewcommand{\epsilon}{\varepsilon}

\newcommand{\A}{{\ensuremath{\mathscr A}}\xspace}
\newcommand{\B}{{\ensuremath{\mathscr B}}\xspace}
\newcommand{\C}{{\ensuremath{\mathscr C}}\xspace}
\newcommand{\D}{{\ensuremath{\mathscr D}}\xspace}
\newcommand{\E}{{\ensuremath{\mathscr E}}\xspace}

\newcommand{\G}{{\ensuremath{\mathscr G}}\xspace}
\renewcommand{\H}{{\ensuremath{\mathscr H}}\xspace}
\newcommand{\I}{{\ensuremath{\mathscr I}}\xspace}
\newcommand{\K}{{\ensuremath{\mathscr K}}\xspace}
\newcommand{\LL}{{\ensuremath{\mathscr L}}\xspace}
\newcommand{\M}{{\ensuremath{\mathcal M}}\xspace}
\renewcommand{\P}{{\ensuremath{\mathscr P}}\xspace}

\newcommand{\V}{{\ensuremath{\mathscr V}}\xspace}
\newcommand{\W}{{\ensuremath{\mathscr W}}\xspace}

\newcommand{\Kf}{\ensuremath{\K_f}\xspace}

\newcommand{\DD}{{\ensuremath{\mathbf{\Delta}}}\xspace}
\newcommand{\aA}{{\ensuremath{\mathbb A}}\xspace}
\newcommand{\BB}{{\ensuremath{\mathbb B}}\xspace}
\newcommand{\EE}{{\ensuremath{\mathbb E}}\xspace}
\newcommand{\NN}{{\ensuremath{\mathbb N}}\xspace}
\newcommand{\TT}{{\ensuremath{\mathbb T}}\xspace}

\newcommand{\Rel}{\textnormal{\bf Rel}\xspace}
\newcommand{\RelE}{\textnormal{\bf Rel}(\EE)\xspace}
\newcommand{\Par}{\textnormal{\bf Par}\xspace}
\newcommand{\Mod}{\textnormal{\bf Mod}\xspace}
\newcommand{\VMod}{\textnormal{\V-{\bf Mod}}\xspace}
\newcommand{\Mat}{\textnormal{\bf Mat}\xspace}
\newcommand{\VMat}{\textnormal{\V-{\bf Mat}}\xspace}
\newcommand{\Span}{\textnormal{\bf Span}\xspace}
\newcommand{\SpanE}{\textnormal{\bf Span}(\EE)\xspace}
\newcommand{\RMod}{\textnormal{\bf $R$-Mod}\xspace}
\newcommand{\SMod}{\textnormal{\bf $S$-Mod}\xspace}

\newcommand{\MonCat}{\textnormal{\bf MonCat}\xspace}
\newcommand{\OpMonCat}{\textnormal{\bf OpMonCat}\xspace}
\newcommand{\Hom}{\textnormal{\bf Hom}\xspace}

\newcommand{\Set}{\textnormal{\bf Set}\xspace}

\newcommand{\SSet}{\textnormal{\bf SSet}\xspace}
\newcommand{\sS}{\SSet}
\newcommand{\Ab}{\textnormal{\bf Ab}\xspace}
\newcommand{\Grp}{\textnormal{\bf Grp}\xspace}
\newcommand{\ordf}{\ensuremath{\textnormal{\bf Ord}_{\textbf{f}}}\xspace}
\newcommand{\Cat}{\textnormal{\bf Cat}\xspace}
\newcommand{\VCat}{\textnormal{\V-{\bf Cat}}\xspace}
\newcommand{\CatE}{\textnormal{\bf Cat}(\EE)\xspace}

\newcommand{\Mnd}{\ensuremath{\textnormal{\bf Mnd}(\K)}\xspace}
\newcommand{\Mndf}{\ensuremath{\textnormal{\bf Mnd}_f(\K)}\xspace}
\newcommand{\Mndfc}{\ensuremath{\textnormal{\bf Mnd}_f(\Cat)}\xspace}
\newcommand{\Mndfcop}{\ensuremath{\textnormal{\bf Mnd}_f(\Cat)\op}\xspace}

\newcommand{\Endf}{\ensuremath{\textnormal{\bf End}_f(\K)}\xspace}

\newcommand{\bicat}{\textnormal{\bf Bicat}\xspace}
\newcommand{\twocatps}{\ensuremath{\textnormal{\bf 2-Cat}_\textnormal{ps}}\xspace}
\newcommand{\twocat}{\textnormal{\bf 2-Cat}\xspace}
\newcommand{\Twocat}{\textnormal{\bf 2-CAT}\xspace}
\newcommand{\NHom}{\textnormal{\bf NHom}\xspace}

\newcommand{\talg}{{\ensuremath{\textnormal{$T$-Alg}}}\xspace}
\newcommand{\talgs}{{\ensuremath{\textnormal{$T$-Alg}_s}}\xspace}
\newcommand{\talgc}{{\ensuremath{\textnormal{$T$-Alg}_c}}\xspace}
\newcommand{\talgl}{{\ensuremath{\textnormal{$T$-Alg}_\ell}}\xspace}
\newcommand{\pstalg}{{\ensuremath{\textnormal{Ps-$T$-Alg}}}\xspace}

\newcommand{\ttalg}{{\ensuremath{\textnormal{$T'$-Alg}}}\xspace}
\newcommand{\salg}{{\ensuremath{\textnormal{$S$-Alg}}}\xspace}
\newcommand{\salgs}{{\ensuremath{\textnormal{$S$-Alg}_s}}\xspace}

\newcommand{\malgs}{{\ensuremath{\textnormal{$M$-Alg}_s}}\xspace}

\newcommand{\nalgs}{{\ensuremath{\textnormal{$N$-Alg}_s}}\xspace}
\newcommand{\fxalgs}{{\ensuremath{\textnormal{$FX$-Alg}_s}}\xspace}
\newcommand{\EM}{\textnormal{\bf EM}\xspace}
\newcommand{\KL}{\textnormal{\bf KL}\xspace}
\newcommand{\alg}{\textnormal{\bf alg}\xspace}
\newcommand{\mnd}{\textnormal{\bf mnd}\xspace}

\newcommand{\pslim}{\textnormal{\sf pslim}\xspace}
\newcommand{\colim}{\textnormal{\sf colim}\xspace}
\newcommand{\Lan}{\textnormal{\sf Lan}\xspace}
\newcommand{\lan}{\textnormal{\sf lan}\xspace}
\newcommand{\ran}{\textnormal{\sf ran}\xspace}
\newcommand{\sem}{\textnormal{\sf sem}\xspace}

\newcommand{\id}{\textnormal{\sf id}\xspace}
\newcommand{\comp}{\textnormal{\sf comp}\xspace}
\newcommand{\cosk}{\textnormal{\sf cosk}\xspace}

\newcommand{\Gray}{\textnormal{\sf Gray}\xspace}
\newcommand{\Psm}{\textnormal{\sf Psm}\xspace}
\newcommand{\Psa}{\textnormal{\sf Psa}\xspace}
\newcommand{\Ps}{\textnormal{\sf Ps}\xspace}
\newcommand{\ob}{\textnormal{\sf ob}\xspace}
\newcommand{\iso}{\textnormal{iso}\xspace}

\newcommand{\hto}{\ensuremath{\,\mathaccent\shortmid\rightarrow\,}}
\newcommand{\mono}{\rightarrowtail}
\newcommand{\too}[1][]{\ensuremath{\overset{#1}{\longrightarrow}}}

\newcommand{\ct}{\pitchfork}
\newcommand{\lhom}{\setminus}

\newcommand{\two}{\ensuremath{{\hbox{\textrm 2}\kern-.25em
        \hbox{\vrule height1.5ex width 0.4pt depth -.2ex}}\kern.2em}\xspace}
\newcommand{\three}{\ensuremath{{\hbox{\textrm 3}\kern-.25em
        \hbox{\vrule height1.5ex width 0.4pt depth -.2ex}}\kern.2em}\xspace}

\newcommand{\<}{\langle}
\renewcommand{\>}{\rangle}

\newcommand{\inv}{\ensuremath{^{-1}}}
\newcommand{\ch}{_{\textnormal{ch}}}
\newcommand{\op}{\ensuremath{^{\textnormal{op}}}}
\newcommand{\co}{\ensuremath{^{\textnormal{co}}}}
\newcommand{\coop}{\ensuremath{^{\textnormal{coop}}}}
\newcommand{\rev}{\ensuremath{^{\textnormal{rev}}}}
\newcommand{\st}{\ensuremath{_{\textnormal{st}}}}

\newcommand{\ot}{\otimes}
\renewcommand{\t}{\times}
\newcommand{\sd}{\rtimes}

\renewcommand{\bar}{\overline}
\newcommand{\Ho}{\textrm{Ho}}
\newcommand{\maps}{\colon}

\newcommand{\gcomp}{\textsf{comp}\xspace}
\newcommand{\garr}{\textsf{arr}\xspace}
\newcommand{\gob}{\textsf{ob}\xspace}
\newcommand{\gtriple}{\textsf{triple}\xspace}
\newcommand{\giso}{\textsf{iso}\xspace}

\newcommand{\fbar}{\ensuremath{\bar{f}}}
\newcommand{\gbar}{\ensuremath{\bar{g}}}

\newcommand{\LaxtwoK}{\ensuremath{\mathbf{Lax}(\two,\K)}\xspace}

\newtheorem{theorem}{Theorem}[section]
\newtheorem{lemma}[theorem]{Lemma}
\newtheorem{corollary}[theorem]{Corollary}
\newtheorem{proposition}[theorem]{Proposition}

\theoremstyle{definition}

\newtheorem{example}[theorem]{Example}
\newtheorem{xca}[theorem]{Exercise}

\theoremstyle{remark}
\newtheorem{remark}[theorem]{Remark}

\newcounter{jump}
\newcommand{\saveenumi}{\setcounter{jump}{\value{enumi}}}
\newcommand{\recall}{\addtocounter{enumi}{\value{jump}}}

\begin{document}

\begin{abstract}
This paper is a rather informal guide to some of the basic theory
of 2-categories and bicategories, including notions of limit and 
colimit, 2-dimensional universal algebra, formal category theory,
and nerves of bicategories. 
\end{abstract}

\maketitle

\section{Overview and basic examples}

This paper is a rather informal guide to some of the basic theory
of 2-categories and bicategories, including notions of limit and 
colimit, 2-dimensional universal algebra, formal category theory,
and nerves of bicategories. As is the way of these things, the choice
of topics is somewhat personal. No attempt is made at either rigour or
completeness. Nor is it completely introductory: you will not find a
definition of bicategory; but then nor will you really need one to
read it. In keeping with the philosophy of category theory, the
morphisms between bicategories play more of a role than the
bicategories themselves.

\subsection{The key players}
\label{sec:introduction}

There are bicategories, 2-categories, and \Cat-categories.  
The latter two are exactly the same (except that strictly speaking
a \Cat-category should have small hom-categories, but that need not
concern us here).
The first two are nominally different --- the 2-categories
are the strict bicategories, and not every bicategory is strict --- 
but every bicategory is {\em biequivalent} to a strict one, and 
biequivalence is the right general notion of equivalence for
bicategories and for 2-categories.
Nonetheless, the theories of bicategories, 2-categories, and
\Cat-categories have rather different flavours. 

An enriched category is a category in which the hom-functors
take their values not in \Set, but in some other category \V. 
The theory of enriched categories is now very well developed, and 
\Cat-category theory is the special case where $\V=\Cat$. In \Cat-category
theory one deals with higher-dimensional versions of the usual notions 
of functor, limit, monad, and so on, without any ``weakening''. 
The passage from category theory to \Cat-category theory is well understood;
unfortunately \Cat-category theory is generally not what one wants to
do --- it is too strict, and fails to deal with the notions that arise
in practice.

In bicategory theory all of these notions are weakened. One never says
that arrows are equal, only isomorphic, or even sometimes only that there is
a comparison 2-cell between them. If one wishes to 
generalize a result about categories to bicategories, it is generally
clear in principle what should be done, but the details can be 
technically very difficult.

2-category theory is a ``middle way'' between \Cat-category theory and
bicategory theory. It {\em uses} enriched category theory, but not in
the simple minded way of \Cat-category theory; and it cuts through some
of the technical nightmares of bicategories. The prefix ``2-'', as in
2-functor or 2-limit, will always denote the strict notion; although
often we will use it to describe or analyze non-strict phenomena.

There are also various other related notions, which will be less important
in this companion. \sS-categories are categories
enriched in simplicial sets; every 2-category induces an \sS-category,
by taking nerves of the hom-categories. Double categories are internal
categories in \Cat. Once again every 2-category can be seen as a double
category. A slight generalization of double categories allows bicategories
to fit into this picture. Finally there are the internal categories in 
\sS; both \sS-categories and double categories can be seen as special cases
of these.

\subsection{Nomenclature and symbols}

In keeping with our general policy, the word {\em 2-functor} is
understood in the strict sense: a 2-functor between 2-categories \A
and \B assigns objects to objects, morphisms to morphisms, and 2-cells
to 2-cells, preserving all of the 2-category structure strictly. 
We shall of course want to consider more general types of morphism
between 2-categories later on.

If ``widget'' is the name of some particular categorical structure,
then there are various systems of nomenclature for weak 2-widgets.
Typically one speaks of {\em pseudo} widgets for the up-to-isomorphism notion, 
{\em lax} widgets for the up-to-not-necessarily-invertible comparison notion, 
and when the direction of the comparison is reversed, either {\em oplax} widget
or {\em colax} widget, depending on the specific case. But there are also
other conventions. In contexts where the pseudo notion is most important,
this is called simply a widget, and then one speaks explicitly of 
{\em strict} widgets in the strict case. In contexts where the lax notion
is most important (such as with monoidal functors), it is this which
has no prefix; and one has {\em strict} widgets in the strict case or
{\em strong} widgets in the pseudo.

As we move up to 2-categories and higher categories, there are various
notions of sameness, having the following symbols:
\begin{itemize}
\item $=$ is equality
\item $\cong$ is isomorphism (morphisms $f$ and $g$ with $gf=1$, $fg=1$)
\item $\simeq$ is equivalence ($gf\cong 1$, $fg\cong 1$) 
\item $\sim$ is sometimes used for biequivalence.
\end{itemize}

In Sections~\ref{sect:examples} and~\ref{sect:examples-bicats} 
we look at various examples of 2-categories 
and bicategories. The separation between the 2-category examples and 
the bicategory examples is not really about strictness but about
the sort of morphisms involved. The 2-category examples involve
functions or functors of some sort; the bicategory examples (except
the case of a monoidal category) involve more general types of
morphism such as relations. These ``non-functional'' morphisms are
often depicted using a slashed arrow ($\hto$) rather than an
ordinary one ($\to$). Typically the functional morphisms can be seen 
as a special case of the non-functional ones. Sometimes it is also
possible to characterize the non-functional ones as a special 
type of functional morphism (with different domain and/or codomain),
and this can provide a concrete construction of a 2-category
biequivalent to the given bicategory. The other special 
type of arrow often used is a ``wobbly'' one ($\rightsquigarrow$);
this denotes a weak (pseudo, lax, etc.) morphism.

\subsection{Contents}
\label{sect:contents}

In the remainder of this section we look at examples of 2-categories
and bicategories. In Section~\ref{sect:formal} we begin the study
of formal category theory, including adjunctions, extensions,
and monads, but stopping short of the full-blown formal theory of
monads. In Section~\ref{sect:morphisms} we look at various types 
of morphism between bicategories or 2-categories: strict, pseudo,
lax, partial; and see how these can be used to describe enriched
and indexed categories. In Section~\ref{sect:2DUA} we begin the 
study of 2-dimensional universal algebra, with the basic definitions
and the construction of weak morphism classifiers. This is continued
in Section~\ref{sect:presentations} on presentations for 2-monads, 
which demonstrates how various categorical structures can be described
using 2-monads. Section~\ref{sect:limits} looks at various 2-categorical
and bicategorical notions of limit and considers their existence in the
2-categories of algebras for 2-monads. Section~\ref{sect:Quillen} is 
about aspects of Quillen model structures related to 2-categories and
to 2-monads. In Section~\ref{sect:ftm} we return to the formal theory
of monads, applying some of the earlier material on limits. 
Section~\ref{sect:pseudomonads} looks at the formal theory of 
{\em pseudomonads}, developed in a Gray-category. Section~\ref{sect:nerves}
looks at notions of nerve for bicategories. There are relatively few
references throughout the text, but at the end of each section there
is a brief commented bibliography.

\subsection{Examples of 2-categories}
\label{sect:examples}

\Cat\ is the mother of all 2-categories, just as \Set\ is the mother of
all categories.  From many points of view, it has all the best
properties as a 2-category (but not as a category: for example
colimits in \Cat are not stable under pullback).

A small category involves a set of objects and a set of arrows, and 
also hom-sets between any two objects. One can generalize the
notion of category in various ways by replacing various of these
sets by objects of some other category.

\begin{enumerate}[(a)]
\item If \V is a monoidal category one can consider the 2-category
\VCat\ of categories enriched in \V%\cite{closed}
; these have \V-valued 
hom-objects rather than hom-sets. The theory works best %\cite{Kelly-book} 
when \V is symmetric monoidal closed, complete, and
cocomplete. As for examples of enriched categories, one has ordinary 
categories ($\V=\Set$), additive categories ($\V=\Ab$), 2-categories 
($\V=\Cat$), preorders ($\V=\two$, the ``arrow category''),
simplicially enriched categories ($\V=\SSet$), and DG-categories (\V
the category of chain complexes).
\item More generally still, one can consider a bicategory \W as a
  many-object version of a monoidal category; there is a corresponding
notion of \W-enriched category: see \cite{BCSW} or 
Section~\ref{sect:lax-functors}. Sheaves on a site can be
described as \W-categories for a suitable choice of \W.
\item If \EE is a category with finite limits, one can consider 
the 2-category \CatE\ of categories internal to \EE; these have
an \EE-object of objects and an \EE-object of morphisms. The theory 
works better the better the category \EE; the cases of a topos or
an abelian category are particularly nice. This includes ordinary
categories ($\EE=\Set$), double categories ($\EE=\Cat$), morphisms
of abelian groups ($\EE=\Ab$), and crossed modules ($\EE=\Grp$).
\saveenumi
\end{enumerate}

There is another class of examples, in which the objects are ``categories 
with structure''. The structure could be something like
\begin{enumerate}[(a)]
\recall
\item category with finite products
\item category with finite limits \label{item:lex}
\item monoidal category 
\item topos
\item category with finite products and coproducts and a distributive law
\saveenumi
\end{enumerate}
For most of these there are also analogues involving enriched or internal 
categories with the relevant structure.

In each case you need to decide which morphisms to use. Normally
you don't want the strictly algebraic ones (preserving the structure 
on the nose): although they
can be technically useful, they are rare in nature.
More common are the ``pseudo'' morphisms: these are functors preserving 
the structure ``up to (suitably
coherent) isomorphism''. In (\ref{item:lex}), for example, this 
would correspond to the usual
notion of finite-limit-preserving functor.

Sometimes, however, it's good to consider an even weaker notion 
of morphism, as in the 2-category 
\MonCat of  monoidal categories, monoidal functors, and 
monoidal natural transformations. Monoidal functors are the
``lax'' notion, involving maps $FA\ot FB\to F(A\ot B)$, coherent, but not
necessarily invertible. 
Here are some reasons you might like this level of generality:
\begin{itemize}
\item Consider the monoidal categories \Ab of abelian groups, with 
the usual tensor product, and \Set of sets, with the cartesian product.
The forgetful functor $U$ from \Ab\ to \Set definitely does not
preserve this structure, but we have the universal bilinear map 
$UG\times UH \to U(G\ot H)$, and this makes $U$ into a monoidal functor.
\item A monoidal functor $\V\to\W$ sends monoids in \V to monoids in \W,
via the rule
$$\xymatrix @R1pc {
M\ot M \ar[dd]_{m} && FM\ot FM \ar[d] \\
& \mapsto & F(M\ot M) \ar[d]^{Fm} \\
M && FM }$$
\item Suppose \V and \W are monoidal categories and $F:\V\to\W$ is a
  left adjoint which does preserve the monoidal structure up to
  coherent isomorphism. There is no reason why the right adjoint $U$
  should do so, but there will be induced comparison maps 
  $UA\ot UB\to U(A\ot B)$ making $U$ a monoidal functor. (Think of the 
  tensor product as a type of colimit, so the
  left adjoint preserves it, but the right adjoint doesn't
  necessarily.) In fact the monoidal functor $U:\Ab\to\Set$ arises in
  this way.
\end{itemize}

The case of monoidal categories is typical. Given an adjunction
$F\dashv U$ between categories \A and \B with algebraic structure, to 
make the right adjoint $U$ a colax morphism is equivalent to making the 
left adjoint $F$ lax, while if the whole adjunction lives within the world
of lax morphisms, then $F$ is not just lax but pseudo. This phenomenon
is called {\em doctrinal adjunction} \cite{Kelly-doctrinal}.

For a further example, consider the structure of categories with
finite coproducts. For a functor $F:\A\to\B$ between categories
with finite coproducts there are canonical comparison maps 
$FA+FB\to F(A+B)$, and these make {\em every} such functor uniquely
into a lax morphism; it is a pseudo morphism exactly when it
preserves the coproducts in the usual sense. Thus in this case
every adjunction between categories with finite coproducts lives
in the lax world, and the fact that the left adjoint is actually 
pseudo reduces to the well known fact that left adjoints preserve
coproducts.

In the case of categories with finite products or finite limits,
however, the lax morphisms are the same as the pseudo morphisms;
they are just the functors preserving the products or limits in 
the usual sense.

\subsection{Examples of bicategories}
\label{sect:examples-bicats}

Any monoidal category \V determines a one-object bicategory 
$\Sigma\V$ whose morphisms are the objects of \V, and whose 2-cells
are the morphisms of \V. The tensor product of \V is the (horizontal)
composition in $\Sigma\V$.

\begin{enumerate}[(a)]
\recall
\item \Rel consists of sets and relations. The objects are sets and
the morphisms $X\hto Y$ are the relations from $X$ to $Y$; that is,
the monomorphisms $R\mono X\t Y$. 
This bicategory is `locally posetal', in the sense that for any two
parallel 1-cells, there is at most one 2-cell between them. There is
a 2-cell from $R$ to $S$ if and only if $R$ is contained in $S$ as 
a subobject of $X\t Y$; in other words, if there is a morphism $R\to S$
making the triangles in 
$$
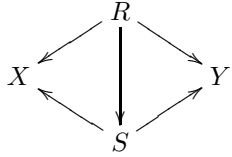
\xymatrix @R 1pc {
& R \ar[dd] \ar[dl] \ar[dr] \\
X && Y \\
& S \ar[ul] \ar[ur] }$$
commute. As usual, $xRy$ means that $(x,y)\in R$.
The composite of $R\mono X\t Y$ and $S\mono Y\t Z$ is the relation
$R\circ S$ defined by 
$$x(R\circ S)z \iff (\exists y) xRySz.$$
We get a
2-category biequivalent to this one by identifying isomorphic 1-cells;
this works for any locally posetal 2-category.

Another 2-category biequivalent to \Rel has sets for objects, and 
as morphisms from $X$ to $Y$ the join-preserving maps from $\P X$ to 
$\P Y$, where $\P X$ denotes the set of all subsets of $X$. Here a
relation $R$ is represented by the function sending a subset $U\subset
X$ to $\{y\in Y: (\exists x\in U)xRy\}$.

\item 
\Par consists of sets and partial functions.  A partial
function from $X$ to $Y$ is a diagram $X\leftarrowtail D \to Y$ in
\Set,
where $D$ is the domain of definition of the partial function;
2-cells and composition are defined as in \Rel.
Again, we get a biequivalent 2-category by identifying isomorphic
1-cells.

Alternatively, this is biequivalent to the 2-category of pointed 
sets and basepoint-preserving maps, with suitably defined (exercise!)
2-cells.

\item 
\Span consists of sets and ``spans'' $X\leftarrow E \to Y$ in
\Set, with composition by pullback, and with 2-cells given by diagrams
such as
$$
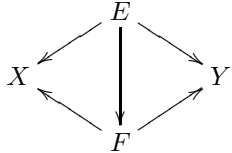
\xymatrix @R 1pc {
& E \ar[dd] \ar[dl] \ar[dr] \\
X && Y \\
& F \ar[ur] \ar[ul]. }$$
Unlike the previous two bicategories, this one
is no longer locally posetal, so to get a biequivalent 2-category we
need to do more than just identify isomorphic 1-cells. There are 
general results asserting that any bicategory is
biequivalent to a 2-category, but in fact naturally occurring bicategories
tend to be biequivalent to naturally occurring 2-categories. In this
case, we can take the 2-category whose objects are sets and whose
morphisms are the left adjoints
$\Set/X \to \Set/Y$. Here the span
$$\xymatrix{
X & E \ar[l]_{u} \ar[r]^{v} & Y }$$
is represented by the left adjoint 
$$\xymatrix{
\Set/X \ar[r]^{u^*} & \Set/E \ar[r]^{v_!} & \Set/Y }$$
given by pulling back along $u$ then composing with $v$.

\item
\Mat has sets as objects, $X\times Y$-indexed families
(``matrices'') of sets as morphisms from $X$ to $Y$, and 2-cells are 
families of functions. Composition of 1-cells is given by matrix 
multiplication: if $A=(A_{xy})$ and $B=(B_{yz})$ then
$$(AB)_{xz} = \sum_y A_{xy}\t B_{yz}.$$
This is biequivalent to \Span, but we'll see below that 
spans and matrices become different when we start to consider 
enrichment and internalization. 
A biequivalent 2-category consists of sets and left adjoints
$\Set^X\to \Set^Y$. (Here $X\times Y\to \Set$ can be seen as
a functor $X\to\Set^Y$, and so, since $\Set^X$ is the free cocompletion
of $X$, as a left adjoint $\Set^X\to\Set^Y$.) This is really just
the same as the construction given for \Span, since 
$\Set/X\simeq\Set^X$; once again, though, when we start to enrich
or internalize, the two pictures diverge.

\item
\Mod has rings as objects, 
left $R$-, right $S$-modules as 1-cells $R\hto S$, and 
homomorphisms as 2-cells.
The composite of modules $R\hto S$ and $S\hto T$ is given
by tensoring over $S$.
A biequivalent 2-category involves adjunctions
$\RMod\rightleftarrows\SMod$.

A ring is the same thing as
an \Ab-category (a category enriched in abelian groups) with only 
one object. The underlying additive group of the ring is the single
hom-object; the multiplication of the ring is the composition.
If we identify rings with the corresponding one-object \Ab-categories,
then a module $R\hto S$ becomes an \Ab-functor $R\to[S\op,\Ab]$

But there is no reason to restrict ourselves to one-object categories,
and there is a bicategory \textrm{\Ab-$\mathbf{Mod}$} whose objects are 
\Ab-categories, and whose 1-cells are \Ab-modules $\A\hto\B$; that
is, \Ab-functors $\A\to[\B\op,\Ab]$.

More generally still, we can replace \Ab by any monoidal category \V
with coequalizers which are preserved by tensoring on either side, and
there is then a bicategory \VMod\ of \V-categories and \V-modules: 
once again, if \A and \B are \V-categories then a \V-module $\A\hto\B$
is a \V-functor $\A\to[\B\op,\V]$, or equivalently a left adjoint
$[\A\op,\V]\to[\B\op,\V]$, (and this last description gives a 2-category).

There's even, if you really want, a version with a bicategory \W rather
than a monoidal category \V.
\saveenumi
\end{enumerate}

Now let's internalize and enrich the other examples.  

\begin{enumerate}[(a)]
\recall
\item If \EE is a {\em regular category}, meaning that any morphism 
factorizes as a strong epimorphism followed by a monomorphism, and the 
strong epimorphisms are stable under pullback, 
then we can form \RelE whose objects are those of \EE\ and 
whose morphisms $X\hto Y$ are monomorphisms $R\mono X\t Y$. To compose
$R:X\hto Y$ and $S:Y\hto Z$ we pullback over $Y$, but the resulting map
into $X\t Z$ need not be monic, so we need the factorization system to 
define composition. It turns out that our assumption that strong epimorphisms
are stable under pullback is precisely what is needed for this composition
to be associative.

\item
Similarly, if \C\ is a category and \M\ is a class of monomorphisms
in \C, then we can look at $\mathbf{Par}(\C,\M)$, defined as above
where the given monomorphism is in \M.  
There are conditions on \M\ you need to make this work well: you 
want to be able to pullback an \M-map by an arbitrary map and obtain an
\M-map, and you want \M to be closed under composition and to contain
the isomorphisms.

\item
If \EE\ has finite limits, we can look at \SpanE\
defined in an obvious way.  You need the pullbacks for composition to
work. You don't need any exactness properties to get a bicategory,
but if you want to get a nice
biequivalent 2-category, you'll need to start making more
assumptions on \EE. It turns out that \SpanE\ plays a crucial role
in internal category: we shall see in Example~\ref{ex:monads-SpanE} below
that an internal category in \EE is the same thing as a monad in \SpanE.

\item \label{item:VMat}
\Mat, on the other hand, gets enriched rather than
internalized.  Then \VMat\ 
%(traditionally, things enriched over go in front, while things 
%internalized in go on the right) 
has \emph{sets} as objects and \V-valued matrices $X\times Y\to \V$
as morphisms.  \VMat\ stands in exactly the same relationship to 
\V-categories as \SpanE\ does to categories in \E.
In the case $\V=\Set$ of course \VMat\ is just \Mat, but there 
is also another special case which we have already seen. Let \V
be the arrow-category \two, consisting of two objects 0 and 1, and 
a single non-identity arrow $0\to 1$. This is cartesian closed (a 
\V-category in this case is just a preorder) and \VMat\ in this 
case is \Rel\ (we identify a subject of $X\t Y$ with its characteristic
function, seen as landing in \two).
\saveenumi
\end{enumerate}

\subsection{Duality}
\label{sec:duality}

A bicategory \B\ has not one but three duals:
\begin{itemize}
\item $\B\op$ is obtained by reversing the 1-cells
\item $\B\co$ is obtained by reversing the 2-cells
\item $\B\coop$ is obtained by reversing both
\end{itemize}

In the case of a monoidal category \V, we can form the monoidal
category $\V\op$ by reversing the sense of the morphisms; this
reverses the 2-cells of the corresponding bicategory $\Sigma\V$,
so $\Sigma(\V\op)=(\Sigma\V)\co$. Reversing the 1-cells of $\Sigma\V$
corresponds to reversing the tensor of \V, denoted $\V\rev$, so 
$\Sigma(\V\rev)=(\Sigma\V)\op$.

\subsection{References to the literature}

The basic references for bicategories and 2-categories are
\cite{bicategories}, \cite{Gray}, \cite{Kelly-Street}, and \cite{FibBic}.
The  basic references for enriched categories are \cite{closed},
\cite{Kelly-book}, and \cite{Lawvere-metric}. For a good
example of simplicially-enriched category theory that is
very close to 2-category theory, see \cite{Cordier-Porter}.
Both 2-categories and double categories were first defined
by Ehresmann (see perhaps \cite{Ehresmann}); bicategories were first defined
by B\'enabou \cite{bicategories}. For (a generalization of) the
fact that every bicategory is biequivalent to a 2-category, see
\cite{MacLane-Pare}.

Categories enriched in a bicategory were first defined by Walters
to deal with the example of sheaves on a space (or site) 
\cite{Walters-sheaves1,Walters-sheaves2}.
A good general reference is \cite{BCSW}. 

The importance of monoidal functors (not necessarily strong)
was observed both by Eilenberg-Kelly \cite{closed} and by 
B\'enabou \cite{bicategories}. 

For doctrinal adjunction see \cite{Kelly-doctrinal}.

\section{Formal category theory}
\label{sect:formal}

One point of view is that a 2-category is a generalized category (add
2-cells).  Another important one is that an \emph{object of} a
2-category is a generalized category (since \Cat\ is the primordial
2-category).  This is ``formal category theory'': think of a 2-category
as a collection of category-like things.

You don't capture all of \V-category theory by thinking of \V-categories
as objects of \VCat, just as you don't capture all of group theory by
thinking of groups as objects of \Grp, but many things do work
out well when we take this ``element-free'' approach.  
In formal category theory you tend to avoid talking about objects of a
category, instead talking about morphisms (functors) into the
category. Thus morphisms become generalized objects (of their 
codomain) in exactly the same way that morphisms in categories are
generalized elements.

One of the starting points of formal category theory was 
Street's beautiful work on the ``formal theory of monads''.  This
was motivated by the desire to develop a uniform approach to universal
algebra for enriched and internal categories.  It uses all four
dualities to incredible effect.

\subsection{Adjunctions and equivalences}

We start here with the notion of adjunction in a 2-category (in other
words, adjunction between objects of a 2-category --- this is not to 
be confused with adjunctions between 2-categories). In ordinary 
category theory there are two main ways to say that a functor $f:A\to B$
is left adjoint to $u:B\to A$. First there is the local approach,
consisting of a bijection between hom-sets
$$B(fa,b)\cong A(a,ub)$$
for each object $a\in A$ and $b\in B$, natural in both $a$ and $b$. 
Alternatively, there is the global approach, involving natural 
transformations $\eta:1_A\to uf$ and $\epsilon:fu\to 1_B$ satisfying
the usual triangle equations. Each can be generalized to the 2-categorical
setting. 

Let \K\ be a 2-category. Everything I'm going to say works for bicategories,
but let's keep things simple; of course you can always replace a bicategory
by a biequivalent 2-category anyway.

An {\em adjunction} in \K\ consists of 1-cells $f:A\to B$ and $u:B\to A$,
and 2-cells $\eta:1_A\to uf$ and $\epsilon:fu\to 1_B$ satisfying the 
triangle equations. This is exactly the global approach to ordinary
adjunctions, with functors replaced by 1-cells, and natural transformations
by 2-cells. In a lot of 2-categories, this is a good thing to study.
We mentioned above the case \MonCat. The study of adjunctions in \textbf{Mod}
is called {\em Morita theory}: it involves adjunctions and
equivalences between categories of the form \RMod for a ring $R$.

In the case where $\eta$ and $\epsilon$ are invertible, we have not just an 
adjunction but an \emph{adjoint equivalence}.

The local approach to adjunctions also works well here, provided that
one uses generalized objects rather than objects. For 
any 1-cells $a:X\to A$ and $b:X\to B$, there is a bijection between
2-cells $fa\to b$ and 2-cells $a\to ub$. One now has naturality with
respect to both 1-cells $x:Y\to X$, and 2-cells $a\to a'$ or $b\to b'$.
This local-global correspondence can be proved more or less as in the 
usual case, or it can be deduced from the usual case using a suitable 
version of the Yoneda lemma. In fact the global-to-local part follows
from the easy fact that \emph{2-functors preserve adjunctions}, so that the 
representable 2-functors $\K(X,-)$ send the adjunction $f\dashv u$ in \K
to an adjunction $\K(X,f)\dashv\K(X,u)$ in \Cat, between $\K(X,A)$ and
$\K(X,B)$, and so the usual properties of adjunctions give the correspondence
between $fa=\K(X,f)a\to b$ and $a\to\K(X,u)b=ub$.

The contravariant representable functors
\[\K(-,X)\maps \K\op\to\Cat\]
also preserve adjunctions. This prepares you for:

\begin{xca}
  $f$ is a left adjoint in \K\ if and only if it is a right adjoint in $\K\co$
  if and only if it is a right adjoint in $\K\op$.
\end{xca}

\begin{xca}
  A morphism $f:A\to B$ in a 2-category \K is said to be an
  equivalence if there exist a morphism $g:B\to A$ and isomorphisms
  $gf\cong 1_A$ and $fg\cong 1_B$. Show that for any equivalence $f$
  these data can be chosen so as to give an adjoint equivalence. Hint:
  you can keep the same $f$ and $g$; you'll need to change at most
  one of the isomorphisms.
\end{xca}

Considering an adjunction $f\dashv u$ in \K as an adjunction in
$\K\op$, and using the local approach, we see that to give a 2-cell 
$s\to tf$ is the same as to 
give a 2-cell $su\to t$.  Even in the case $\K=\Cat$ this is not as 
well known as it should be.

More generally, given a pair of adjunctions $f\dashv u$ and 
$f'\dashv u'$, we have bijections between 2-cells $f'a\to bf$, 
2-cells $a\to u'bf$, and 2-cells $af'\to u'b$:  squares
\[\xymatrix @C1pc {
A \ar[rr]^f\ar[d]_a  & {}\dtwocell\omit{^\beta} & B \ar[d]^b\\
  A'\ar[rr]_{f'} & {} & B'}\]
correspond to squares
\[\xymatrix @C1pc {A \ar[d]_a & {}\dtwocell\omit{^\alpha} & B \ar[ll]_u\ar[d]^b\\
  A' &{}& B'\ar[ll]^{u'}}\]
These pairs of 2-cells are called \emph{mates}. To pass between
$\alpha$ and $\beta$ one pastes with the unit and counit:
$$\xymatrix @C1pc {
& {}\ddtwocell\omit{^\epsilon} & B \ar[dll]_{u} \ar[dd]^{1} & 
&&&
   A \ar[dd]_{1} \ar[drr]^{f}  & {}\ddtwocell\omit{^\eta} \\
A \ar[drr]^{f} \ar[d]_{a} & {}\ddtwocell\omit{^\beta} &&
&&&
   & {}\ddtwocell\omit{^\alpha} & B \ar[dll]_{u} \ar[d]^{b} \\
A' \ar[dd]_{1} \ar[drr]^{f'} & {}\ddtwocell\omit{^\eta} & B' \ar[d]^{b} &
&&&
   A \ar[d]_{a} & {}\ddtwocell\omit{^\epsilon} & B' \ar[dd]^{1}\ar[dll]_{u'} \\
& {}& B' \ar[dll]^{u'} &
&&&
   A' \ar[drr]_{f'} & {}\\
A' & {} && 
&&&
   & {} & B' }$$

\subsection{Extensions}
\label{sec:extensions}

Extensions generalize Kan extensions. They provide limit 
and colimit notions for {\em objects of a 2-category}, generalizing
the usual notions for categories.

Let \K\ be a 2-category.  What is the universal solution to extending
$f$ along $j$?
\[\xymatrix{ B \ar@{.>}@/^/[dr]\\
  A\ar[u]^j \ar[r]_f \rtwocell\omit{^<-2>} & C}\]
Such a universal solution is denoted $\lan_j f$; by universal we 
mean that it induces a bijection
\[\frac{f\too gj}{\lan_j f \too g}\]
for any $g:B\to C$.
When such a $\lan_j f$ exists in \K, it is  called a \emph{left
extension} of $f$ along $j$.

A colimit is called {\em absolute} if it is preserved by any functor; similarly
we say that the left extension $\lan_j f$ is absolute if composing
with any $h:C\to D$ gives another extension, so that $h\lan_j f=\lan_j(hf)$.

Consider the case $\K=\Cat$. There would be such a bijection if
$\lan_j f$ were the left Kan extension $\Lan_j f$ of $f$ along $j$, 
as indeed the notation is supposed to suggest.  In the case of
(pointwise) left Kan extensions, we have a coend formula
\[(\lan_j f)b = \int^a B(ja,b) \cdot fa.\]
Alternatively the right hand side can be expressed using colimits: 
given $b$ we can form the comma category $j/b$, with pairs $(a\in A,ja\to b)$
as objects, and the canonical functor $d:j/b\to A$, then the coend
on the right hand side is (canonically isomorphic to) the colimit 
of $fd:j/b\to C$. 

Kan extensions which are not `pointwise' --- in other words, which
don't satisfy this formula --- can exist if $C$ is not cocomplete, 
but should be regarded as somewhat pathological.

How might we express this formula so that it makes sense in an arbitrary
2-category? Once again, the answer will involve generalized objects.
Staying for a moment in the case of \Cat, consider an object $b\in B$
as a morphism $b\maps 1\to B$, and then consider the diagram
\[\xymatrix{
1 \ar[r]^b & B \ar@/^/[dr]^{\lan_j f}  \\
j/b \ar[r]_d\ar[u]^c \rtwocell\omit{^<-3>} & 
A \ar[u]^j \ar[r]_f \rtwocell\omit{^<-3>} & C }\]
in which $j/b$ is the comma category. The coend $\int^a B(ja,b)\cdot fa$  
is isomorphic to the colimit of $fd$, as we saw, but the colimit of $fd$
is itself isomorphic to the left Kan extension of $fd$ along the 
unique map $j/b\to 1$. A careful calculation of the isomorphisms involved
reveals that the coend formula amounts to the assertion that the 
diagram above is a left extension.

This motivates the definition of pointwise extension in a
general 2-category \K with comma objects. We say that the left 
extension $\lan_j f$ is {\em pointwise} if, for any $b:X\to B$, when 
we form the comma object the 2-cell
\[\xymatrix{
X \ar[r]^b & B \ar@/^/[dr]^{\lan_j f}  \\
j/b \ar[r]_d\ar[u]^c \rtwocell\omit{^<-3>} & 
A \ar[u]^j \ar[r]_f \rtwocell\omit{^<-3>} & C }\]
exhibits $(\lan_j f)b$ as $\lan_c(fd)$.

This agrees with the usual definition in the case $\K=\Cat$,
works perfectly in the case of \CatE, and captures many but 
not all features in \VCat. The problem is that for \V-categories $A$
and $B$, the (\V-)functor category $[A,B]$ should really be regarded
as a \V-category, but the 2-category \VCat can't see this extra
structure. There are ways around this if $B$ is sufficiently complete
or cocomplete.

Let's leave the pointwise aspect aside and go back to extensions.
\begin{itemize}
\item A left extension in $\K\co$ (reverse the 2-cells) is called a
  \emph{right extension}.
\item A left extension in $\K\op$ (reverse the 1-cells) is called a
  \emph{left lifting}.
\item A left extension in $\K\coop$ (reverse both) is called a
  \emph{right lifting}.
\end{itemize}

The right lifting $r:X\to A$ of $b:X\to B$ through $f:A\to B$
is characterized by a bijection
\[\frac{fa \too b}{a\too r}\]
which is a sort of internal-hom; indeed, in the one-object case,
where the composite $fa$ is given by tensoring, it really is an
internal hom. 
Some people use the notation $r=f\lhom b$ for this lifting.

A special case is adjunctions. Given $f\dashv u:B\to A$, we have a bijection
\[\frac{fa\too b}{a\too ub}\]
and so $ub=f\lhom b$
is the right lifting of $b$ through $f$. In particular, $u$
is the right lifting of the identity $1_B$ through $f$. Conversely,
a right lifting $u$ of the identity through $f$ is a right adjoint
if and only if it is absolute; in other words, if $ub$ is the right
lifting of $b$ through $f$ for all $b:X\to B$; in symbols 
$f\lhom b=(f\lhom 1)b$.
%In particular, $u = f\rhd 1$.  Thus every adjunction gives a right
%lifting, and a right lifting $u = f\rhd 1$ is an adjunction iff it is
%``respected'' by any $b$, in the sense that $ub = f\rhd b$.  

Dually, given an adjunction $f\dashv u:B\to A$ we have a bijection 
\[\frac{xu\too y}{x\too yf}\]
and so $yf=\ran_u y$ and $f=\ran_u 1_B$; while in general a right
extension $f=\ran_u 1_B$ of the identity is a left adjoint of $u$
if and only if it is absolute.

A bicategory is said to be \emph{closed} if it has
right extensions and right liftings. In the one-object case, this
means that the endofunctors $-\ot c$ and $c\ot -$ of the monoidal
category have right adjoints for any object $c$.

We saw that pointwise left extensions in \Cat are given by colimits.
Thus the existence of left extensions is some kind of internal 
cocompleteness condition. So in 2-categories like \CatE\ or \VCat\
they will exist only in some cases. In bicategories like \textbf{\V-Mod},
on the other hand, all extensions exist (provided that \V is itself
complete and cocomplete).

Let me point out a little lemma which everyone knows for \Cat, but
which is true for 2-categories basically because everything is
representable. A morphism $f:A\to B$ in a 2-category \K is said to 
be {\em representably fully faithful} if $\K(X,f):\K(X,A)\to\K(X,B)$ 
is a fully faithful functor for all objects $X$ of \K. For $\K=\Cat$
this is equivalent to $f$ being fully faithful.

\begin{lemma}
Let $f\dashv u$ be an adjunction in a 2-category \K for which 
the unit $\eta:1\to uf$ is invertible. Then $f$ is representably
fully faithful.
\end{lemma}

Similarly, under the same hypotheses, $u$ will be (representably)
``co-fully-faithful'', in the sense that each
$\K(u,X):\K(B,X)\to\K(A,X)$
is fully faithful.

\subsection{Monads}
\label{sec:monads}

Just as in ordinary category theory, an adjunction $f\dashv u:B\to A$ 
in a 2-category induces a 1-cell $t=uf$, with 2-cells $\eta:1\to
uf=t$, given by the unit of the adjunction, and a multiplication 
$\mu=u\epsilon f:t^2=ufuf\to uf=t$, where
$\epsilon:fu\to 1$ is the counit. This $\eta$ and $\mu$ make $t$
into a monoid in the monoidal category $\K(A,A)$. 

More generally, a monad in a 2-category \K on an object $A\in\K$ consists
of a 1-cell $t:A\to A$ equipped with 2-cells $\eta:1\to t$ and $\mu:t^2\to t$
satisfying the usual (associative and identity) equations; the situation
of the previous paragraph is a special case. One often speaks simply of
a monad $(A,t)$, when $\eta$ and $\mu$ are understood.

The case $\K=\Cat$ is just the usual notion of monad on a category $A$.
(This is sometimes called a monad {\em in} $A$, but this usage is to be 
avoided: it is {\em in} \K and {\em on} $A$.)

\begin{example}
  Monads in \Cat are the usual monads. Monads in \VCat\ or
  \CatE\ are the obvious notion of enriched or internal monad.
  Monads in \MonCat\ are called monoidal monads. Monads in the 
  2-category \OpMonCat\ of monoidal categories, opmonoidal functors,
and opmonoidal natural transformations are called opmonoidal monads,
or sometimes Hopf monads (see \cite{Moerdijk-Hopf-monads}).
\end{example}

\begin{example}
  Monads in the one-object 2-category $\Sigma\V$ are monoids in the 
  strict monoidal category \V. Conversely, a monad in an arbitrary 
  2-category \K, on an object $X$ of \K, is a monoid in the (strict)
  monoidal category $\K(X,X)$. There are analogous facts for
  bicategories and (not necessarily strict) monoidal categories.
\end{example}

\begin{example}
  Monads in \Rel. We have a set $E_0$; a relation $t\maps E_0\hto E_0$,
  in the form of a subset $R$ of $E_0\t E_0$; the ``identity'' $1\to t$
  amounts to the assertion that the relation $R$ is reflexive, and the
  multiplication to the fact that $R$ is transitive. The associative
  and unit laws are automatic.
\end{example}

\begin{example}\label{ex:monads-SpanE}
  Monads in \SpanE.  We have an object $E_0$, a 1-cell
  $t\maps E_0\hto E_0$, as in
  \[\xymatrix @R1pc { & E_1 \ar[dl]_d \ar[dr]^c\\ E_0 && E_0}\]
  (a directed graph in \EE), with a multiplication
  \[\mu\maps E_1\times_{E_0} E_1 \to E_1\]
  from the object of composable pairs to the object of morphisms, giving a
  composite; associativity of the monad multiplication is precisely
  associativity of the composition. Similarly the unit
  $1\too[\eta] t$ gives $E_0\to E_1$ since the identity span is
  \[\xymatrix @R1pc {& E_0\ar[dl] \ar[dr] \\ E_0 && E_0 ~,}\]
  and the unit laws for the monad are precisely the identity laws for
  the internal category.
  Thus \emph{a monad in {\em\SpanE} is the same as
    an internal category in \EE}.

  This is one of the main reasons for considering
  the span construction.
\end{example}

\begin{example}
  Monads in \VMat.  We have an object $X$, which is just a set, a
  1-cell $X\hto X$, in the form of a matrix $X\times X \to \V$, 
  which we think of as sending $(x,y)$ to a hom-object  $\C(x,y)$.
  The multiplication map goes from the matrix product, as in
  \[\sum_y \C(y,z) \ot \C(x,y) \too \C(x,z)\]
  and gives a composition map.  Once again the associative and identity
  laws for the composition are precisely the associative and unit laws
  for the monad, and we see that \emph{a monad in \VMat\ is the same as
    a category enriched in \V}. 

  In the special case $\V=\two$ we have $\VMat=\Rel$, and so we recover
  the observation, made in Example~(\ref{item:VMat}) above, that
  a category enriched in \two\ is just a preorder (a reflexive and transitive
  relation).
\end{example}

A morphism of monads from $(A,t)$ to $(B,s)$ consists of a
1-cell $f:A\to B$ equipped with a 2-cell $\phi:sf\to ft$, satisfying
two conditions: see \cite{ftm} or Section~\ref{sect:ftm} below.
A morphism of monads in \SpanE\ is  \emph{not} an internal functor,
since it would involve a 1-cell (a span) $E_0\hto F_0$ between the 
objects of objects, rather than a morphism in \EE. In order to 
get internal functors, we need to consider not \SpanE\ itself, but rather
\SpanE\ equipped with the class of ``special'' 1-cells consisting of those
spans whose left leg is the identity; these can of course be identified
with the 1-cells in \EE.
An internal functor will turn out to be a monad morphism, for which the span
$E_0\hto F_0$ is ``special''.

The case of enriched functors is similar: one needs to keep track of
which 1-cells in \VMat\ are really just functions.

To get (enriched or internal) natural transformations, you do not use
the obvious notion of monad 2-cells as in \cite{ftm}, but rather those 
of \cite{ftm2}; once again see Section~\ref{sect:ftm} below.

\subsection{References to the literature}

For adjunctions in 2-categories and the calculus of mates see
\cite{Kelly-Street} or \cite{Gray}. For monads in 2-categories see the 
classic \cite{ftm}. 
For extensions and liftings see \cite{yoneda}. The idea that
categories can be seen as monads in \Span comes from \cite{bicategories}.

\section{Morphisms between bicategories}
\label{sect:morphisms}

\subsection{Lax morphisms}
\label{sect:lax-functors}

We talked before about the virtues of {\em monoidal functors} between
monoidal categories. The corresponding morphisms between bicategories
are the {\em lax functors} (originally just called morphisms of bicategories
by B\'enabou). A lax functor $\A\to\B$ sends 
objects $A\in\A$ to objects $FA\in\B$, has functors $F\maps \A(A,B)\to
\B(FA,FB)$ (thus preserving 2-cell composition in a strict way),
and has comparison maps $\phi\maps Fg\cdot Ff \to F(gf)$ and
$\phi_0\maps 1_{FA} \to F(1_A)$ and some coherence conditions, which
are formally identical to those for monoidal functors.

All the good things that happen for monoidal functors happen for
lax functors. For example, monoidal functors take monoids to monoids,
and lax functors take monads to monads. (Recall that a monad in \B\
on an object $X$ is the same as a monoid in the monoidal category
$\B(X,X)$.) 

As a very special case, consider the terminal 2-category $1$. 
This has a unique object $*$, and a unique monad on $*$ (the identity
monad). Then for any lax functor $1\to \B$, the object $*$ gets sent
to $F*= A$, the identity $1$ is sent to $F1=t$, the comparison maps
become $\mu\maps tt\to t$ and $\eta\maps 1\to t$, and the coherence 
conditions make this precisely a monad.  In fact, 
\emph{monads in \B are the same as lax functors $1\to \B$}.  
For B\'enabou, this was a key reason to consider lax morphisms of 
bicategories, rather than the stronger version.  

In particular, \V-categories are the same as monads in \VMat, and
so also the same as lax functors $1\to\VMat$.  This is the same as a 
set $X$ together with a lax functor
\[X\ch \too \Sigma\V\]
where $X\ch$ is $X$ made into a chaotic bicategory (also called
indiscrete: every hom-category $X\ch(x,y)$ is trivial).
Why?  We send each $x$ to $*$, we have a
functor 
\[1=X\ch(x,y) \to \Sigma\V(*,*)=\V\]
picking out the hom-object $\C(x,y)\in\V$, and the lax
comparison maps $\phi$ become the composition and identity maps.

If we replace $\Sigma\V$ by an arbitrary bicategory \W, we get the
notion of a \W-enriched category: a set $X$ with a lax functor
\[X\ch \too \W\]

Another way to think about $X\ch$, as a bicategory, is to say that
the unique map $X\to 1$ is fully faithful.  But we can also consider,
more generally, a pair of bicategories with a partial map
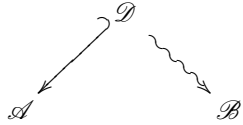
\[\xymatrix{& \D\ar@{_(->}[dl]  \ar@{~>}[dr] \\ \A && \B}\]
where the hooked arrow $\hookrightarrow$ denotes a fully faithful strict
morphism, and the wobbly map $\rightsquigarrow$ denotes a lax functor.
This partial map is called a \emph{2-sided enrichment} or a
\emph{category enriched from \A\ to \B}.  If \A\ is 1, it's
just a category enriched over \B.  

Using the notion of composition
for these things is very helpful in analyzing the change of base
between different bicategories. For example, a
\B-category is a partial map from 1 to \B; this can be composed
with a partial map from \B to \C to get a \C-category. As a
(better-known) special case, lax functors from \B to \C send 
\B-categories to \C-categories; as a still more special case, 
monoidal functors from \V to \W send \V-categories to \W-categories.

\subsection{Pseudofunctors and 2-functors}
\label{sec:pseudofunctors}

A pseudofunctor (or homomorphism of bicategories) is a lax functor for 
which $\phi$ and $\phi_0$ are invertible.

\begin{example}
  For a bicategory \B, the representables
  \[\xymatrix @C3pc {\B \ar[r]^{\B(B,-)} &  \Cat}\]
  are pseudofunctors, not strict in general.
\end{example}

\begin{example}[Indexed categories]
  A pseudofunctor $\B\op\to\Cat$ is sometimes called a \B-indexed category.
  Often \B itself will just be a category (no non-identity 2-cells), in  
  which case such a pseudofunctor corresponds to a fibration $\E\to\B$ in
  the Grothendieck picture.
\end{example}

An important property of pseudofunctors not shared by lax functors
is that they preserve adjunctions. Consider a pseudofunctor $F:\A\to\B$,
and an adjunction $f\dashv u:B\to A$ in \A, with unit $\eta:1_A\to uf$
and $\epsilon:fu\to 1_B$. We may apply $F$ to $f$ and $u$ to get 
$Ff:FA\to FB$ and $Fu:FB\to FA$, and now the composite 2-cells
$$\xymatrix{
Ff.Fu \ar[r]^{\phi} & F(fu) \ar[r]^{F\epsilon} & F1_B
\ar[r]^{\phi\inv_0} & 1_{FB} \\
1_{FA} \ar[r]^{\phi_0} & F1_A \ar[r]^{F\eta} & F(uf) \ar[r]^{\phi\inv} &
Fu.Ff }$$
provide the unit and counit for an adjunction $Ff\dashv Fu$. This
fails for a general lax functor $F$.

If $\phi$ and $\phi_0$ are not just invertible, but in fact identities,
then one speaks of a strict homomorphism; or, in the case of 2-categories,
of a 2-functor. Note that in the bicategory case the associativity and identity
constraints must still be preserved: this is the content of the
coherence condition for $\phi$ and $\phi_0$.

2-functors are much nicer to work with, but often it is the pseudofunctors
which arise in nature. One reason you might prefer 2-functors is so as not 
to have to worry about coherence.  Furthermore, 2-functors have better 
properties than pseudofunctors:  for example, the category \twocat of 
2-categories and 2-functors has limits and colimits, but the category 
\twocatps of 2-categories and pseudofunctors does not. 
For example the diagram
$$
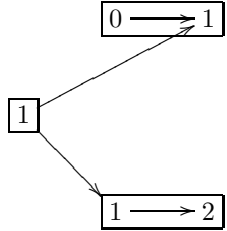
\xymatrix @R2pc { 
& 0 \save [].[r]*[F]\frm{}\restore \ar[r] & 1 \\
**[F]{1} \ar[urr] \ar[dr] \\
&  1 \save [].[r]*[F]\frm{}\restore \ar[r] & 2 
}$$
has no pushout: such a pushout would have to have morphisms
$0\to 1\to 2$ 
and a composite, but in some other cocone we have no way to decide where
to send the composite.  
If, however, we made \twocatps into a tricategory, then it would
have trilimits (the relevant ``weak'' notion of limits for tricategories).

On the other hand, even if you start in the world of 2-categories and
2-functors, you may be forced out of it.  A 2-functor
$\A\too[F]\B$ is a
\emph{biequivalence} if $\A(A,B)\to \B(FA,FB)$ are equivalences
and it is ``bi-essentially surjective'', in the sense that for all $X\in \B$,
there exists an $A\in\A$ and an equivalence $FA\simeq X$ in \B.  This is
the ``right notion'' of equivalence for 2-functors.

The point is that you'd like something going back the other way from
\B to \A. Well you do have {\em something}, but it's just not a 2-functor 
in general.  Given
$X\in\B$, pick $A\in\A$ and $FA\simeq X$ and let $GX=A$.  Given $X\too[x] Y$,
we can bring it across the equivalences $FA\simeq X$ and $FB\simeq Y$ to
get $\bar{x}\maps FA\to FB$, and since $F$ is locally an equivalence,
$\bar{x} \cong Fa$ for some $a\maps A\to B$; let $Gx = a$.  This all works,
but since everything is only defined up to isomorphism, there's no way
you can possibly hope for $G$ to preserve things strictly.

There is a Quillen model structure on \textbf{2-Cat} --- see
Section~\ref{sect:qmc2cat} below --- for which the weak
equivalences are the biequivalences, and clearly getting a 2-functor
$\B\to\A$ is going to have something to do with \B\ being cofibrant.

\subsection{Higher structure}

As well as lax (and other) morphisms between bicategories, there is
higher structure. Given morphisms $F,G:\A\to\B$, one can consider
families $\alpha A:FA\to GA$ of morphisms in \B indexed by the objects 
of \A, and subject to (lax, oplax, pseudo, or strict) naturality 
conditions. There is even a further level of structure, consisting 
of morphisms between such transformations: these are called 
modifications.

\subsection{References to the literature}
The importance of lax functors,
especially lax functors with domain 1, was observed by B\'enabou 
\cite{bicategories}.

Categories enriched in a bicategory were first defined by Walters
to deal with the example of sheaves (on a space or site) 
\cite{Walters-sheaves1,Walters-sheaves2}.
A good general reference is \cite{BCSW}. Two-sided enrichments (although
not from the point of view of partial morphisms of bicategories) were
defined in \cite{KLSS}, in order to deal with change of base issues.

\section{2-dimensional universal algebra}
\label{sect:2DUA}

There are various categorical approaches to universal algebra: theories,
operads, sketches, and others, but I'll mostly talk about monads, although
you may see parallels with operads and with theories if you know about those.

The ordinary universal algebra picture you might have in mind is
monoids (or groups, rings, etc.) living over sets. But our algebras don't
have to be single-sorted; they could live over some power of sets.  Abstractly,
of course, we could be living over almost everything.  A good many-sorted
example to have in mind is the functor category $[\C,\Set]$ living over
$[\ob\C,\Set]$, for a small category \C. If \C\ has one object, then 
we may identify \C\ with the monoid $M$ of its arrows, and the functor
category is then the category of $M$-sets.

When we come to 2-categories, we might generalize monoids over sets to
monoidal categories over categories; or (also living over categories) 
categories with finite products, or with finite coproducts, or with both, or
with finite products and finite coproducts and a distributive law.

For an example of the many-sorted case, let \B be a small bicategory. 
There is a 2-category $\Hom(\B,\Cat)$ of homomorphisms from \B to \Cat
(\B-indexed categories), whose morphisms and 2-cells are the pseudonatural
transformations and modifications. The domain of the forgetful 2-functor
\[\xymatrix{ \Hom(\B,\Cat) \ar[d]\\ [\ob\B,\Cat]}\]
is an example of the sort of algebraic structure we have in mind.

In the next two sections there is a lot of interplay between 2-category
theory and \Cat-category theory. Since I don't want to assume enriched
category theory, I'll tend to describe the ordinary (unenriched) setting,
take it for granted that one can modify this to get a \Cat-enriched version,
and concentrate more on how to modify this to do the proper 2-categorical
one.  

\subsection{2-monads}
\label{sec:2-monads}

We continue to follow the convention that the prefix 2- indicates a strict 
notion. Thus a 2-monad consists of a 2-category \K\ equipped with a
2-functor $T:\K\to\K$, and 2-natural transformations $m:T^2\to T$ and
$i:1\to T$, satisfying the usual equations for a monad. In other words,
this is a monad in the (large) 2-category of 2-categories, 2-functors,
and 2-natural transformations. (This could be made into a 3-category,
but we don't need to do so for this observation.) 

There is a good theory of enriched monads --- this was one of the 
motivations of the formal theory of monads --- and 2-monads are 
just \V-monads in the case $\V=\Cat$.

A (strict) \emph{$T$-algebra} is the usual thing, an object $A\in \K$ with
a morphism $a:TA\to A$ satisfying the usual equations, written $(A,a)$. 
Once again, this is the strict (or \Cat-enriched) notion. 

\begin{remark}\label{rmk:pseudo-algebras}
There are pseudo and lax notions of monad and of algebra, but they seem
to be less important in universal algebra than the strict ones. The
main reason for this 
is that the actual structures one wants to describe using 2-dimensional
monads are the strict algebras for strict monads in a fairly straightforward
way --- an example is given below --- whereas identifying the structures of 
interest with pseudoalgebras is rather more work. A secondary reason is
that in reasonable cases a pseudomonad $T$ can be replaced by a strict monad 
$T'$ whose strict algebras are the pseudoalgebras of $T$.
\end{remark}

It is when we come to the {\em morphisms} of algebras that
we are forced to depart from the strict setting.  A \emph{lax
  $T$-morphism} $(A,a)\to (B,b)$ is a morphism $f\maps A\to B$ in \K,
equipped with a 2-cell
\[\xymatrix{TA \ar[r]^{Tf} \ar[d]_a  & TB \ar[d]^b \\
  A\ar[r]_{f} \rtwocell\omit{<-4>\bar{f}} & B}\]
satisfying two coherence conditions:
\[\xymatrix{ 
T^2A \ar[r]^{T^2f}\ar[d]_{Ta}  & T^2B
  \ar[d]^{Tb} && T^2A \ar[r]^{T^2f}\ar[d]_{mA} & T^2B \ar[d]^{mB} \\
  TA \ar[r]_{Tf}\ar[d]_a \rtwocell\omit{<-4>~~T\fbar} & TB \ar[d]^b &=&
  TA \ar[r]_{Tf}\ar[d]_a  \rtwocell\omit{<-4>\fbar} & TB \ar[d]^b  \\
  A \ar[r]_f & B &&   A \ar[r]_f \rtwocell\omit{<-4>\fbar} & B}
\]
\[\xymatrix{A \ar[r]^{f} \ar[d]_i & B \ar[d]^{iB} && 
A \ar[r]^{f} \ar[dd]_{1} & B \ar[dd]^{1} \\
  TA \ar[r]^{Tf} \ar[d]_{a} & TB \ar[d]^{b} &=& \\
  A\ar[r]_{f} \rtwocell\omit{<-4>\bar{f}}& B && A \ar[r]_{f} & B.}\]  
Note that the outer 1-cells are the same (I wouldn't write this equation down
if they weren't), and that empty regions commute, and are deemed to
contain the relevant identity 2-cell.

Let's do a baby example: $\K=\Cat$ and $TA=\sum_n A^n$ the
usual free monoid construction.  The $T$-algebras are strict monoidal
categories, and a lax morphism involves 2-cells 
\[\xymatrix{{\sum_n A^n} \ar[r]\ar[d]_{\ot} 
& {\sum_n B^n} \ar[d]^{\ot} \\
  A \ar[r] \rtwocell\omit{<-4>\bar{f}} & B}\]
so we have transformations
\[f(a_1)\ot\dots\ot f(a_n) \too f(a_1\ot\dots\ot a_n).\]
for each $n$. The definition of
monoidal functor only mentions the cases $n=0$ and $n=2$,
but all the others can be built up from these in an obvious way; the 
coherence conditions for lax $T$-morphisms say that you {\em do} build 
them up in this sensible way, and that the coherence conditions for
monoidal functors are satisfied. 

So for this $T$, the lax morphisms are precisely the monoidal functors.
This provides a practical motivation for the definition of lax $T$-morphism. 
Here's a theoretical one. There's a 2-category $\mathbf{Lax}(\two,\K)$ where
$\two$ is the arrow category.  In detail:
\begin{itemize}
\item An object is an arrow $a:A'\to A$ in \K
\item A 1-cell is a square
  \[\xymatrix{A' \ar[r]^{f'} \ar[d]_{a} & B' \ar[d]^{b}\\
    A\ar[r]_{f} \rtwocell\omit{<-4>\phi} & B}\]
\item A 2-cell $(f,\phi,f')\to(g,\psi,g')$ consists of 2-cells 
$\alpha:f\to g$ and $\alpha':f'\to g'$ satisfying the equation
$$\xymatrix{
A' \ar[d]_{a} \rtwocell^{f'}_{g'}{\alpha'} & B' \ar[d]^{b} \\
A \rlowertwocell_g<-2>{<-2>\psi} & B } =
\xymatrix{A' \ar[d]_{a} \ruppertwocell^{f'}<3>{<2>\phi} & B' \ar[d]^b \\
A \rtwocell^{f}_{g}{\alpha} & B. }$$
\end{itemize}
Since this is functorial in \K, the 2-monad $T$ induces a 2-monad
$\mathbf{Lax}(\two,T)$ on $\mathbf{Lax}(\two,\K)$.  Then a (strict)
$\mathbf{Lax}(\two,T)$-algebra is precisely a lax $T$-morphism. 
The coherence conditions for lax morphisms become the usual 
axioms for algebras.

Similarly, a \emph{$T$-transformation} between lax $T$-morphisms
$(f,\fbar),(g,\gbar)\maps (A,a)\to (B,b)$ is a 2-cell $\rho\maps f\to g$ in
\K\ such that
\[\xymatrix{TA \ar[d]_a \rtwocell{T\rho}  & TB \ar[d]^b \\
  A \rlowertwocell_g<-2>{<-3>\gbar} & B} = 
\xymatrix{TA \ruppertwocell^f<3>{<3>\fbar} \ar[d]_a & TB \ar[d]^b \\
  A \rtwocell{\rho} & B}\]

In the baby example, for $n=2$ this says that
\[\xymatrix{fa_1\ot fa_2 \ar[r]\ar[d]_{\rho a_1\ot \rho a_2} & 
f(a_1\ot a_2) \ar[d]^\rho\\
  ga_1 \ot ga_2\ar[r] & g(a_1\ot a_2)}\]
which is exactly the condition for $\rho:f\to g$ to be a monoidal natural
transformation.

\begin{xca}
Play the $\mathbf{Lax}(\two,\K)$ game with $T$-transformations: find
a 2-category $\K'$ and a 2-monad $T'$ on $\K'$ whose algebras are
the $T$-transformations.
\end{xca}

There is a 2-category \talgl of $T$-algebras, lax
$T$-morphisms, and $T$-transformations, and a forgetful 2-functor
\[\talgl \too[U_\ell] \K \]
and in some cases, such as that of monoidal categories, this is the
2-category of primary interest, but often the pseudo case is more
important (and of course strong monoidal functors are themselves important).
If $\fbar$ is invertible, we say that
$(f,\fbar)$ is a \emph{pseudo $T$-morphism} or just a
\emph{$T$-morphism} (privileging these over the strict or the lax).
These are the morphisms of the 2-category \talg of $T$-algebras,
pseudo $T$-morphisms, and $T$-transformations; it has a forgetful
2-functor 
\[\talg \too[U] \K.\]

When $\fbar$ is an identity we have a \emph{strict} $T$-morphism. 
Of course this just means that the square commutes, and we have
a morphism in the usual unenriched sense, but it is still useful
to think of the identity 2-cell as being ``an $\fbar$'', since it
is used in the condition on 2-cells. The $T$-algebras, strict
$T$-morphisms, and $T$-transformations form a 2-category \talgs
with a 2-functor
\[\talgs \too[U_s] \K.\]
Each of these 2-categories has the same objects, and we have a
diagram
\[\xymatrix{
  && \talgc \\
  \talgs \ar[r]_J \ar@/^5mm/[rr]^{J_\ell} \ar[dr]&
  \talg \ar[r] \ar[d] \ar[ur] &
  \talgl \ar[dl]
  \\
  & \K}\]
of 2-categories and 2-functors, 
where \talgc is the 2-category of {\em colax} morphisms, defined like
lax morphisms except that the direction of the 2-cell is reversed.
I won't worry too much about them since they can be treated as the 
lax morphisms for an associated 2-monad on $\K\co$.

At this point we need to start making some assumptions. To start with,
suppose that \K is cocomplete, and that $T$ has a \emph{rank}, which means
that $T\maps\K\to\K$ preserves $\alpha$-filtered colimits for 
some $\alpha$.  For ordinary monads on categories, it says that
we can describe the structure in terms of operations which may not be
finitary, but are at least $\alpha$-ary for some regular cardinal $\alpha$.
The famous example of a
monad on \Set\ which is not $\alpha$-filtered for any $\alpha$ is the 
covariant power set monad.

Under these conditions
\[\talgs \too[J]  \talg\]
\[\talgs \too[J_\ell]  \talgl\]
have left adjoints.  What does this mean?  Among other things it means 
that for each algebra $A$ there is an algebra $A'$ and bijections
\[\frac{A\rightsquigarrow B}{A'\to B}\]
where the wobbly arrow denotes a weak morphism and the normal arrow
a strict one. Here ``weak'' might mean either pseudo or lax, depending
on the context; of course there will be a different $A'$ depending on 
whether we consider the pseudo or the lax case.

These are 2-adjunctions, so these bijections are just part
of isomorphisms of categories
\[\talgs(A',B) \cong \talg(A,JB)\]
2-natural in $A$ and $B$.  We usually omit writing the $J$,
since it is the identity on objects. We say that such an $A'$ 
{\em classifies weak morphisms out of $A$}. 
From this we get a unit
\[p\maps A\rightsquigarrow A'\]
and counit
\[q\maps A' \to A\]
and one of the triangle equations tells you that $qp = 1$.  
An unfortunate consequence of the (otherwise reasonable) notation
$A'$, is that the left adjoint to $J:\talgs\to\talg$ is sometimes
saddled with the rather embarrassing name $(~)'$; I shall call it $Q$
instead.

\subsection{Sketch proof of the existence of $A'$}
\label{sec:sketch-proof-exist}

\subsubsection*{Step 1}
\label{sec:step-1}

\textbf{\talgs is cocomplete.}

This part is entirely ``strict'': it is really an
enriched category phenomenon, and not really 
any harder than the corresponding fact for ordinary categories.
It is here that you use the assumptions on \K and $T$. 

Colimits of algebras, as we know, are
generally hard.  The problem is essentially that algebras are
a ``quadratic'' notion, involving  $a:TA\to A$ with two copies of $A$.  
We ``linearize'' and it becomes easy.  What does that mean?

Take the $T$-algebra $(A,a:TA\to A)$, forget the axioms, and also forget that
the two $A$'s are the same, so consider it only as a map $a:TA\to A_1$.
This defines the objects of a new category, whose morphisms are squares
of the form
\[\xymatrix{TA \ar[r]^{Tf} \ar[d]_a & TB \ar[d]^b \\
  A_1\ar[r]_{f_1} & B_1. }\]
With the obvious notion of 2-cell this becomes a 2-category; in fact
it is just the comma 2-category $T/\K$. The point is that we have a
full embedding
$$\talgs\hookrightarrow T/\K$$
since by the unit condition for algebras, any morphism in $T/\K$ 
between algebras must have $f=f_1$ and so be a strict $T$-morphism. 
It is this $T/\K$ which is the ``linearization''
of \talgs, and colimits in it are easy. Say we have a diagram
of things $TA_i\to B_i$.  Take the colimits in \K\ and take the
pushout
\[\xymatrix{\colim TA_i \ar[r]\ar[d] & \colim B_i \ar[d]\\
  T\colim A_i \ar[r] & B}\]
to get the colimits in $T/\K$.  

The hard bit, which I'll leave out, is the construction of a 
reflection $T/\K\to \talgs$ (a left
adjoint to the inclusion). This is where we use the assumption on $T$.
There are some transfinite calculations, as you might expect given
the condition on $\alpha$-filtered colimits.

Note, however, that should $T$ preserve all colimits, then this Step~1
becomes easy: the colimits are constructed pointwise. In particular,
this is true in the case of categories of diagrams ($\talgs=[\B,\Cat]$).

In fact when we come to step 2, we'll see that only 
certain (finite) colimits in \talgs are actually needed, and if $T$
should preserve these colimits, as does sometimes happen, then once
again the proof simplifies.

\subsubsection*{Step 2}
\label{sec:part-2}

Let $(A,a)$ be an algebra; we want to construct the pseudomorphism
classifier $A'$ using colimits in \talgs. A lax $T$-morphism
$(A,a)\to (B,b)$ consists of various data in \K, and we want to translate 
all these data into \talgs.  

A lax $T$-morphism $A\to B$ consists of
\begin{itemize}
\item A morphism   $f\maps A\to B$ in \K, which becomes
  a morphism $g\maps TA\to B$ in \talgs, where $g=b\cdot Tf$.
\item A 2-cell
  \[\xymatrix{TA \ar[r]^{Tf}\ar[d]_a & TB \ar[d]^b\\
    A\ar[r]_f \rtwocell\omit{<-4>\fbar}& B}\]
  in \K, which becomes a 2-cell
  \[\xymatrix{T^2A \ar[r]^{mA} \ar[d]_{Ta} & TA  \ar[d]^g\\
    TA \ar[r]_g \rtwocell\omit{<-4>\zeta} & B }\]
  in \talgs, since
  \begin{gather*}
    b. T(fa) = b.Tf.Ta = g.Ta\\
    b.T(b.Tf) = b.Tb.T^2f = \dots = g.mA.
  \end{gather*}
\item The condition $\fbar .iA = \id$ corresponds to saying that
  $\zeta . TiA = \id$
\item The other condition becomes
$$\xymatrix @R1pc @C1pc {
& T^2A \ar[rr]^{mA} \ar[ddr]^{Ta} && TA \ar[ddr]^{g} & &&
& T^2A \ar[rr]^{mA} && TA \ar[ddr]^{g} \\
&& {}\rtwocell\omit{\zeta} & {} \\
T^3A \ar[uur]^{mTA} \ar[ddr]_{T^2a} && TA \ar[rr]^{g} && B &=&
T^3A \ar[rr]^{TmA} \ar[uur]^{mTA} \ar[ddr]_{T^2a} && 
T^2A \ar[uur]^{mA} \ar[ddr]_{Ta} \rrtwocell\omit{\zeta} && B \\
&& {}\rtwocell\omit{\zeta} & {} \\
& T^2A \ar[rr]_{Ta} \ar[uur]_{mA} && TA \ar[uur]_{g} & &&
& T^2A \ar[rr]_{Ta} && TA \ar[uur]_{g} }$$
\end{itemize}

We have a truncated simplicial object:
\[ \xymatrix{
  T^3A\ar@<-4mm>[r]^{mTA} \ar@<4mm>[r]^{T^2a} \ar[r]^{TmA} &
  T^2A  \ar@<-4mm>[r]^{mA} \ar@<4mm>[r]^{Ta} &
  A\ar[l]_{TiA}
}\]
We now form a 2-categorical colimit, called the \emph{codescent object}, of 
this truncated simplicial object, and the result is the desired $A'$. 
Alternatively, we can break this up into bite-sized pieces. We first 
construct the {\em coinserter} of $mA$ and $Ta$: this is the universal 
$p:TA\to A_1$ equipped with a 2-cell $\rho:p.mA\to p.Ta$.
To give a map $A_1\to B$ in \talgs is equivalently to give 
a map $f:A\to B$ in \K and a 2-cell $\bar{f}:b.Tf\to fa$, without
any coherence conditions. To capture the coherence conditions, we 
have to perform a special sort of quotient, called a {\em coequifier},
which universally makes equal a parallel pair of 2-cells. We'll 
talk about 2-categorical limits and colimits later.

If we used the ``pseudo'' version of weak morphisms, then we'd use
a {\em co-isoinserter} instead of an coinserter, which is the obvious
analogue in which $\rho$ is invertible.

\subsection{Consequences of the pseudomorphism classifier}
\label{sec:more}

Recall that we have
\[\xymatrix{&A' \ar[dr]^q \\  A\ar@{~>}[ur]^p \ar[rr]_1 && A}\]
with $qp = 1$.  It's also true that $pq \cong 1$, so that this is an
equivalence, and thus $A\simeq A'$ in \talg, although generally not in 
\talgs. If however $q$ has a section $s$ in \talgs, so that $qs=1$, then
$s\cong p$, so $sq\cong 1$, and $q$ is an equivalence in \talgs. When $q$ 
does have such a section, the algebra $A$ is said to be {\em flexible}.

We'll see in Section~\ref{sect:model-talgs} that there is a 
model structure on \talgs for which $A'$ is
a cofibrant replacement of $A$. The weak equivalences are the strict
morphisms which become equivalences in \talg, or equivalently in \K;
the cofibrant objects are precisely the flexible algebras. 

\begin{xca}
  If $A$ is flexible, then any pseudo $A\rightsquigarrow B$ is isomorphic to a
  strict $A\to B$.
\end{xca}

The equivalence $A\simeq A'$ is a kind of coherence result for morphisms.
There are also coherence results for algebras. Consider the composite
\[\talgs \to \talg \to \pstalg\]
To give a left adjoint is to construct a pseudo morphism classifier 
$A'\in\talgs$ not just for each strict $T$-algebra, but also for 
pseudo-$T$-algebras. This can still be done; rather than a truncated
simplicial object one has a truncated pseudosimplicial object (some of
the simplicial identities are satisfied only up to isomorphism), but 
we can still form the codescent object $A'$ and obtain an isomorphism
of categories
$$\talgs(A',B)\cong\pstalg(A,B)$$
for any strict algebra $B$, natural in $B$ with respect to strict maps.
This time we have a counit $q:B'\to B$ only when $B$ is strict, and a
unit $A\rightsquigarrow A'$ for any pseudo algebra $A$. For a general
pseudo algebra $A$, there seems no way to construct a map from $A'$
back to $A$, and so no way to show that $p$ is an equivalence. 
In some cases, however $p$ is an equivalence. In particular it is so
if $T$ preserves the relevant codescent objects, since then one can
construct the codescent object in \K, and get the inverse-equivalence
down there. There are various other sufficient conditions for this to work.

The existence of $A'$ for each pseudoalgebra $A$, along with the fact that
the unit $A\rightsquigarrow A'$ 
is an equivalence is sometimes called the ``full coherence result''.  

There are 2-monads for which not every pseudoalgebra is equivalent
to a strict one, but the only examples I know involve horrible 
2-categories \K. I don't know of an example satisfying the assumptions
made in this section (\K cocomplete and $T$ preserving
$\alpha$-filtered colimits).

\subsection{References to the literature}
2-monads were first considered in \cite{Lawvere:ordinal-sum}.
The basic reference is \cite{BKP}, although many of the 
key ideas go all the way  back to Kelly's \cite{Kelly-LaxAlg}, including
the constructions $\<A,A\>$ and $\{f,f\}$ which allow one to describe
algebras in terms of monad morphisms. For the latter, see also
\cite{property}. The accessibility issues in \cite{BKP}
were treated in Blackwell's (unpublished) thesis, and later in the
monumental (and somewhat impenetrable) \cite{Kelly-transfinite}, which
builds on many earlier papers, in particular \cite{Barr:free-triples}.
Anyway, \cite{BKP} contains the results about limits and (bi)colimits 
in $\textrm{T-Alg}$, biadjoints to algebraic 2-functors, and the 
left adjoint to $\textrm{T-Alg}_s\to\textrm{T-Alg}$. The proof given here
for the existence of this left adjoint follows \cite{codescent}. See 
\cite{codescent} and the references therein for a discussion of proofs
that pseudo algebras are equivalent to strict ones. (But I can't omit
explicit mention of one those references: the short and beautiful paper
\cite{GeneralCoherence} of Power.) Kelly's work \cite{Kelly:clubs} on clubs
was not at first appreciated, but it has been influentially recently in
monad-theoretic approaches to higher categories.

\section{Presentations for 2-monads}
\label{sect:presentations}

Presentations involve {\em free} gadgets and {\em colimits}.  
Both are defined in terms of a universal property involving maps
{\em out of} the constructed gadget.
Why are these important in the case of 2-monads (or monads)?
It turns out that one can use colimits to build up 2-monads out
of free ones exactly as one builds up algebraic structure using
basic operations, derived operations, and equations. Both the 
colimits and the freeness will involve the world of strict 
morphisms of monads. Exactly what this world might be is 
discussed below, but to start with we indicate why (strict) maps out
of a given monad are important.

\subsection{Endomorphism monads}
\label{sec:endomorphism-monads}

Let $T$ be a monad on a complete category \K. Everything
works without change for 2-categories, or indeed
for \V-categories.
For objects $A,B\in\K$, the right Kan extension
\[\xymatrix{ & \K \ar[dr]^{\< A,B\>}\\
  1 \ar[rr]_B \ar[ur]^A \rrtwocell\omit{<-3>} && \K}\]
can be computed as
\[\< A,B\> C = \K(C,A) \ct B\]
where $\ct$ means the cotensor, defined by
\[\K(D,X \ct B) \cong \Cat(X,\K(D,B))\]
for a set (or category or object of \V, as the case may be) $X$, and 
objects $B$ and $D$ of \K. The universal property of the right Kan
extension implies in particular that we have bijections of natural
transformations.
\[\frac{T\too \< A,B\>}{TA \too B}\]
(This is starting to look like something you might want to do if $T$ is
a monad.)

We have natural ``composition'' and ``identity'' maps
\[\< B,C\> \< A,B\> \too \< A,C\>\]
\[1 \to\<A,A\>\]
which provide \K\ with an enrichment over $[\K,\K]$ with internal-hom
$\<A,B\>$.  (Writing down where the composition and identity come from
is a good exercise.) Thus $\<A,A\>$ becomes a monoid in $[\K,\K]$;  
that is, a monad. This can be regarded as the monoid of endomorphisms
of $A$ in the $[\K,\K]$-category \K.

The important thing about this monad is that the bijection
\[\frac{T\too\<A,A\>}{TA\too B}\]
restricts to a bijection
\[\frac{T\too[\text{monad}] \< A,A\>}{TA \too[\text{alg. str.}] A}\]
between monad maps into $\<A,A\>$ and algebra structures on $A$.

This is exactly like the endomorphism operad of an object, except that
instead of an object of $n$-ary operations for each $n\in\NN$, we have
an object ``$C$-ary operations'' 
\[\<A,A\> C = \K(C,A)\ct A\]
for each object $C\in \K$.

Thus colimits of monads are interesting.  For simple example,
algebras for $S+T$ (coproduct as monads) are objects with an algebra
structure for $S$ and an algebra structure for $T$, with no particular 
relationship between the two. 

We can play the same game with morphisms. First observe that $\<A,B\>$
is functorial (covariant in $B$, contravariant in $A$), and so for
any $f:A\to B$ we can form the solid part of 
\[\xymatrix{T \ar@{.>}[r]^{\beta} \ar@{.>}[d]_{\alpha} & 
\<B,B\> \ar[d]^{\<f,B\>} \\
  \<A,A\> \ar[r]_{\<A,f\>} & \<A,B\>}\]
and now if we have monad maps $\alpha:T\to \<A,A\>$ and $\beta:T\to \<B,B\>$, 
then the square commutes if and only if $f$ is a strict map between 
the corresponding algebras.  

It is at this point that we want to make things 2-categorical, and 
allow for pseudo or lax morphisms. So suppose that \K is a (complete)
2-category, and that $T$ is a 2-monad on \K. To give a 2-cell
\[ \xymatrix{
TA \ar[r]^{Tf} \ar[d]_{a} & TB \ar[d]^{b} \\
  A \ar[r]_{f} \rtwocell\omit{<-4>\fbar} & B}  \]
is equivalently to give a 2-cell
\[\xymatrix{T \ar[r]^{\beta} \ar[d]_{\alpha}  & \<B,B\> \ar[d]^{\<f,B\>} \\
  \<A,A\> \ar[r]_{\<A,f\>} \rtwocell\omit{<-4>\tilde{f}} & \<A,B\>. }\]
The \emph{comma object} 
\[\xymatrix{\{f,f\}_\ell \ar[r]^{d} \ar[d]_{c} & \<B,B\> \ar[d]^{\<f,B\>} \\
  \<A,A\> \ar[r]_{\<A,f\>} \rtwocell\omit{<-4>}  & \< A,B\> }\]
is the universal diagram of this shape, so to give $\tilde{f}$ as
above is equivalent to giving a 1-cell $\phi:T\to\{f,f\}_\ell$
with $d\phi=\beta$ and $c\phi=\alpha$.

Now $\{f,f\}_\ell$ becomes a monad: this can be seen via a routine argument
using pasting diagrams; or one can get more sophisticated, and show
that \LaxtwoK is enriched over $[\K,\K]$, and now regard $\{f,f\}_\ell$ as
the endomorphism monoid. The important thing is that 
\[\xymatrix{\{f,f\}_\ell \ar[r]^{d} \ar[d]_{c}  & \<B,B\> \\
  \<A,A\>}\]
are monad maps (although $\< A,B\>$ is not a monad), and that
a map $T\to \{f,f\}_\ell$ is a monad map if and only if the corresponding
$(f,\fbar)$ is a lax $T$-morphism.  

Of course there is also a pseudo version of this: use the
{\em iso-comma object} $\{f,f\}$ rather than the comma object $\{f,f\}_\ell$;
this is the evident analogue in which the 2-cell is required to be invertible.

Thus we can work out the algebras and the (strict, pseudo, or lax) 
morphisms for a monad just by looking at monad morphisms out of $T$,
and {\em that} shows why free monads and colimits of monads should
be important.

\begin{xca}
Describe the $T$-transformations in this way.
\end{xca}

\subsection{Pseudomorphisms of monads}
\label{sect:monad-pseudomorphisms}

In addition to strict monad maps, where the good colimits live, there
are also pseudo maps of monads.  A \emph{pseudomorphism} of 2-monads
on \K\ is a 2-natural transformation which preserves the multiplication
and unit in exactly the sense that strong monoidal functors preserve
the tensor product and unit of monoidal categories. Thus there
are invertible 2-cells
$$\xymatrix{
1 \ar[r]^{i} \ar[dr]_{j} \ar@{}[drr]|(0.3){\cong} & T \ar[d]^{f} &
T^2 \ar[l]_{m} \ar[d]^{f^2} \ar@{}[dl]|{\cong} \\
& T & S \ar[l]^{n} }$$
satisfying the same usual coherence conditions.  

We have seen that to give an arbitrary map $\alpha:T\to\<A,A\>$ is
equivalent to giving $a:TA\to A$ in \K, and that $\alpha$ is a strict
map of monads if and only if $a$ makes $A$ into a strict algebra; it turns 
out that to make $\alpha$ into a pseudomorphism of monads 
\[\alpha:T\rightsquigarrow \< A,A\>\] 
is precisely equivalent to making 
$$a:TA\to A$$ 
into a pseudoalgebra.

\subsection{Locally finitely presentable 2-categories}
\label{sect:lfp}

For a large 2-category \K, the 2-category \Mnd of 2-monads on \K
has all sorts or problems: its hom-categories are large, it is
not cocomplete, and free monads don't exist. We shall therefore pass 
to a smaller 2-category of 2-monads.

Assume that \K\ is a locally finitely presentable
2-category. If you know what a locally finitely presentable category
is then this is just the obvious 2-categorical analogue. If not, then
here are some ways you could think about them:
\begin{itemize}
\item The formal definition (which you don't need to know because I'm not
  going to prove anything): a cocomplete 2-category with a small full
  subcategory which is a strong generator and consists of finitely
  presentable objects.
\item A 2-category which is complete and cocomplete and in which transfinite 
  arguments are more inclined to work then is usually the case.
\item A 2-category of all finite-limit-preserving 2-functors from \C to
  \Cat, where \C is a small 2-category with finite limits; you can
  take \C\ to be $\Kf\op$ where \Kf is the full subcategory of
  finitely presentable objects.
\item Full reflective sub-2-categories of presheaf 2-categories which are
  closed under filtered colimits.
\item A 2-category which is complete and cocomplete, and is the free 
  cocompletion under filtered colimits of some small 2-category (an 
  Ind-completion). In fact you don't need to suppose both completeness
  and cocompleteness: for an Ind-completion, either implies the other.
\end{itemize}

Examples include the presheaf 2-category $[\A,\Cat]$ for any small 2-category
\A, or $\Cat^X$ for any set $X$. The 2-category of groupoids is
another example.

Once again, this is really an enriched categorical notion: there
is a notion of locally finitely presentable \V-category, provided
that \V itself has a good notion of finitely presentable object: more
precisely, provided that \V is a locally finitely presentable category
and the full subcategory of finitely presentable objects is closed
under the monoidal structure. 

Because \K is the free completion of \Kf under filtered colimits,
to give an arbitrary 2-functor $\Kf\to\K$ is equivalent to giving
a finitary (that is, filtered-colimit-preserving) 2-functor
$\K\to\K$. We write \Endf for the monoidal 2-category of
finitary endo(-2-)functors on \K.  Unlike $[\K,\K]$ this is locally
small, since $\Kf$ is small.

A 2-monad is said to be finitary if its endo-2-functor part is so.
Then the 2-category \Mndf of finitary 2-monads on \K is the 2-category
of monoids in \Endf, and the forgetful 2-functor $U:\Mndf\to\Endf$
does indeed have a left adjoint, so in this world we do have free
monads. 

Moreover, the adjunction is monadic; there is a 2-monad on
\Endf for which \Mndf is the strict algebras and strict
morphisms.  We can drop down even further to get
$$\xymatrix{
\Mndf \dtwocell^{W}_{H}{'\dashv} \ar@/^10mm/[dd]^{U} \\
\Endf \dtwocell^{V}_{G}{'\dashv} \\
[\ob\Kf,\K] \ar@/^10mm/[uu]^{F} }$$
and go back up (along $G$) by left Kan extension along the inclusion
$\ob\Kf\to\K$.  The lower adjunction is also
monadic, as indeed is the composite, although this does not
follow from monadicity of the two other adjunctions.

Thus \Mndf is monadic both over \Endf and over $[\ob\Kf,\K]$, and the 
choice of which base 2-category to work over affects
what the pseudomorphisms and pseudoalgebras will be.  Dropping
down one level, the transformations are 2-natural ones, as in 
Section~\ref{sect:monad-pseudomorphisms}; while if we
drop down the whole way, they will be only pseudonatural.

The induced monads on $[\ob\Kf,\K]$ are finitary, and so it follows
that \Endf and \Mndf are themselves locally finitely presentable,
and in particular are complete and cocomplete.
In fact slightly more is true, since the inclusion of \Mndf in \Mnd has
a right adjoint, and so preserves colimits. 
\Mnd does not have colimits in general, but it does
have colimits of finitary monads, and these are finitary. Free monads
on arbitrary endo-2-functors may not exist, but free monads on finitary
endo-2-functors do, and they are themselves finitary. This is useful
since the $\<A,A\>$ are {\em not} finitary, although we can use the
coreflection of \Mnd into \Mndf to obtain a finitary analogue.

Everything in this section remains true if you replace ``finite'' by 
some regular cardinal $\alpha$.

\subsection{Presentations}
\label{sec:presentations}

The most primitive generator for a 2-monad is an object
of $[\ob\Kf,\K]$: a family $(X_c)_{c\in\ob\Kf}$ of objects of \K, indexed
by the objects of \Kf. This then generates a free 2-monad $FX$.
What is an $FX$-algebra?  A monad map
\[FX\to \< A,A\>\]
which is the same as 
\[X\to U\< A,A\>.\] 
This just means that for each $c$, we have
\[Xc\to \< A,A\> c\]
which unravels to a functor
\[\K(c,A) \to \K(Xc,A)\]
between hom-categories.  Since \K\ is cocomplete,
this is the same as a map
\[\sum_c \K(c,A)\cdot Xc \too A\]
where $\K(c,A)\cdot Xc$ denotes the tensor of $Xc\in\K$ by the
category $\K(c,A)$.
Thus we can think of $Xc$ as the ``object of all $c$-ary operations''.

\begin{example}
  Let $\K=\Cat$, so \Kf is the finitely presentable categories, and
  $X$ assigns to every such $c$ a category $Xc$ of $c$-ary operations.
  We take
  \[Xc = 
  \begin{cases}
    1 & c=0,2 \\
    0 & \text{otherwise}
  \end{cases}
  \]
  where 2 is the discrete category $1+1$.  Thus we have one binary operation
  and one nullary operation.  An $FX$-algebra is then a category $A$
  with maps as above.  If $Xc$ is empty, then $\K(Xc,A)$ is terminal,
  so there's nothing to do.  In the other cases, we get maps
  \[A^2 \to A\] when $c=2$ and
  \[A^0 = 1 \to A\] when $c=0$.  This is the first step along the path
  of building up the 2-monad for monoidal categories.  The pseudo (or
  lax) morphisms can be determined using $\{f,f\}$ or $\{f,f\}_\ell$:
  they will preserve $\ot$ and $I$ in the pseudo or the lax
  sense, as the case may be, but without coherence conditions.
\end{example}

\begin{example}
  Again let $\K=\Cat$, and let
  \[X c = 
  \begin{cases}
    \two & c=1\\
    0 & \text{otherwise}
  \end{cases}
  \]
  Then an $FX$ algebra is a category with a map
  \[A\to A^\two\]
  in other words, a pair of maps with a natural transformation
  \[\xymatrix{A \rtwocell & A}.\]
  This is an example in which $Xc$ is not discrete. We say that $X$
  specifies a ``basic operation of arity $1$ (unary) and type arrow''.
\end{example}

In the case of monoidal categories, you can use operations of ``type
arrow'' to provide the associativity and unit isomorphisms, but I'll
take a different approach.

\subsection{Monoidal categories}
\label{sec:monoidal-categories}

Actually, let's forget about the units, just worry about the binary
operation.  Then $Xc$ is $1$ if $c=2$ and $0$ otherwise, so an
$FX$-algebra is a category with a single binary operation.  Then we
have a (non-commutative) diagram of 2-categories and 2-functors

$$
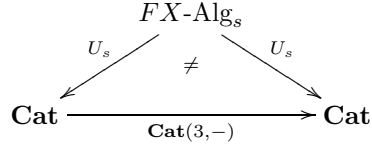
\xymatrix{
& {\fxalgs} \ar[dl]_{U_s} \ar[dr]^{U_s} \ar@{}[d]|{\neq} \\
\Cat \ar[rr]_{\Cat(3,-)} && \Cat}$$
which act on an $FX$-algebra $(C,\ot)$ by
$$
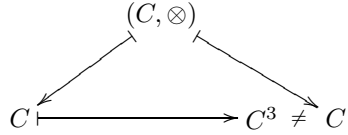
\xymatrix @C1pc { 
&& (C,\ot) \ar@{|->}[dll] \ar@{|->}[drr] \\
C \ar@{|->}[rrr] &&& C^3 \ar@{}[r]|{\neq} & C }$$
and now we have the two maps 
$$\xymatrix{
C^3 \ar@<1ex>[r]^{\ot(\ot 1)} \ar@<-1ex>[r]_{\ot(1\ot)} & C }$$
which are natural in $(C,\ot)$, and so induce two natural
transformations $\Cat(3,-)U_s\to U_s$ in the previous triangle.
We can take their mates under the adjunction $F_s\dashv U_s$ to get 2-cells in
$$
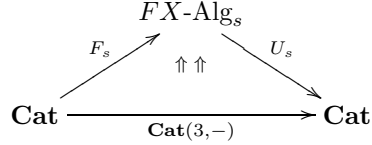
\xymatrix{
& {\fxalgs} \ar[dr]^{U_s} 
\dtwocell\omit{\omit \Uparrow \;\Uparrow } & {} \\
\Cat \ar[ur]^{F_s} \ar[rr]_{\Cat(3,-)} && \Cat}$$
with two 2-cells in the middle.  Note that $U_s F_s = FX$ is the
monad, so that we have two natural transformations
\[\Cat(3,-) \rightrightarrows FX,\]
which are morphisms of endofunctors. We can now construct the free
2-monad $H\Cat(3,-)$ on $\Cat(3,-)$ and the induced monad morphisms
\[\kappa_1,\kappa_2\maps H\Cat(3,-) \rightrightarrows FX.\]

Consider now an $FX$-algebra $(C,\ot)$, and the corresponding 
monad map $\gamma:FX\to\<C,C\>$. Then $(C,\ot)$ is strictly
associative if and only if $\gamma\kappa_1=\gamma\kappa_2$, while
to give an isomorphism $\ot(1\ot)\cong\ot(\ot 1)$ is equivalent
to giving an isomorphism $\gamma\kappa_1\cong\gamma\kappa_2$.
In the 2-category \Mndf construct the universal map $\rho:FX\to S$
equipped with an isomorphism $\rho\kappa_1\cong\rho\kappa_2$:
this is called a {\em co-iso-inserter}, and it's a (completely strict)
2-categorical colimit, which we'll meet later on.

Now, an $S$-algebra is a category $C$ equipped with a monad map
$S\to\<C,C\>$, or equivalently a monad map $\rho:FX\to\<C,C\>$ and
an isomorphism $\rho\kappa_1\cong\rho\kappa_2$, or equivalently with
a functor $\ot\maps C^2\to C$ and a natural isomorphism 
$\alpha:\ot(1\ot)\cong \ot(\ot 1)$.
You can also write down what it means to be a pseudo or lax morphism of
such algebras, and it's what you want it to be: the tensor-preserving
isomorphisms must be compatible with the associativity constraints.

The coherence condition states that a pair of 2-cells 
$$\xymatrix @C4pc {
C^4 \rtwocell^{\ot(\ot1)(\ot11)}_{\ot(1\ot)(11\ot)}{\omit\Downarrow\;\Downarrow} & C}$$
built up using $\alpha$ are equal. Much as before, these are all
natural in $(C,\ot,\alpha)$, so $\ot(\ot1)(\ot11)$ corresponds to one monad map
$\lambda_1:H\Cat(4,-)\to S$, and $\ot(1\ot)(11\ot)$ to another
$\lambda_2:H\Cat(4,-)\to S$. Furthermore the two isomorphisms 
$\ot(\ot1)(\ot11)\cong\ot(1\ot)(11\ot)$ correspond to invertible
monad 2-cells 
$$\xymatrix @C4pc {\relax\llap{$H\Cat(4,-)$} 
\rtwocell^{\lambda_1}_{\lambda_2}{\omit\Lambda_1\Downarrow\;\Downarrow\Lambda_2}
& S.~~~}$$
Now we form the \emph{coequifier} $q\maps S\to T$, in the category of
monads, of these two 2-cells $\Lambda_1$ and $\Lambda_2$: this is the 
universal map $q$ with the property
that $q\Lambda_1 = q\Lambda_2$.

Then the 2-category \talg is the 2-category of ``semigroupoidal
categories'' and strong morphisms (we can get the strict and lax
morphisms in the obvious way too).  All this follows from the
universal property of the monad $T$.

Often, as here, we build up structure in a particular order, starting with the 
operations of type object, then those of type arrow or
isomorphism, and finally impose equations on these arrows
or isomorphisms.

\subsection{Terminal objects}
\label{sect:terminal-objects}

Consider the structure of \emph{category with terminal object}. This 
is a baby example, but you can do any limits you like once you
understand this example. 

How do you say algebraically that a category $A$ has a terminal
object? You give an object
\[1\too[t] A\]
with a natural transformation
\[\xymatrix{A \rrtwocell\omit{<3>\tau} \ar[rr]^{1} \ar[dr]_{!} && A\\ 
& 1 \ar[ur]_{t} }\] 
such that the component 
\[\xymatrix{1 \ar[r]^t & A \rrtwocell\omit{<3>\tau} \ar[rr]^{1} \ar[dr] && A\\
  && 1 \ar[ur]}\]
of $\tau$ at $t$
is the identity. This last condition plus naturality of $\tau$ guarantees
that $\tau a:a\to t$ is the only map from $a$ to $t$, and so that $t$
is terminal.

Let's give a presentation for it.  
First we have the nullary operation $t$, which takes the form
\[\Cat(0,A) \to \Cat(1,A)\]
or equivalently
\[\Cat(0,A)\cdot X0 \to A\]
where $X0=1$, or equivalently
\[\sum_c \Cat(c,A)\cdot Xc \to A\]
where now $Xc$ is 0 unless $c=0$. Thus an object $A$ with nullary
operation $t:1\to A$ is precisely an $FX$-algebra, where 
$$Xc=\begin{cases}1 & \text{if $c=0$} \\ 0 & \text{otherwise}\end{cases}$$

For any $FX$-algebra $(A,t)$, there are two canonical maps
from $A\to A$, given by 
$$\xymatrix @R1pc {
A \ar[rr]^{1} \ar[dr]_{!} && A \\
& 1 \ar[ur]_{t} }$$
and these are clearly natural in $(A,t)$; in other words,
they define a pair of natural transformations from 
$U_s:\fxalgs\to\Cat$ to itself. Taking mates
under the adjunction $F_s\dashv U_s$ gives a pair of natural
transformations
$1_\Cat\to U_sF_s$. Now $U_sF_s$ is just $FX$, so forming the
free monad $H1$ on the identity $1_\Cat$, we get a pair of 
monad maps 
$\kappa_1,\kappa_2:H1\to FX$.
We now form the {\em coinserter} $\rho:FX\to S$ of $\kappa_1$ 
and $\kappa_2$. This is another 2-categorical colimit; it is the
universal $\rho$ equipped with a 2-cell $\rho\kappa_1\to\rho\kappa_2$.
An $S$-algebra is now an $A$ equipped with an object $t:1\to A$, and
a natural transformation $\tau:1_A\to t\circ !$, as in our earlier
description of terminal objects. Finally one can construct a suitable
coequifier $q:S\to T$ to obtain the 2-monad $T$ for categories with
terminal objects.

Here's a different presentation: it starts as before by putting in
a nullary operation
\[\Cat(0,A) \too[t] \Cat(1,A)\]
but then adds a unary operation of type arrow:
\[\Cat(1,A) \too[\tau] \Cat(\two,A)\]
which specifies two endomorphisms of $A$ and a natural transformation
between them:
\[\xymatrix{A \rtwocell^f_g{\tau} & A.}\]
Later we'll introduce equations $f=1$, $g=t\circ !$, and 
$\tau t=\id$.

To specify $t$ and $\tau$, define 
$$Xc = \begin{cases} 1 & \text{if $c=0$} \\
                     \two & \text{if $c=1$} \\
                     0 & \text{otherwise}
       \end{cases} $$
so that $FX$-algebra structure an a category $A$ amounts to 
$$\sum_c \Cat(c,A)\cdot Xc \to A$$
or equivalently $t:1\to A$ and $\tau:f\to g:A\to A$.

Now we turn to the equations. Consider the (non-commuting) diagram
$$\xymatrix{
& \fxalgs \ar[dr]^{U_s} \ar[dl]_{U_s} \ar@{}[d]|(0.5){\neq} \\
\Cat \ar[rr]_{2\cdot-+\two} && \Cat }$$
which acts on an $FX$-algebra $(A,t,\tau)$ by
$$\xymatrix{ 
& (A,t,\tau) \ar@{|->}[drr] \ar@{|->}[dl] \\
A \ar@{|->}[rr] && A+A+\two \ar@{}[r]|-{\neq} & A.}$$
There is a map $\alpha_{(A,t,\tau)}:A+A+\two\to A$ whose components 
are $f:A\to A$, $g:A\to A$, and the functor $\two\to A$ corresponding 
to $\tau\circ t$. This is natural in $(A,t,\tau)$.

There is another map $\beta_{(A,t,\tau)}:A+A+\two\to A$ whose
components are $1:A\to A$,
$t\circ !:A\to A$, and the functor $\two\to A$ corresponding to the 
identity natural transformation on $t$. Once again this is natural
in $(A,t,\tau)$. 

A category with terminal object is precisely an $FX$-algebra
$(A,t,\tau)$ for which $\alpha_{(A,t,\tau)}=\beta_{(A,t,\tau)}$.

Now $\alpha$ and $\beta$ live in the diagram
$$\xymatrix{
& \fxalgs \ar[dr]^{U_s} \ar[dl]_{U_s} 
\dtwocell\omit{\omit \alpha\Uparrow\; \Uparrow\beta} \\
\Cat \ar[rr]_{E} && \Cat }$$
where $EC=C+C+\two$,
and we can take their mates under the adjunction $F_s\dashv U_s$ to 
obtain natural transformations
$$\alpha',\beta':E \to U_sF_s$$
and now $U_sF_s$ is the monad $FX$, so there are induced monad maps
$$\bar{\alpha},\bar{\beta}:HE\to FX$$
from the free monad $HE$ on $E$,
and the required 2-monad $T$ for categories-with-terminal object is
obtained as the coequalizer 
$$\xymatrix{
HE \ar@<1ex>[r]^{\bar{\alpha}} \ar@<-1ex>[r]_{\bar{\beta}} & 
FX \ar[r]^{q} & T. }$$

\begin{remark} Whichever approach we take, the algebras will be the
categories with a \emph{chosen}
terminal object.  This may seem strange, but is not really a problem. The
strict morphisms preserve the chosen terminal object strictly, which
is probably not what we really want, but the pseudo morphisms preserve 
it in the usual sense.
\end{remark}

\subsection{Bicategories}
\label{sect:bicategories}

There are two reasons for including this example: first of all it's a
fairly easy case with $\K\neq\Cat$, and second it's important for 2-nerves.
I won't give all the details.

Let $\K=\Cat\text{-}\mathrm{Grph}$, the 2-category of
category-enriched graphs.  A \Cat-graph consists of objects
$G,H,\dots$ and hom-categories $\G(G,H)\in\Cat$.  (Of course one could 
do this for any \V in place of \Cat.)  A morphism is a function 
$G\mapsto FG$ on objects, along with functors $\G(G,H)\to\H(FG,FH)$ between
hom-categories. One might hope that the 2-cells would
be some sort of natural transformations, but since \Cat-graphs have no 
composition law, there is no way to assert that a square in a \Cat-graph 
commutes, and so no way to state naturality. Instead, we use a special
sort of lax naturality. We only allow 2-cells
\[\xymatrix{\G \rtwocell^F_{F'} & \H}\]
to exist when $F$ and $F'$ agree on objects, and then the 2-cell consists of
natural transformations
$$\xymatrix{\relax\llap{$\G(G,H)$} \rtwocell ^F_{F'} &
\relax\rlap{$\H(FG,FH)$} }$$
on all hom-categories. 

Now, given a \Cat-graph, what do you need to do to turn it into a
bicategory?  To start with, you have to give compositions
\[\G(H,K) \times \G(G,H) \too \G(G,K).\]
Let $\gcomp$ and $\garr$ be the \Cat-graphs $\cdot\to\cdot\to\cdot$ 
and $\cdot\to\cdot$ (no 2-cells).  
Then 
\[\K(\gcomp,\G) = \sum_{G,H,K} \G(H,K) \times \G(G,H).\]
\[\K(\garr,\G) = \sum_{G,K} \G(G,K)\]
so if we define 
$$Xc = \begin{cases}\garr & \text{if $c=\gcomp$} \\
                    0 & \text{otherwise}\end{cases}$$ 
then an $FX$-algebra structure on \G amounts to a map
$$\sum_c \K(c,\G)\cdot Xc \to \G$$
and so to a map
$$M:\sum_{G,H,K}\G(H,K)\t\G(G,H) \to \sum_{G,K}\G(G,K).$$
We need to make sure that the restriction
$$M_{G,H,K}:\G(H,K)\t\G(G,H)\to\sum_{G,K}\G(G,K)$$
to the $(G,H,K)$-component lands in the $(G,K)$-component: this can be
done by constructing a quotient of $FX$.

Define
$$Yc = \begin{cases}\gob & \text{if $c=\gcomp$} \\
                    0   & \text{otherwise} \end{cases} $$
where \gob denotes the \Cat-graph $\cdot$ (no 1-cells or 2-cells).
An $FY$-algebra structure on a \Cat-graph \G is a map 
$$\sum_{G,H,K}\G(H,K)\t\G(G,H)\to\sum_{G}1$$
where $1$ denotes the terminal category.

Suppose now that $(\G,M)$ is an $FX$-algebra. There are many 
induced $FY$-algebra structures on \G; in particular, there
are the following two:
$$\xymatrix @C3pc @R1pc {
{\sum_{G,H,K} \G(H,K)\t\G(G,H)} \ar[r]^-{M} & 
{\sum_G\sum_K\G(G,K)} \ar[r]^-{\sum_G !} & {\sum_G 1} \\
{\sum_{G,H,K}\G(H,K)\t\G(G,H)} \ar[r]^-{!} & 1 \ar[r]^-{\text{inj}_G} & 
{\sum_G 1} }$$
Each is functorial, and so each induces a monad map $FY\to FX$; we
form their coequalizer $q_1:FX\to S_1$, and now an $S_1$-algebra
is a \Cat-graph \G equipped with a composition $M$ such that
$M_{G,H,K}$ lands in $\sum_K \G(G,K)$. A further quotient forces
$M_{G,H,K}$ to land in $\G(G,K)$ as desired.

One now introduces an associativity isomorphism. This has the form
of a map 
$$\K(\gtriple,\G) \to \K(\giso,\G)$$
where $\gtriple$ is the \Cat-graph $\cdot\to\cdot\to\cdot\to\cdot$ and 
$\giso$ is $\xymatrix{{\cdot} \rtwocell{\omit\cong} & {\cdot} }$.
There are also left and right identity isomorphisms, and various
coherence conditions to be encoded, but I'll leave all that as an exercise.
The result of the exercise is:
\begin{itemize}
\item An algebra is a bicategory.
\item A lax morphism is a lax functor.
\item A pseudo morphism is a pseudo functor.
\item A strict morphism is a strict functor.
\item A 2-cell is an \emph{icon}.  This is an oplax natural
  transformation (which we haven't officially met yet) for which the
  1-cell components are identities.  ICON stands for ``Identity
  Component Oplax Natural-transformation''. An icon $F\to G$ can exist
only if $F$ and $G$ agree on objects, in which case it consists of a 
2-cell 
\[\xymatrix @C1pc {FA \ar@{=}[rr]\ar[d]_{Ff} && GA \ar[d]^{Gf} \\
  FB \ar@{=}[rr] & {}\utwocell\omit & GB}\]
for each $f:A\to B$ in \G, subject to conditions expressing compatibility
with respect to composition of 1-cells and identities, and naturality in $f$
with respect to 2-cells.
In the case of one-object bicategories these are precisely the 
monoidal natural transformations.
\end{itemize}

These icons are just nice enough to give us a 2-category of
bicategories.  In general, lax natural transformations between lax
functors can't even be whiskered by lax functors --- the composite
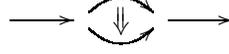
\[\xymatrix{{\ar[r]} & {\rtwocell} & {\ar[r]} & }\] 
isn't well-defined. In the pseudo case it is defined, but not associative,
and so we are led into the world of tricategories. But with just icons, we do
get a 2-category, which is moreover the category of algebras for the
2-monad just described.  

For example, in this 2-category, it's true that every bicategory is
equivalent (in the 2-category) to a 2-category; this works because in
replacing a bicategory by a biequivalent 2-category you don't have to 
change the objects of the bicategory. The 2-category \NHom of bicategories,
normal homomorphisms, and icons, is a full sub-2-category of the the 
2-category $[\DD\op,\Cat]$ of simplicial objects in \Cat, via a
``2-nerve'' construction. In order to deal with normal homomorphisms 
(which preserve identities strictly) rather than general ones, it's
convenient to start with {\em reflexive} \Cat-graphs rather than
\Cat-graphs.

The choice of direction of the 2-cell in lax transformations and 
oplax transformations goes back to B\'enabou. It seems that the
oplax transformations are generally more important than the lax 
ones.

\subsection{Cartesian closed categories}
\label{sect:ccc}

The comments in this section apply equally to monoidal closed categories,
symmetric monoidal closed categories, and toposes.

There is no problem constructing a monad for categories with finite 
products, similarly to the constructions given above. When we come 
to the closed structure, however, things are not so straightforward.
The internal hom is a functor
\[A\op\times A \to A\]
and we're not allowed to talk about $A\op$ the way we're doing
things: our operations are supposed to be of the form $A^c\to A$.
How can we deal with this?

In fact, it's a theorem that cartesian closed categories \emph{don't} have 
the form \talg\ for a 2-monad $T$ on \Cat. What you can do, however, is
change the base 2-category \K to the 2-category $\Cat_g$ of categories,
functors, and natural 
\emph{isomorphisms}.  Recall that $\Cat(2,A)$ is just $A\times A$, but
in $\Cat_g(2,A)$ we have only $A_{\iso} \times A_{\iso}$, where 
$A_{\iso}$ is the subcategory of $A$ consisting of the isomorphisms.  The
internal-hom \emph{does} give us a functor
$$\xymatrix @R0.5pc { 
A_{\iso} \times A_{\iso} \ar[r] & A_{\iso} \\
(a,b) \ar@{|->}[r] & [a,b] }$$
which has the form
\[\Cat_g(2,A) \too \Cat_g(1,A)\]
since we can turn around an isomorphism in the first variable to make
everything covariant.  This gives a new problem; the product
is now only given as a functor $A_{\iso}\t A_{\iso}\to A_{\iso}$, we have
to put in the rest of the functoriality separately ``by hand'', using
an operation
$$\xymatrix @R0.5pc {
{}\Cat_g(\two + \two, A) \ar[r] &  \Cat_g(\two,A) \\
(f:a\to a',g:b\to b') \ar@{|->}[r] & (f\t g:a\t b\to a'\t b') }$$ 
subject to various equations. You also have to relate the product
to the internal hom.

Any 2-monad on  \Cat\ induces monads on
$\Cat_g$ and on the 1-category $\Cat_0$ (since things are stable under
change of base monoidal category: categories to groupoids to sets).
But at each stage, to present the same structure becomes harder. In
the groupoid enriched stage we can still talk about pseudomorphisms, although
at this stage every lax morphism is pseudo; by the time we get to 
the \Set-enriched stage there is no longer any genuine pseudo notion 
at all --- everything is strict.

\subsection{Diagram 2-categories}
\label{sec:diagram-2-categories}

The first version of this is not really an example of a presentation
at all, since the 2-monad pops out for free.  Let \C\ be a small
2-category, and consider the 2-category $[\C,\Cat]$ of (strict)
2-functors, 2-natural transformations, and modifications.  This is the
\Cat-enriched functor category.  The forgetful 2-functor has both adjoints
$$\xymatrix @R3pc {
[\C,\Cat] \ar[d] \ar@{}[d]^{~~~\dashv}_{\dashv~~} \\
[\ob\C,\Cat] \uuppertwocell{\omit} \ulowertwocell{\omit} }$$ 
given by left and right Kan extension.  The
existence of the right adjoint tells us that the forgetful functor
preserves all colimits.  In this case $U_s$ is strictly monadic
as is easily proved using the enriched version of Beck's
theorem.  The induced monad $T$ then preserves \emph{all} colimits,
and we can write, using the Kan extension formula,
\[(TX)c = \sum_d \C(d,c)\cdot Xd.\]

It's now a long, but essentially routine, exercise to check that
\begin{itemize}
\item pseudo $T$-algebras are pseudo-functors,
\item lax algebras are lax functors,
\item pseudo morphisms are pseudo-natural transformations,
\end{itemize}
and so on. When you write down the coherence conditions for a lax morphism
it will tell you more than is in the \emph{definition} of a lax
functor: it will also include a whole lot of consequences of the 
definition. Notice, by the way, that this is a 2-monad for which
every pseudoalgebra is equivalent to a strict one, and so every 
pseudofunctor from \C to \Cat is equivalent to a strict one.

Now let \C\ be a bicategory.  If we tried the same game, we
wouldn't get a 2-monad, since the associativity of the multiplication
for the monad corresponds to the associativity of composition in \C,
so we'd just get a pseudo-monad.  We could just go ahead and do this, 
but we've been avoiding pseudo-monads, and there is an alternative.
One can give a 
presentation for a 2-monad $T$ on $[\ob\C,\Cat]$ whose
\begin{itemize}
\item (strict) algebras are pseudofunctors $\C\to\Cat$,
\item pseudomorphisms of algebras are pseudonatural transformations,
\end{itemize}
and so on. You start with a family $(X_c)_{c\in\ob\C}$, then introduce 
operations
$$\C(c,d)\t X_c\to X_d$$
and so on.
The target doesn't really need to be \Cat, although it would need
to be cocomplete.

\subsection{References to the literature}
For locally finitely presentable categories
see \cite{Gabriel-Ulmer} or \cite{AR}; and for the enriched version see
\cite{Kelly-amiens}. 
A more formal approach to presentations for 2-monads can be found in
\cite{KP}, which in turn builds upon ideas of \cite{Dubuc-Kelly} 
See \cite{mnd} for the monadicity of \Mndf over $[\ob\Kf,\K]$.

\section{Limits}
\label{sect:limits}

We'll begin with some concrete examples of limits, looking in particular
at limits in \talg, for a finitary 2-monad $T$ on a complete and 
cocomplete 2-category \K (you could get by with much
less for most of this). Recall that \talg is the 2-category of 
strict algebras and pseudomorphisms. A good example to bear in mind
would have $\K=\Cat$, and \talg the 2-category of categories with 
chosen limits of some particular type, and functors which preserve
these limits in the usual, up-to-isomorphism, sense.

\subsection{Terminal objects}
\label{sect:terminal-objects-1}

Let's start with something really easy: terminal objects.  Let $1$ be
terminal in \K; we have a unique map $T1\to 1$, making $1$ a
$T$-algebra, and then for any $T$-algebra $(A,a:TA\to A)$ we have a 
unique $!:A\to 1$, and
\[\xymatrix{TA \ar[r]^{T!}\ar[d] & T1 \ar[d]\\
  A\ar[r]_{!} & 1}\]
commutes strictly, so there's a unique \emph{strict} algebra morphism
$A\to 1$.  Moreover, by the 2-universal property of $1$, there's a
unique isomorphism in the above square, which happens to be an
identity; thus there is only one pseudo morphism as well (which
happens to be strict).  A similar argument works for endomorphisms of
this morphism; thus
\[\talg((A,a), (1,!)) \cong 1\]
so $(1,!)$ is a terminal object in \talg.

\subsection{Products}
\label{sec:products}

Similarly for products: given $T$-algebras $(A,a)$ and $(B,b)$, and 
a product $A\times B$ in \K, there is 
an obvious map $\<a,b\>$ as in 
\[T(A\times B) \to TA \t TB \to A\times B\]
which makes $A\t B$ into a $T$-algebra (exactly as for ordinary monads:
nothing 2-categorical going on here).  The point is that if we have 
\emph{pseudo} morphisms
\[\xymatrix{TC \ar[r] \ar[d]_{c} \ar@{}[dr]|{\cong} & TA \ar[d]^{a} &&
TC \ar[r] \ar[d]_{c} \ar@{}[dr]|{\cong} & TB \ar[d]^{b} \\
  C \ar[r] & A && C \ar[r] & B }\] \\
we get a unique induced pseudo morphism
\[\xymatrix{TC \ar[r]\ar[d]_{c} \ar@{}[dr]|{\cong} & 
T(A\times B)\ar[d]^{\<a,b\>} \\
  C\ar[r] & A\times B}\]
and indeed there is a natural isomorphism (of categories)
\[\talg(C,A\times B) \cong \talg(C,A)\times \talg(C,B).\]
Thus $A\times B$ is a product in \talg in the strict \Cat-enriched
sense.

Note that the projections $A\times B\to A$ and $A\times B\to B$ are
actually strict maps, by construction.  Moreover, they jointly ``detect
strictness'': a map into $A\times B$ is strict if and only if its
composites into $A$ and $B$ are strict.  This is a useful technical
property.

Actually, we didn't really need to check anything, since we've already
seen that $\talgs\hookrightarrow \talg$ has a left adjoint, hence preserves 
all limits, and in the case of terminal objects and products the diagram
of which we are taking the limit consists only of objects, so already exists
in the strict world. (On the other hand, the explicit argument works for
any 2-monad on any 2-category with the relevant products, whereas the 
adjunction needs a transfinite argument, and much stronger assumptions on
$T$ and \K.)

\subsection{Equalizers}
\label{sect:equalizers}

Now let's look at equalizers.  Here it's different, because the
morphisms whose equalizer we seek may not be strict.  If they {\em are}, 
then the equalizer exists in $\talgs$ and is preserved, but if they aren't,
the adjunction doesn't help.  In fact, in general equalizers of pseudo
morphisms need \emph{not} exist.

For example, let $T$ be the 2-monad on \Cat\ for categories with a
terminal object.  Let $1$ be the terminal category and let $\I$ be
the free-living isomorphism, consisting of two objects and a single
isomorphism between them.  Clearly both categories have a terminal object, and
both inclusions are pseudo morphisms.  But any functor which equalizes
them has to have empty domain, and no category with an empty domain
has a terminal object.

Thus \talg\ is not complete, but we can look at some
of the limits that it does have.  

\subsection{Equifiers}
\label{sect:equifiers}

Consider 
a parallel pair of 1-cells in \talg with a parallel pair of 2-cells
between them:
\[\xymatrix{ A \rrtwocell^{f}_{g}{\omit \alpha \Downarrow \;\Downarrow \beta} 
& & B. }\]
The {\em equifier} of these 2-cells, is the universal 1-cell 
$k:C\to A$ with $\alpha k=\beta k$. Here universality means that $\K(D,C)$ 
is \emph{isomorphic} (not just equivalent) to the category of morphisms 
$D\too[h] A$ with $\alpha h = \beta h$. Equifiers do lift from \K to \talg:
if $(A,a)$ and $(B,b)$ are $T$-algebras, $(f,\bar{f})$ and $(g,\bar{g})$
are $T$-morphisms, and $\alpha$ and $\beta$ are $T$-transformations, then 
the composites
$$\xymatrix @C1.5pc {
TC \ar[r]^{Tk} & TA \rrtwocell^{Tf}_{Tg}{~~T\alpha} && TB \ar[r]^{b} & 
B \ar@{}[r]|{=} &  
TC \ar[r]^{Tk} & TA \rrtwocell^{Tf}_{Tg}{~~T\beta} && TB \ar[r]^{b} &B}$$
are equal. Paste the isomorphism $\bar{g}:b.Tg\cong ga$ on the bottom
of each side and the isomorphism $\bar{f}:fa\cong b.Tf$ on the top, and 
use the $T$-transformation condition for $\alpha$ and $\beta$ to get
the equation
$$\xymatrix @C1.5pc {
TC \ar[r]^{Tk} & TA \ar[r]^{a} & A \rrtwocell^{f}_{g}{\alpha} && 
B \ar@{}[r]|{=} & TC \ar[r]^{Tk} & TA \ar[r]^{a} & 
A \rrtwocell^{f}_{g}{\beta} && B }$$
and now by the universal property of the equifier $C$ there is a unique
$c:TC\to C$ satisfying $kc=a.Tk$. Two applications of the universal
property show that $c$ makes $C$ into a $T$-algebra, and so clearly 
$k$ becomes a strict $T$-morphism $(C,c)\to(A,a)$. Further judicious
use of the universal property shows that $k:(C,c)\to(A,a)$ is indeed
the equifier in \talg.

Observe that once again, the projection map $k$ of the limit is
actually a strict map, and detects strictness of incoming maps.

Why does the analogous argument for equalizers fail?  Given 
pseudo morphisms $(f,\fbar)$ and $(g,\gbar)$ from $(A,a)$ to $(B,b)$,
we could form the equalizer $k:C\to A$ of $f$ and $g$, and then hope
to make $C$ into a $T$-algebra using the universal property of $C$,
but we'd need to show that $fa.Tk=ga.Tk$. All we actually know is
that $c.Tf.Tk=c.Tg.Tk$ while $c.Tf.Tk\cong fa.Tk$ and $c.Tg.Tk\cong ga.Tk$,
which just isn't good enough.

The moral is that in forming limits in \talg, we can ask for existence
or invertibility of 2-cells, and equations between them, but we can't
generally force equations between 1-cells.

\subsection{Inserters}
\label{sect:inserters}

There is a sort of lax version of an equalizer, called an inserter.
Rather than making 1-cells equal, you put a 2-cell in between them.
The {\em inserter} of a parallel pair of arrows $f,g:A\to B$ is the 
universal $k:C\to A$ equipped with a 2-cell $\kappa:fk\to gk$. 
More precisely, the universal property states that $\K(D,C)$ should be 
isomorphic to the category whose objects are morphisms $\ell:D\to A$ equipped 
with a 2-cell
$\lambda:f\ell\to g\ell$, and whose morphisms 
$(\ell,\lambda)\to(m,\mu)$ are 2-cells
\[\xymatrix{D\rtwocell^\ell_m{\alpha} & A}\]
such that
\[\xymatrix{f\ell \ar[r]^\lambda \ar[d]_{f\alpha} & g\ell \ar[d]^{g\alpha}\\
  fm\ar[r]_\kappa & gm}\]
commutes. Thus for every pair $(l:D\to A, \lambda:fl\to gl)$, there
is a unique $l':D\to C$ with $kl'=l$ and $\kappa
l'=\lambda$. Furthermore, given $l,\lambda, m,\mu,\alpha$ as above,
there is a unique $\alpha':l'\to m'$ with $k\alpha'=\alpha$.

Once again, inserters in \K\ lift to \talg, where they have strict
projections and detect strictness.  Given a pair
\[\xymatrix{(A,a) \ar@/^4mm/[r]^{(f,\fbar)} \ar@/_4mm/[r]_{(g,\gbar)} & (B,b)}\]
of pseudo morphisms, we construct the inserter $(k:C\to A,\kappa:fk\to gk)$ of 
$f$ and $g$ in \K, and want to make it an algebra.  We need a 2-cell
$f.a.Tk\to g.a.Tk$ to induce $c:TC\to C$, so we follow our nose:
$$\xymatrix{
f.a.Tk \ar[r]^{\bar{f}^{-1}.Tk} &  b.Tf.Tk \ar[r]^{b.T\kappa} &
 b.Tg.Tk \ar[r]^{\bar{g}.Tk} & g.a.Tk. }$$ 
This composite
must be $\kappa c$ for a unique $c$, by the universal property of
the inserter in \K. Now check that $c$ makes $C$ into an algebra, and so
on; everything goes through just as before.

Observe that an inserter in a (2-)category with no non-identity 2-cells is
just an equalizer.

\subsection{PIE-limits}
\label{sect:pie-limits}

Thus \talg\ has Products, Inserters, and Equifiers, and many 
important types of limit can be constructed out of these. A limit
which can be so constructed is called a {\em PIE-limit}, so clearly
\talg\ has all PIE-limits, and equally clearly equalizers are not
PIE-limits. Some other examples of PIE-limits are:
\begin{itemize}
\item \emph{iso-inserters}, which are inserters where we ask the
  2-cell to be invertible.  {\em Insert} 2-cells in each direction,
  then {\em equify} their composites with identities. (Of course you can't
  go the other way: iso-inserters don't suffice to construct inserters.)
\item \emph{inverters}, where we start with a 2-cell $\alpha$ and make
  it invertible: we want the universal $k$ such that $\alpha k$ is
  invertible. {\em Insert} something going back the
  other way, then {\em equify} composites with the identities.  
\item \emph{cotensors} by categories. Cotensors by discrete categories
can be constructed using {\em products}. Any category can be constructed from
discrete ones using coinserters (to add morphisms) and coequifiers (to
specify composites). So cotensors by arbitrary categories can be constructed
from cotensors by discrete categories using {\em inserters} and {\em
equifiers}.  
\end{itemize}

The dual (colimit) notions of {\em coinserter}, {\em coequifier}, 
and {\em co-iso-inserter}  were important in 
giving presentations of monads. The dual of inverter is the 
{\em coinverter}. The coinverter of a 2-cell $\alpha:f\to g:A\to B$ 
is the universal $q:B\to C$ with $q\alpha$ invertible. In \Cat,
this is just the category of fractions $B[\Sigma\inv]$, where
$\Sigma$ consists of all arrows in $B$ which appear as components
of $\alpha$. Of course the dual of cotensor is tensor, not cocotensor!

\subsection{Weighted Limits}
\label{sec:weighted-limits}

In this section we briefly review the general notion of weighted
limit, before turning in the next section to the case $\V=\Cat$, 
where we shall see how the various examples of the previous 
section arise.

Let $S:\C\to\K$ be a functor between, say, ordinary categories.
The limit is supposed to be defined by the fact that 
$$\K(A,\lim S) \cong \mathrm{Cone}(A,S)$$
where the right hand side is the set of cones under $S$ with vertex $A$.
This is typically defined as the hom-set $[\C,\K](\Delta A,S)$, where
$\Delta A$ denotes the constant functor at $A$, but it 
can also be expressed as $[\C,\Set](\Delta 1,\K(A,S))$. It is
this last description of cones which forms the basis for the 
generalization to weighted limits; we're going to replace $\Delta 1$
by some more general functor $\C\to\Set$.

\begin{example}
  No one really uses this in practice, but it's useful to think
  about, and motivates the name ``weighted'' in ``weighted limit''.  
  Let $\C=2$ have two objects, so a functor $S\maps \C\to \K$
  is a pair of objects $B$ and $C$, and a weight is a functor $J\maps
  \C\to \Set$, say it sends one to 2 and the other to 3.  Then
  \[[\C,\Set](J,\K(A,S))\] consists of functions $2\to \K(A,B)$ and $3\to
  \K(A,C)$, or equivalently two arrows $A\to B$ and three arrows $A\to C$,
  so that the ``weighted product'' is $B^2\times C^3$.
\end{example}

For general \V, we start with \V-functors 
$S\maps \C\to\K$ and $J\maps \C\to\V$ and consider
\[[\C,\V](J,\K(A,S)).\]
If this is representable as a functor of $A$, the representing object
is called the \emph{$J$-weighted limit} of $S$ and written $\{J,S\}$.
Thus we have a natural isomorphism
\[\K(A,\{J,S\}) \cong [\C,\V](J,\K(A,S)).\]
which defines the limit.

\begin{xca}
If $\K=\V$, then $\{J,S\}$ is the \V-valued hom $[\C,\V](J,S)$.
\end{xca}

When $\V=\Set$, weighted limits don't
give you any \emph{new} limits: if
\K\ is an ordinary category which is complete in the usual sense of
having all conical limits ($J=\Delta 1$), then it also has all
weighted limits. More precisely, for any weight $J:\C\to\Set$
and any diagram $S:\C\to\K$, there is a category \D and a diagram
$R:\D\to\K$, such that the universal property of $\{J,S\}$ is precisely
the universal property of the usual limit of $R$.

But the weighted ones are more expressive, so
it's still useful to think about them. In particular, you might want
to talk about all limits indexed by a particular weight $J:\C\to\Set$;
this class is not so easy to express using only conical limits.

When $\V\neq\Set$ it's not longer true that all limits can be
reduced to conical ones.  But if you have all
conical limits \emph{and} cotensors, you can construct all weighted
limits.  

\begin{remark}
There is a slight subtlety here. In the case $\V=\Set$, the conical
limit of a functor $S:\C\to\K$ is just the limit of $S$ weighted by
$\Delta 1:\C\to\Set$. But for a \V-category \C, the ``constant functor
at 1'' (from \C to \V) is usually not what you want to look at, and 
indeed may fail to exist. What you really want, to get the right
universal property, is the constant functor $\Delta I$ at the unit object
$I$ of \V. But even this may not exist, unless \C\ is the free \V-category
on an ordinary category \B. So this is the right general context for
conical limits in enriched category theory.
\end{remark}

That's all I want to say about general \V.

\subsection{\Cat-weighted limits}
\label{sect:cat-weighted-limits}

Here I describe the weights for some of the limit notions introduced
earlier.

\begin{example}[Inserters]
  Let \C be the 2-category $\cdot\rightrightarrows\cdot$, so $S$ is 
  determined by a parallel
  pair of arrows $A\rightrightarrows B$.  The weight $J\maps\C\to
  \Cat$ has image $(1\rightrightarrows \two)$.  Then a
  natural transformation $J\to \K(C,S)$ gives has two components. The
first is a functor $1\to \K(C,A)$, or equivalently a morphism
$h\maps C\to A$, while the second is a functor $\two\to \K(C,B)$, or 
equivalently a 2-cell
  \[\xymatrix{C\rtwocell^u_v{\beta} &B}.\]
Naturality of these components means precisely that $u=fh$ and $v=gh$, 
so the data consists of 1-cell $h\maps C\to A$ and a 2-cell 
$\beta\maps fh\to gh$.
To give $h$ and $\beta$ is just to give a map from $C$ into the inserter.
\end{example}

This is the 1-dimensional aspect of the universal property, which 
characterizes the 1-cells into $C$; there is also a 2-dimensional aspect
characterizing the 2-cells, since the limit is defined in terms of an
isomorphism of categories, not just a bijection between sets. In general,
this 2-dimensional aspect must be checked, but if the 2-category \K should
admit {\em tensors}, the 2-dimensional aspect follows from the 1-dimensional
one. Similar comments apply to all the examples.

\begin{example}[Equifiers]
  Here, our 2-category \C\ is
  \[\xymatrix{{\phantom{1}} \rrtwocell{\omit \alpha\Downarrow\;\Downarrow\beta}
  && {\phantom{\two}} }\]
  and our weight is
  \[\xymatrix{1 \rrtwocell{\omit \Downarrow\;\Downarrow} && \two}\]
  in which $\alpha$ and $\beta$ get mapped to the same 2-cell in \Cat.
\end{example}

\begin{example}[Comma objects]
  \C\ is the same shape as for pullbacks
  \[\xymatrix{ &  \ar[d]\\ \ar[r] & }\]
  and $J$ is
  \[\xymatrix{ &  1\ar[d]^1\\ 1\ar[r]_0 & \two}.\]
  There is no 2-cell in \C, since we don't \emph{start} with a 2-cell,
  we only add one universally.
\end{example}

\begin{example}[Inverters]
  Recall, this is where we start with a 2-cell and universally make it 
  invertible.  Then \C\ is
  \[\xymatrix @C3pc { {\phantom{1}}\rtwocell & {\phantom{\I}} }\]
  and $J$ is
  \[\xymatrix @C3pc {1 \rtwocell & \I}\]
  where \I\ is the ``free-living isomorphism'' $\cdot\rightleftarrows\cdot$.
\end{example}

\subsection{Colimits}
\label{sect:colimits}

Colimits in \K\ are limits in $\K\op$.  That's really all you have to
say, but I should show you the notation. As usual, we rewrite things
so as to refer to \K. In this case, it's also convenient to replace 
\C\ by $\C\op$, so that we start with 
\begin{align*}
  S &\maps \C\to \K\\
  J &\maps \C\op \to \V
\end{align*}
and now the weighted colimit is written $J\star S$ and
defined by a natural isomorphism
$$\K(J\star S,A)\cong[\C\op,\V](J,\K(S,A)).$$

One form of the {\em Yoneda lemma} says that
\[J\cong J\star Y\]
where $Y\maps \C\to [\C\op,\V]$ is the Yoneda embedding and 
$J:\C\op\to\V$ is arbitrary.  

Here's an application.
Suppose you have some ``limit-notion'' which you know in
advance is a weighted limit, but you don't know what the weight is.
Thus you know $\{J,S\}$ given $S$, but you don't know $J$ itself.
Consider the version of the Yoneda embedding $Y:\C\to{} [\C,\V]\op$ and 
take its ``limit'', for the notion of limit we're interested in; equivalently,
take the relevant {\em colimit} of $Y:\C\to[\C\op,\V]$. This is 
$J\star Y$ for our as yet unknown $J$; but by the Yoneda lemma this
$J\star Y$ is itself the desired weight. This can be used to 
calculate the weights for all the concrete examples of \Cat-weighted limits
discussed here.

\subsection{Pseudolimits}
\label{sec:pseudolimits}

The pseudolimit of a 2-functor $S:\C\to\K$ is defined an object 
$\pslim S$ of \K equipped with an isomorphism
\[\K(A,\pslim S) \cong \Ps(\C,\Cat)(\Delta 1, \K(A,S))\]
of categories natural in $A$,
where $\Ps(\A,\B)$ is the 2-category of 2-functors,
pseudonaturals, and modifications from \A\ to \B.  The right
side is what we mean by a \emph{pseudo-cone}.  Note that this is still
an \emph{isomorphism} of categories, not an equivalence, so such 
pseudolimits are determined up to isomorphism not just equivalence.

\begin{example}[Pseudopullbacks]
  Again we take \C\ to be
  
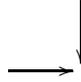
\[\xymatrix{ &  \ar[d]\\ \ar[r] & }\]
  A pseudo-cone then consists of
  
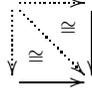
\[\xymatrix{ \ar@{.>}[r]\ar@{.>}[d] \ar@{.>}[dr]^\cong_\cong
    &  \ar[d]\\
    \ar[r] & }\]
  with isomorphisms in each triangle.  We have made the cones commute
  only up to isomorphisms, but the universal property and
  factorizations are still strict.  Note that the pseudopullback is
  \emph{equivalent} (not isomorphic) to the \emph{iso-comma object}
  (assuming both exist).  In the latter, we specify $fa\cong gb$
  without specifying the middle diagonal arrow.  Of course, we can
  take it to be $fa$, or $gb$, so we get ways of going back and forth.

  The pseudopullback is {\em not} in general equivalent to the
  pullback, although
  it is possible to characterize when they are \cite{JS-pseudopullbacks}. 
  This situation is entirely analogous to homotopy pullbacks, and indeed
  it can be regarded as a special case, via the ``categorical'' Quillen
  model structure on \Cat (see Section~\ref{sect:Quillen}).
\end{example}

Again, given a weight $J\maps \C\to \Cat$, the \emph{weighted
  pseudolimit} is defined by
\[\K(C,\{J, S\}_{ps}) \cong \Ps(\C,\Cat)(J, \K(C,S)).\]
I don't really want to do any examples of this one,  I want to do some
general nonsense instead.

Recall that $\Ps(\C,\Cat) = \talg$ for a 2-monad $T$ on $[\ob\C,\Cat]$,
while $\talgs = [\C,\Cat]$, so the inclusion
\[ [\C,\Cat] \to \Ps(\C,\Cat) \]
has a left adjoint $Q$, with $QJ=J'$. Thus
\[\Ps(\C,\Cat)(J,\K(C,S)) \cong [\C,\Cat](J',\K(C,S))\]
which just defines the universal property for the $J'$-weighted limit.
In other words, \emph{pseudolimits are not some more general thing,
but a special case of ordinary (weighted) limits}.  Thus we say that
a weight ``is'' a pseudolimit if it has the form $J'$ for some $J$.

\begin{remark}
This sort of phenomenon is common.
Recall, for example, that pseudo-algebras for monads
are strict algebras over a cofibrant replacement monad.  Thus talking
about things of the form \pstalg is actually \emph{less} general
than things of the form \talg, since everything of the former form
has the latter form, but not conversely.
\end{remark}

\subsection{PIE-limits again}
\label{sec:pie-limits-1}

Recall that PIE-limits are the limits constructible from products,
inserters, and equifiers. We can now make this more precise.
A weight $J:\C\to\Cat$ is a (weight for a) PIE-limit if and only
if the following conditions hold:
\begin{itemize}
\item any 2-category \K with products, inserters, and equalizers has 
      $J$-weighted-limits;
\item any 2-functor $F:\K\to\LL$ which preserves products, inserters, and
      equalizers (and for which \K has these limits) also preserves 
      $J$-weighted limits.
\end{itemize}

There is a characterization of such weights. Given a 2-functor $J:\C\to\Cat$,
first consider the underlying ordinary functor $J_0:\C_0\to\Cat_0$ obtained
by throwing away all 2-cells. Now compose this with the functor 
$\ob:\Cat_0\to\Set$ which throws away the arrows of a category, leaving
just the set of objects. This gives a functor $j:\C_0\to\Set$. Then $J$
is a PIE-weight if and only if $j$ is a coproduct of representables; and $j$
will be a coproduct of representables if and only if each connected component
of the category of elements of $j$ has an initial object.
 
Pseudolimits are also PIE-limits, as we shall now see.
For a general $T$-algebra $A$, the pseudomorphism classifier $A'$ 
was constructed from free algebras using coinserters and coequifiers.
Thus for a general weight $J:\C\op\to\Cat$ we can construct $J'$ from
``free weights'', using coinserters and coequifiers. Free weights, in this
context, are coproducts of representables, thus $J'$ can be constructed
from representables using coproducts, coinserters, and coequifiers. 
It will follow that pseudolimits can indeed be constructed using
products, inserters, and equifiers, and so that they are PIE-limits.

Now a limit weighted by a representable $\C(C,-)$ is given by evaluation
at $C$; a limit weighted by a coproduct of representables is given by 
the product of the evaluations; a limit weighted by a coinserter of 
coproducts of representables is given by an inserter of products of 
evaluations, and so on. This is part of a general result, not needed
here, that colimits of weights give iterated limits, as in the formula
$$\{J\star H,S\}\cong\{J,\{H,S\}\}$$
reminiscent of a tensor-hom situation. In other words, the 2-functor
$$\{-,S\}:[\C,\Cat]\op\to\K,$$
sending a weight $J$ to the $J$-weighted limit of $S$, sends colimits
in $[\C,\Cat]$ to limits in \K. In any case we can conclude in
the current context that $J'$-weighted limits can be constructed using
products, inserters, and equifiers, and so we have:

\begin{proposition}
Pseudolimits are PIE-limits.
\end{proposition}

The converse is false: for example inserters are not pseudolimits.
Neither are iso-comma objects, although they're pretty close (as we saw above).

Remember that $\talg$ has all PIE-limits.  It therefore has all
pseudolimits as well.  But consider the class of all limits (weights)
which are equivalent (in $[\C,\Cat]$, so that the equivalences are
2-natural) to pseudolimits.  It is not the case that \talg\ has all of
those limits.  So equivalence of limits is not always totally trivial.

For example, consider splitting of idempotent equivalences, which
seems like a very benign thing to do. If we split an idempotent 
equivalence 
\[\xymatrix{TA \ar[r]^{Te} \ar[d]_a \ar@{}[dr]|{\cong} & TA \ar[d]^a \\
  T \ar[r]_e & T}\]
in \talg, we won't necessarily get a $T$-algebra back, only a pseudo-algebra.  

As an example, let \C be a non-strict monoidal category, and $\C\st$ its
strictification. Then there is an idempotent equivalence on $\C\st$, which
when split, gives \C. This shows that \talg does not have splittings 
of idempotents when $T$ is the 2-monad for strict monoidal categories.

\subsection{Bilimits}
\label{sec:bilimits}

I'm going to write down all the same symbols, but they'll just mean
different things! So now \C and \K are bicategories, while $S:\C\to\K$ and 
$J:\C\to\Cat$ are now homomorphisms (pseudofunctors). 
The \emph{weighted bilimit} is defined by an \emph{equivalence}
\[\K(C,\{J,S\}_b) \simeq \Hom(\C,\Cat)(J,\K(C,S)).\]
Now our limits are determined only up to equivalence, instead of up to
isomorphism.

In the case when \C\ and \K\ are 2-categories and $J$ and $S$ are
2-functors, then the right hand side is equal to the right hand side
for pseudolimits, just by
definition (since $\Ps(\C,\Cat) \hookrightarrow \Hom(\C,\Cat)$ is locally an
isomorphism).  Thus \emph{every pseudolimit is a bilimit}.

On the other hand, if just \K is a 2-category, then you can replace \C 
by a 2-category $\C'$ such that homomorphisms out of \C are the same
as 2-functors out of $\C'$. Now for any $A\in\K$ we have
$$\Hom(\C\op,\Cat)(J,\K(S,A))\simeq\Ps(\C\op,\Cat)(\tilde{J},\K(\tilde{S},A))$$
where $\tilde{J}$ and $\tilde{S}$ are the 2-functors corresponding to $J$ 
and $S$. Thus {\em a 2-category with all pseudolimits also has all bilimits.}

As we shall see in the following section, the converse is false: there
are 2-categories with bilimits which do not have pseudolimits, so the
definition of pseudolimit is logically harder to {\em satisfy} than that of
bilimit. On the other hand in concrete examples it is often much easier
to {\em verify} the definition using pseudolimits than bilimits.
This is certainly the case for pseudolimits in \talg. It's also the 
case for the opposite of the 2-category of Grothendieck toposes.

\subsection{Bilimits and bicolimits in \talg}
\label{sec:talg-again}

Suppose once more that 
\K\ is locally finitely presentable and $T$ is finitary, and consider
the 2-category \talg.  It has PIE-limits, as we saw, and so has pseudo-limits, 
and so has bilimits.  So from a bicategorical perspective, we have all the
limits we might want.

\talg also has {\em bicolimits}, although not in general pseudocolimits or
PIE-colimits. Thus $\talg\op$ is an example of a 2-category with bilimits
but not pseudolimits. Just as in the case of ordinary monads, the 
(bi)colimits in \talg are not generally constructed as in \K.

The existence of bicolimits follows from:

\begin{theorem}
  Suppose we have a 2-functor $G:\talg\to\LL$ such that the composite $GJ$ in
  \[\talgs \too[J] \talg \too[G] \LL\]
  has a left adjoint $F$.  Then $JF$ is
  left biadjoint to $G$.
\end{theorem}

We start with a left 2-adjoint $F$ to $GJ$ but end up with only 
a left biadjoint to $G$. Here's the idea of the proof.
The biadjunction amounts to a (pseudonatural) equivalence
\[\talg(JFL,A) \simeq \LL(L,GA).\]
Since \talgs and \talg have the same objects, we may write $A$ as $JA$. 
Now the adjunction $F\dashv GJ$ gives an isomorphism of categories
$$\talgs(FL,A)\cong\LL(L,JGA)$$
so it suffices to show that 
$$\talgs(FL,A)\simeq\talg(JFL,JA)$$
which in turn amounts to the fact that every pseudomorphism from $FL$ 
to $A$ is isomorphic to a strict one. This will hold if we know that
$(FL)'\simeq FL$. Writing $Q$ for the left adjoint to $J:\talgs\to\talg$,
we have a pseudomorphism $p:JFL\rightsquigarrow JQJFL$ (unit of $Q\dashv J$), 
and a map $n:L\to GJFL$ (unit of $F\dashv GJ$), so we can form the composite
$$\xymatrix{ L \ar[r]^-{n} & GJFL \ar[r]^-{Gp} & GJQJFL}$$
and the corresponding {\em strict} map
$$\xymatrix{FL \ar[r]^-{r} & {QJFL} }$$
under the adjunction $F\dashv GJ$ 
provides the desired inverse-equivalence to $q:QJFL\to FL$ 
(the counit of $Q\dashv J$).

\begin{corollary}~
\begin{enumerate}[\rm ($i$)]
\item \talg\ has bicolimits;
\item for any monad morphism $f:S\to T$, the induced 2-functor
$f^*:\talg\to\salg$ has a left biadjoint.
\end{enumerate}
\end{corollary}

Part ($ii$) is easier: we have a commutative diagram of 2-functors
$$\xymatrix{
\talgs \ar[r]^{J} \ar[d]_{f^*_s} & \talg \ar[d]^{f^*} \\
\salgs \ar[r]_{J} & \salg }$$
in which the left hand map has a left adjoint, by a general
enriched-category-theoretic fact (no harder than the corresponding
fact for ordinary categories), and the bottom map has a left 
adjoint (the pseudomorphism classifier for $S$-algebras). Thus
the composite has a left adjoint, and so $f^*$ has a left 
biadjoint. (The argument as stated uses the pseudomorphism classifier 
for $S$-algebras, and so requires $S$ to have rank, but this can be avoided.)

What about part ($i$)? For any $S:\C\to\talg$, we can form the diagram 
$$\xymatrix{
\talgs \ar[r]^{J} \ar[dr]_{\talg(S,J)} & \talg \ar[d]^{\talg(S,1)} \\
& \Hom(\C\op,\Cat) }$$
and now the existence of bicolimits in \talg amounts to the existence
of left biadjoints for all such $\talg(S,1)$. So it will suffice
to show that the composite $\talg(S,J):\talgs\to\Hom(\C\op,\Cat)$ has
a left adjoint. But $\talg(S,J)\cong\talgs(QS,1)$, where $Q\dashv J$, and
$\talgs(QS,1)$ has a left adjoint provided that \talgs has pseudocolimits.
Finally since $T$ is finitary, \talgs is cocomplete (by a general 
enriched-category-theoretic fact no harder than the corresponding 
fact for ordinary categories) and so in particular has pseudocolimits.

A {\em direct} proof that \talg\ has bicolimits would be a nightmare,
but using pseudocolimits it becomes manageable.

\subsection{References to the literature}

Many people came up with some notion of weighted limit at about the
same time. But I guess the main reference for general $\mathcal V$ is now 
just \cite{Kelly-book}. On the other hand, for various limit notions
for 2-categories, \cite{Kelly-limits} is very readable. Once you've
got through that, you should turn to \cite{BKPS}, \cite{Street-Cat-limits},
and \cite{Gray}. For the beautiful
theory of PIE-limits, see \cite{PIE}. For the connection between
pullbacks and pseudopullbacks, see \cite{JS-pseudopullbacks}.
Section~\ref{sec:talg-again} is based on \cite{BKP}. For relationships
between pseudolimits and bilimits, see \cite{Power-bilimits}.

\section{Model categories, 2-categories, and 2-monads}
\label{sect:Quillen} 

This section involves Quillen model categories, henceforth called
model categories, or model structures on categories. There are various
connections between model categories, 2-categories, and 2-monads
which I'll discuss.

\begin{enumerate}[(i)]
\item Model structures \emph{on} 2-categories: a model 2-category
is a category with a model structure and an enrichment over \Cat, 
with suitable compatibility conditions between these structures. 
Any 2-category with finite limits and colimits has a 
``trivial'' such structure, in which the weak equivalences are
the categorical equivalences. These trivial model 2-categories are not
so interesting in themselves, but can be used to generate other more
interesting model 2-categories.
\item Model categories \emph{for} 2-categories: There's a model
structure on the category of 2-categories and 2-functors, and one
for bicategories too.
\item Model structures {\em induced by} 2-monads. If $T$ is a 2-monad
on a 2-category \K, we can lift the trivial model 2-category structure on 
\K coming from the 2-category structure to get a model structure on 
\talgs.
\item Model structures {\em for} 2-monads: the 2-category \Mndf of 
finitary 2-monads on \K is also a model 2-category.
\saveenumi
\end{enumerate}

One thing which I won't discuss, but deserves further study:

\begin{enumerate}[(i)] 
\recall
\item ``Many-object monoidal model categories''. There's a notion of
monoidal model category: this is a monoidal category with a model
structure, suitably compatible with the tensor product. The
many-object version of this would involve a bicategory (or 2-category)
with a model structure on each hom-category, subject to certain 
conditions (somewhat more complicated than those for monoidal model
categories).
\end{enumerate}

\subsection{Model 2-categories}

There's a model structure on the category $\Cat_0$ of categories
and functors in which the weak equivalences are the equivalences
of categories, and the fibrations are the functors $f:A\to B$
such that for any object $a\in A$ and any isomorphism $\beta:b\cong fa$ 
in $B$, there is an isomorphism $\alpha:a'\cong a$ in $A$ with 
$fa'=b$ and $f\alpha=\beta$. This is sometimes called the
``categorical model structure'' or ``folklore model structure''.  
(There are other model structures on $\Cat_0$, in particular the 
famous one due to Thomason that gives you a homotopy theory equivalent 
to simplicial sets.) 

As mentioned above, a category with a monoidal structure and a model 
structure satisfying certain compatibility conditions is called a
monoidal model category.
The cartesian product makes $\Cat_0$ into a monoidal model category. 

If we now consider a category that has both a model structure and an
enrichment over \Cat, there is a notion of compatibility between these
structures, which can be expressed in terms of the monoidal model structure
on $\Cat_0$. We call this notion a {\em model 2-category}.

First of all the 2-category \K\ is required to have finite limits and 
colimits in the 2-categorical sense; if the underlying ordinary category 
$\K_0$ already has finite limits and colimits, then the extra requirement
is that \K have tensors and cotensors with \two. A model structure on 
$\K_0$ makes \K\ into a model 2-category if two new axioms hold for any
 cofibration $i:A\to B$ and fibration $p:C\to D$:
\begin{enumerate}[(a)]
\item Given morphisms $x:A\to C$, $y:B\to C$, and $z:B\to D$, with
$px=zi$, and invertible 2-cells $\alpha:x\cong yi$ and $\beta:z\cong py$
with $p\alpha=\beta i$, there exist a morphism $y':B\to C$ and an
isomorphism $\gamma:y'\cong y$ with $p\gamma=\beta$ and $\gamma i=\alpha$;
\item If either $i$ or $p$ is trivial, then for any morphisms 
$x,y:B\to C$ and any 2-cells $\alpha:xi\to yi$ and $\beta:px\to py$ 
with $\beta i=p\alpha$, there exists a unique 2-cell $\gamma:x\to y$
with $p\gamma=\beta$ and $\gamma i=\alpha$.
\end{enumerate}
It follows that every equivalence is a weak equivalence, and that any
morphism isomorphic to a weak equivalence is itself a weak equivalence.

\subsection{Trivial model 2-categories}
\label{sec:model-structures-2}

Let \K\ be a 2-category with finite limits and colimits. The most
important limit here will be the \emph{pseudolimit of an arrow} $f:A\to B$.
Ordinarily we don't talk about limits of an arrow, since the ordinary
limit of an arrow is just its domain, but 
the pseudolimit is only equivalent to the domain, not equal.  It's 
the universal diagram
$$\xymatrix{ \arrowobject{\cong}
& A \ar[dd]^{f} \\
L \ar[ur]^{u} \ar[dr]_{v} \ar@{}[r]|{\cong\lambda} & \\
& B }$$
such that given $a:X\to A$, $b:X\to B$, and $\beta:fa\cong b$,
there is a unique $c\maps X\to L$ with
$\lambda c=\beta$ (and so also $uc=a$ and $vc=b$).  In this case, $u$
is an equivalence, because $\id:f1\cong f$ factors through by a 
$d\maps A \to L$ with $ud=1$, and one can also check that $du\cong 1$.
The technique of Section~\ref{sect:colimits} can be used to calculate
the weight for pseudolimits of arrows.

The model structure on \K is:
\begin{itemize}
\item The weak equivalences are the equivalences;
\item The fibrations are the \emph{isofibrations}, the maps such that each 
  invertible 2-cell
  
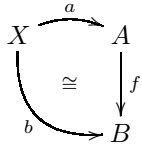
\[\xymatrix{ X \ar@/^/[r]^{a} \drlowertwocell_{b}{\omit}
               \ar@{}[dr]|(0.5){\cong} & A \ar[d]^f \\ 
  & B }\]
  lifts to an invertible 2-cell
  
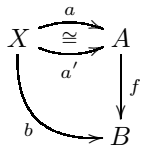
\[\xymatrix{ X \ar@/^/[r]^{a} \ar@/_/[r]_{a'} \ar@{}[r]|{\cong} 
               \drlowertwocell_b{\omit} & A \ar[d]^f \\ & B}\]
\item The cofibrations have the left lifting property with respect to the
trivial fibrations.
\item It follows that the trivial fibrations are the 
\emph{surjective equivalences}: the
$p$ for which there exists an $s$ with $ps=1$ and $sp\cong 1$.
\end{itemize}
We call such a model 2-category \K a trivial model 2-category. 
When $\K=\Cat$ this is just the folklore structure.
When \K has no non-identity 2-cells, then the equivalences are the 
isomorphisms, and all maps are isofibrations, so this agrees with the
usual notion of trivial model category. For a general 2-category \K,
however, there will be weak equivalences which are not invertible.

The pseudolimit of $f$ gives us, for any $f$, a factorization $f=vd$
where $v$ is a fibration (which follows from the universal property of
the pseudolimit) and $d$ is an equivalence.  In the case of \Cat, you
could stop there and $d$ would already be a trivial cofibration, but
in general there's more work to do, although we have reduced 
the problem to factorizing an equivalence.

The way you do that is also the way you get the other factorization:
use the dual construction. Form the pseudo\textit{co}limit of the arrow $f$,
as in the diagram below, and let $e$ be the unique map with $ei=f$, 
$ej=1$, and $e\phi=\id$.
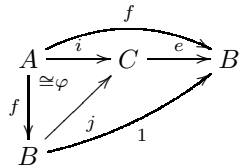
$$\xymatrix{A \ar[r]^i \ar[d]_f \ar@/^4mm/[rr]^f \ar@{}[dr]|(0.25){\cong\phi}& 
  C \ar[r]^e & B\\
  B\ar[ur]_j \ar@/_2mm/[urr]_1 & }$$
This time $i$ is a cofibration and $e$ is a trivial fibration, and if
$f$ itself is an equivalence, then $i$ has the left lifting property
with respect to the fibrations (so it's what's going to become a
trivial cofibration).

That's all I'll say about the proof.  There is, of course, a dual model
structure in which the cofibrations are characterized and the
fibrations are defined by a right lifting property.  For \Cat, these
coincide, in general they don't.

When \K\ is arbitrary, there is no reason why the model structure should be 
cofibrantly generated.  Certainly for \Cat\ it is, but even for such a
simple 2-category as $\Cat^\two$ it is not. From the homotopical point
of view the trivial model structure is trivial in several ways, including:
\begin{itemize}
\item All objects are cofibrant and fibrant;
\item The morphisms in the homotopy category $\Ho(\K_0)$ are the
  isomorphism classes of 1-cells in \K.
\end{itemize}

In the case $\K=\Cat(\EE)$, one typically considers different model
structures:
an internal functor $F\maps \aA\to \BB$ is usually called a weak equivalence 
if it's full and faithful and
essentially surjective in an internal sense.  For \Cat\ this is
equivalent to the usual notion (by the axiom of choice), but in general
it won't be.  It's the weak equivalences in this sense that people
tend to use as their weak equivalences for $\Cat(\EE)$.  When \EE\
is a topos, this was studied by Joyal and Tierney, and there's been
recent work on other cases, when \EE\ is groups (so that internal 
categories are crossed modules) or abelian groups.

\subsection{Model structures for $T$-algebras}
\label{sect:model-talgs}

Now let $T$ be a (finitary) 2-monad on (a locally finitely presentable) 
2-category \K, and \talgs the 2-category of strict algebras and strict 
morphisms. In the usual way, one can lift the model structure on \K to get one
on \talgs: a strict $T$-morphism $f:(A,a)\to(B,b)$ is a weak equivalence 
or fibration if and only if $U_sf:A\to B$ is one in \K; the cofibrations
are then defined via a left lifting property.

Now the lifted model 2-category structure on \talgs is {\em not} trivial.
In general, if $(f,\fbar):(A,a)\to (B,b)$ is
a pseudomorphism of $T$-algebras and $f:A\to B$ is an equivalence,
then any inverse-equivalence $g:B\to A$ naturally becomes an equivalence
upstairs in \talg.  This is a 2-categorical analogue of the fact that if an
algebra morphism is a bijection, its inverse also preserves the
algebra structure.
But if $f$ is strict ($\fbar$ is an identity), there is no
reason why its inverse equivalence should also be strict, and thus
no reason why $f$ should be an equivalence in \talgs.  For
example, a strict monoidal functor which is an equivalence of
categories has an inverse which is strong monoidal, but which need not be
strict.

Recall the adjoint to the inclusion
\[\xymatrix @C3pc {\talgs \rtwocell^{Q}_{J}{`\bot} &\talg}\]
where $QA=A'$, so that we have a bijection
\[\frac{A\rightsquigarrow B}{A' \to B}.\]
This fits into the model category framework very nicely.  The counit of
this adjunction
\[A' \too[q] A\] is a cofibrant replacement: a trivial fibration with
$A'$ being cofibrant.  So we see that \talg, which is the thing
we're more interested in, is starting to come out of the picture: a weak
morphism out of $A$ is {\em the same thing as} a strict morphism out
of the ``special cofibrant replacement'' $A'$ of $A$.
This is much tighter than the general philosophy that ``we
should think of maps in the homotopy category as maps out of a
cofibrant replacement.''

An algebra turns out to be cofibrant if and only if $q:A'\to A$ has a section
\emph{in} \talgs.  (There's always a section in \talg, using a pseudomorphism.)
Since $q$ is a 
trivial fibration, there will certainly be a section if $A$ is cofibrant.
Conversely, if there is a section, then $A$ is a retract of $A'$; but
$A'$ is always cofibrant, and so then $A$ must be cofibrant too.
In 2-categorical algebra, the word {\em flexible} is used in place of
cofibrant. 

\subsection{Model structures for 2-monads}
\label{sect:model-structures-2-1}

Recall now that we have adjunctions
\[\xymatrix{ \Mndf \dtwocell^W_H{'\dashv}  \\ 
  \Endf \dtwocell^V_G{'\dashv}\\
  [\ob \K_f, \K]
}\]
both of which are monadic, as is the composite.  Thus \Mndf is
both \malgs and \nalgs where $M$ is the induced monad on
\Endf and $N$ is the induced monad on $[\ob \K_f, \K]$.

Thus \Mndf has \emph{two} lifted model structures, coming from the 
trivial structures on \Endf and on $[\ob\Kf,\K]$.  They're not
the same, since something can be an equivalence all the way downstairs
without being one in \Endf (which is itself the 2-category of algebras for
another induced monad on $[\ob \K_f, \K]$). 

A monad map $S\too[f] T$ is a 2-natural transformation compatible with
the unit and multiplication.  If the 2-natural transformation is an
equivalence in \Endf, it is a weak equivalence for the
$M$-model structure; if the \emph{components} of the 2-natural
transformation are equivalences, it is a weak equivalence for
the $N$-model structure.

It's the $M$-model structure (the one lifted from \Endf) which seems to 
be more important, and we'll only consider that one here.
The corresponding prime construction classifies
pseudomorphisms of monads.  These are precisely the things that arise
when talking about pseudoalgebras: recall that a pseudo-$T$-algebra
was an object $A$ with a pseudo-morphism
\[T\rightsquigarrow \langle A,A\rangle \] into the ``endomorphism 2-monad'' of
$A$, corresponding to maps $TA\too[a] A$ which are associative and
unital up to coherent isomorphism.

This corresponds to a strict map $T'\to \langle A,A\rangle$, so that
$\ttalg = \pstalg$.  (This is the part of the
justification for working with strict algebras that people tend to
understand first, but it's the less important one: see 
Remark~\ref{rmk:pseudo-algebras} above.)

If $q:T' \to T$ has a section in $\malgs=\Mndf$, then $T$ is said to be 
\emph{flexible} (=\emph{cofibrant}).  This was the context in which the 
notion of flexibility was first introduced.  {\em Any monad that you can 
give a presentation for without having to use equations between objects
  is always flexible.}  For example, the monad for monoidal categories
is flexible, but the monad for strict monoidal categories, which 
involves the equation $(X\ot Y)\ot Z=X\ot(Y\ot Z)$, is not.

Flexible monads have the property that every pseudo-algebra is (not
just equivalent but) {\em isomorphic} to a strict one; in fact isomorphic
via a pseudomorphism whose underlying \K-morphism is an identity!  
Remember that the importance of
pseudo-algebras is \emph{not} for describing concrete things, but for
the theoretical side, since various constructions don't preserve
strictness of algebras.  For particular structures like monoidal
categories, you're better off choosing the ``right'' monad to start with:
the one for which monoidal categories are the strict algebras.

\subsection{Model structure on 2-Cat}
\label{sect:qmc2cat}

\twocat is the category of 2-categories and 2-functors.  It
underlies a 3-category, and a 2-category, and perhaps more importantly
a Gray-category.  But we want to describe a model structure on the mere
category \twocat, analogous to the one above for \Cat.

The weak equivalences will be the \emph{biequivalences}. 
Recall that $F\maps \A\to\B$ is a
biequivalence if
\begin{itemize}
\item each $F:\A(A,B)\to \B(FA,FB)$ is an equivalence of
  categories; and
\item $F$ is ``biessentially surjective'' on objects: if
  $C\in\B$, there exists an $A\in \A$ with $FA\simeq C$ in \B.
\end{itemize}

Every equivalence has an inverse-equivalence, going back the other way.
For a biequivalence $F:\A\to\B$ you can build a thing
$G\maps \B\to \A$ with $GF\simeq 1$ and $FG\simeq 1$.  You can make
$G$ a pseudofunctor, but generally not a 2-functor, even when $F$ is one.
That's somehow the whole point of
the model structure. Similarly the equivalences $FG\simeq 1$ and $GF\simeq 1$
will generally only be pseudonatural.

Clearly biequivalence is the right notion of ``sameness'' for
bicategories, or 2-categories, but there is this stability (under 
biequivalence-inverses)  problem,
if you want to work entirely within \twocat. If you allow pseudofunctors,
and so move to \twocatps, then as we have seen, you lose completeness
and cocompleteness.

The fibrations are similar to the case of categories. Fibrations for
the model structure on \Cat involved lifting invertible 2-cells;
here we lift equivalences: a 2-functor $F\maps\A\to\B$ is a fibration if
\begin{itemize}
\item given an object $A$ upstairs and equivalence downstairs, we have 
a lift as in 
  \[\xymatrix{ A' \ar@{|.>}[d] \ar@{.>}[r]^\simeq & A \ar@{|->}[d] & \A \ar[d]^F \\
    B \ar[r]^\simeq & FA & \B}\]
\item given a 1-cell $f$ upstairs and an invertible 2-cell downstairs, we have 
a lift as in 
  \[\xymatrix{ A' \ar@{|->}[d] \ruppertwocell<2>{\omit\cong} \ar@/_2mm/@{.>}[r]
  & A \ar@{|->}[d] & \A \ar[d]^F \\
    B \rtwocell{\omit\cong} & FA & \B}\]
Equivalently, each of the functors 
 $\A(A_1,A_2)\to\B(FA_1,FA_2)$ is an (iso)fi\-bra\-tion in \Cat.
\end{itemize}

Note that the notion of biequivalence is \emph{not} internal to the
3-category or Gray-category of 2-categories, 2-functors, and so on, 
which speaks against the existence of a general model structure on 
an arbitrary 3-category or Gray-category which would reduce 
to this one.

There's an equivalent way of characterizing the fibrations which is useful.  
Keep the iso-2-cell lifting property as is, but modify the 
equivalence-lifting to deal with adjoint equivalences rather than
equivalences. Here it is not just the 1-cell, but also the equivalence-inverse,
and the invertible unit and counit which must be lifted. 

In the presence of the iso-2-cell lifting property, these two types
of equivalence-liftings are equivalent:
clearly the lifting of adjoint equivalences implies the lifting of 
equivalences, since we can complete
any equivalence to an adjoint equivalence, but the converse is also
true provided that we can lift 2-cell isomorphisms.

This is related to a mistake I made in my first paper on this topic,
where I used a condition like this on lifting equivalences that aren't
necessarily adjoint equivalences.  Regard ``being an equivalence'' as a
property, and ``an adjoint equivalence'' as a structure, but be wary of
regarding ``a not-necessarily-adjoint equivalence'' as a structure.
Adjoint equivalences are now completely algebraic, classified by maps
out of ``the free-living adjoint equivalence'', which is biequivalent to
the terminal 2-category 1.  A ``free-living not-necessarily-adjoint 
equivalence'' would \emph{not} be biequivalent to 1.

The trivial fibrations, which are the things which are both fibrations
and weak equivalences, can be characterized as the 2-functors that
\begin{itemize}
\item are surjective on objects; and
\item have each $\A(A_1,A_2)\to \B(FA_1,FA_2)$ a surjective equivalence
  (a trivial fibration in \Cat).
\end{itemize}
Note that the trivial fibrations don't use the 2-category structure;
you don't need anything about the composition to know what these
things are, only the ``2-graph'' structure.  So they're much simpler to
work with.

There's an $\omega$-categorical analogue to these things,
which permeates Makkai's work on $\omega$-categories.  You don't
need the $\omega$-category structure, only a globular set, to say what
this means.  The corresponding notion of ``cofibrant object'' is then
what he calls a ``computad''.

It's a bit less trivial than with the other model structures to prove 
that this all works, but it's not really hard.  Everything is directly a 
lifting property (once you use the version with adjoint equivalences), so
finding generating cofibrations and trivial cofibrations is easy.

All objects are fibrant, but \emph{not} all objects are
cofibrant.  We have a ``special'' cofibrant replacement $q:\A'\to \A$
with the property that pseudofunctors out of \A are the same
as 2-functors out of \A':
\[\frac{\A\rightsquigarrow \B}{\A'\to \B}\]
and $\A$ is cofibrant
(flexible) if and only if the trivial fibration $q$ has a section in \twocat.
This happens exactly when the underlying category $\A_0$ is free on
some graph (you haven't imposed any equations on 1-cells, but you may 
have introduced isomorphisms between them). In principle, a cofibrant
$\A_0$ could be a retract of something free,
but it turns out that this already implies that it is free.

There are three main things of interest to me in relation to the 
model structure on \twocat.  The first is the equation  
``cofibrant = flexible''. The second involves the
monoidal structures.  The
model structure is \emph{not} compatible with the cartesian product
$\times$.  The thing to have in mind is that the locally
discrete 2-category $\two=(0\to 1)$ is cofibrant, but $\two\times\two$ is
not, since the commutative square
\[\xymatrix{(0,0) \ar[r]\ar[d] & (0,1) \ar[d]\\
(1,0)  \ar[r] & (1,1) }\]
is not free.  There are various tensor products you can put on
\twocat.  The cartesian product is also called the \emph{ordinary
  product} (since it is also a special case of the tensor product of
\V-categories), but I like to call it the \emph{black product} since
the square is ``filled in'', in the sense that the square commutes.
(Think of nerves of categories: 2-simplices, or solid triangles, 
represent  commutative triangles in a category, whereas the boundary
of a 2-simplex represents a not necessarily commutative triangle.)

There's also the \emph{white} or \emph{funny} product, in which the
square does not commute (think of it as just the boundary, not filled
in). It's a theorem that on \Cat\ there
are exactly 2 symmetric monoidal closed structures: the ordinary one
and the ``funny'' one.  The closed structure corresponding to the funny
product is the not-necessarily-natural transformations (just
components).  Enriching over this structure gives you a
``sesquicategory'' (perhaps an unfortunate name, but you can see how it
came about), which has hom-categories and whiskering, but no
middle-four interchange, hence no well-defined horizontal composition
of 2-cells. This funny tensor product can also be defined on \twocat,
or indeed on \VCat for any \V.

In the case of \twocat, there's an intermediate possibility:
the \emph{Gray} or \emph{grey} tensor product, due to John Gray, in
which you put an isomorphism in the square, so it's ``partially filled
in''.
\[\xymatrix{(0,0) \ar[r]\ar[d] \ar@{}[dr]|{\cong} & (0,1) \ar[d]\\
(1,0)  \ar[r] & (1,1)}\]
(This is the ``pseudo'' version of the Gray tensor product; there's
also a ``lax'' version: different shade of grey!) 

The black and white tensor product 
make sense for any \V\ at all, but the grey one doesn't.
%The reason I started thinking about this at all was to think about the
%Gray tensor product from various different points of view.  
There's a canonical comparison from the funny/white product to
the ordinary/black one, and the Gray/grey tensor product is a sort of 
``cofibrant 
replacement'' in between.  

The Gray tensor product $\two\ot\two$ is cofibrant, and more
generally, the model structure is compatible with the Gray tensor
product.

The third thing of interest is the connection between 2-categories and
bicategories. There is a model structure on the category \bicat of
bicategories and strict homomorphisms. The notion of biequivalence
still makes perfectly good sense, and these are the weak equivalences.
The fibrations are once again the maps which are isofibrations on the
hom-categories, and have the equivalence lifting property. Thus the full
inclusion $\twocat \hookrightarrow \bicat$ preserves and reflects weak 
equivalences and fibrations.  This inclusion has a left adjoint 
``free strictification'' \[\xymatrix{\twocat \rtwocell{`\bot} & \bicat}\]
given by universally making the associativity and identity isomorphisms
into identities.
This left adjoint is not a pseudomorphism classifier, since we are using
strict morphisms of bicategories and of 2-categories,
and it's not the usual strictification functor ``st'' either: in general the 
unit is \emph{not} a biequivalence. But the component of the unit at a
{\em cofibrant} bicategory {\em is} an equivalence. This fits well into
the model category picture: it's part of what makes this adjunction a
Quillen equivalence (one that induces an equivalence of homotopy 
categories). In fact the usual strictification can be seen as a derived 
version of the free strictification.

There exist bicategories (even monoidal categories) for
which there does not exist a strict map into a 2-category which is a
biequivalence, although we know that any \B\ has a pseudofunctor $\B
\rightsquigarrow \B_{st}$ which is a biequivalence.  The point about 
cofibrant bicategories \B
is that ``there aren't any equations between 1-cells'', and this is
what makes the unit at such a \B an equivalence.

Just as for \twocat, we have pseudomorphism classifiers in \bicat, 
which serve as ``special cofibrant replacements''.

The model structure on \twocat\ is proper: showing that it's left
proper (biequivalences are stable under pushout along cofibrations)
is harder than any of the other results mentioned above, but the
fact that it is right proper
is an immediate consequence of the fact that every object is fibrant.

\subsection{Back to 2-monads}
\label{sec:back-2-monads}

There's a connection between the model structure on \Mndf and that
on \twocat. There's a 2-functor
\[\sem\maps\Mndfcop \to \Twocat/\Cat\]
which you might call \emph{semantics}, defined by:
\[ T \mapsto (\talg \too[U] \Cat).\]
and
\[ (S\too[f] T) \mapsto (\talg \too[f^*] \salg)\]
since if $a:TA\to A$ is a $T$-algebra, then its composite with $fA:SA\to TA$ 
makes $A$ into an $S$-algebra.

In the ordinary unenriched case or the \V-enriched case, or even here, if we 
used $T\mapsto \talgs$ rather than $T\mapsto \talg$, the semantics
functor would be fully faithful. But the semantics functor defined
above, using \talg, is not: to give a map $\sem(T)\to\sem(S)$ in 
$\Twocat/\Cat$ corresponds to giving a weak morphism from $S$ to $T$,
but not in the sense of pseudomorphisms of monads, considered above;
rather in a still broader sense, in which the
$fA:SA\to TA$ need not even be natural.

Now the definitions of fibration, weak equivalence, and trivial fibration
in \twocat have nothing to do with smallness, and make perfectly good
sense in the category \Twocat of not-necessarily-small 2-categories.
We can therefore define a morphism in $\Twocat/\Cat$ to be a fibration,
weak equivalence, or trivial fibration if the underlying 2-functor
in $\Twocat$ is one.

Under these definitions,  \sem preserves limits, fibrations, 
and trivial fibrations, as one verifies using the 2-monads 
$\langle A,A\rangle$, $\{f,f\}_\ell$, and so on.  Limits, fibrations, and trivial fibrations in
\Mndfcop, correspond to colimits, cofibrations, and
trivial cofibrations in \Mndfc.  Thus, it should in principle be the right
adjoint part of a Quillen adjunction.  It's not, of course, because of
size problems: $\Twocat/\Cat$ has large hom-categories, and \sem lacks
a left adjoint.

The assertion that \sem preserves the weak equivalence $q:T'\to T$ is 
equivalent to the assertion that every pseudo $T$-algebra is equivalent
to a strict one. More generally, \sem preserves all weak equivalences 
if and only if pseudo algebras are equivalent to strict ones for every 
$T$. Whether or not this is the case is an open problem in the current 
generality, but it is true that \sem preserves weak equivalences between 
cofibrant objects (flexible monads).

\subsection{References to the literature}

Model categories go back to \cite{qmc}. There are now several modern
treatments: \cite{Hovey-book} is one which emphasizes the compatibility
between model structures and monoidal structures.
The ``categorical'' model structure on \Cat seems to be folklore; the
first reference I know is \cite{Joyal-Tierney-CatE}. The Thomason
model structure on \Cat comes from \cite{Thomason-Cat}, which does
have an error in the proof of properness, corrected in
\cite{Cisinski-correction}. For the model structures on \textrm{2-Cat}
and \textrm{Bicat} see \cite{qm2cat,qmbicat}. There is also a 
``Thomason-style'' model structure on \textrm{2-Cat} \cite{Thomason-2-Cat}.
The theory of model 2-categories, the model structure on the category of 
monads, and its relation to structure
and semantics, all come from  \cite{hty2mnd}.

\section{The formal theory of monads}
\label{sect:ftm}

In this section we return to formal category theory; in fact, to
one of its high points: the formal theory of monads.

\subsection{Generalized algebras}
\label{sect:generalized-algebras}

Let's start by thinking about ordinary monads.
Let $A$ be a category, $t = (t,\mu,\eta)$ a monad on $A$. Write
$A^t$ for the Eilenberg-Moore category (the category of algebras). The
starting point is to think about the universal property of this
construction.  What is it to give a functor $C\too[a] A^t$?  We give
an algebra $ac$ for each $c\in C$, and use $ac$ also for the name of
the underlying object, with structure map $tac\too[\alpha c] ac$.  And
for every $\gamma\maps c\to d$, we have an $a\gamma\maps ac\to ad$
with a commutative square
\[\xymatrix{tac \ar[r]^{\alpha c}\ar[d]_{ta\gamma} & ac \ar[d]^{a\gamma}\\
  tad\ar[r]_{\alpha d} & ad.}\]
This square awfully like a naturality square; it wants to
say that $\alpha$ is natural with respect to $\gamma$, and indeed this
is in fact the case.
What we're actually doing is giving a functor $C\too[a] A$ and a
natural transformation
\[\xymatrix{ C \ar[r]^a \drlowertwocell<-4>_a{<-1>\alpha} & A \ar[d]^t\\
  & A}\]
with equations of natural transformations
\[\xymatrix{t^2a \ar[r]^{\mu a}\ar[d]_{t\alpha} & ta \ar[d]^\alpha & a
  \ar[l]_{\eta a} \ar[dl]^1\\
  ta\ar[r]_\alpha & a}\]
which just says that on components, it makes each $ac$ into a
$t$-algebra.

You might call this a \emph{generalized algebra}, or a
\emph{$t$-algebra with domain $C$}.  Think of a usual algebra as an
algebra with domain $1$.

Similarly, you can look at natural transformations.  To give a natural
transformation
\[\xymatrix{ C\rtwocell^a_b & A^t}\]
amounts to giving
\[\xymatrix{ C\rtwocell^a_b{\phi} & A}\]
which is suitably compatible, in the sense that
\[\xymatrix{ta \ar[r]^{t\phi}\ar[d]_\alpha & tb \ar[d]^\beta\\
  a \ar[r]_\phi & b.}\]
This is the universal property of the Eilenberg-Moore construction,
and the starting point of the theory.

I've been talking all along about categories, but once we've moved beyond
algebras with domain 1, there's no
reason to restrict in that way, so we can talk instead about a monad on
an object $A$ in any 2-category \K.  (The notion of monad has not been 
weakened in any way.  The 2-category \K\ 
might be \Cat, or \twocat, or \VCat, but we use the same definition.)

We can't just {\em construct} $A^t$ as we
did before, but we can \emph{ask} whether there exists an object $A^t$
with the universal property. A slick way to do this is as follows. The
hom-category $\K(C,A)$ has a monad $\K(C,t)$ on it (since 2-functors
take monads to monads), and this is the ordinary type of monad in
\Cat.  The endofunctor part of this monad sends $a\maps C\to A$ to
$ta\maps C\to A$.  This generalized notion of algebra is then nothing
but the usual sort of algebra for the ordinary monad $\K(C,t)$.  So
what we want is an isomorphism
\[\K(C,A^t) \cong \K(C,A)^{\K(C,t)}\]
naturally in $C$ (where the right hand side means the ordinary Eilenberg-Moore
category of algebras for the ordinary monad $\K(C,t)$).
We call $A^t$ the Eilenberg-Moore object of $t$, or EM-object for short.

It turns out that in some 2-categories, such as \Cat, it's enough to check
the universal property for $C=1$, since $1$ generates \Cat, in a 
suitable sense; but in an abstract 2-category there may
not be a 1, and if there is one, it may not be enough to get the full
universal property. Of course there are similar phenomena in ordinary
category theory.

The universal property of $A^t$ makes it look like a limit, and indeed
it is one, but we'll look at some other points of view first.

\subsection{Monads in \K}
\label{sec:monads-k}

Let \K be a 2-category. Previously we looked at the 2-category \Mndf
of finitary 2-monads {\em on} \K\ (as a fixed base object). We now
consider the 2-category $\mnd(\K)$ of all the (internal) monads {\em in}
\K, with variable base object.
\begin{itemize}
\item Its objects are monads in \K.
\item Its 1-cells correspond to morphisms which lift
  to the level of algebras:
  $$\xymatrix{A^t \ar[r]^{\bar{m}}\ar[d]_{u^t} & B^s \ar[d]^{u^s}\\
    A \ar[r]_m & B}$$
  where we assume temporarily that $A^t$ and $B^s$ exist,
  and we can think of this commutative diagram as an identity 2-cell and 
  then take its mate, since the $u$'s are right adjoints:
  \[\xymatrix{A \ar[r]^m \ar[d]_{f^t} \drtwocell\omit & B \ar[d]^{f^s}\\
    A^t \ar[r]^{\bar{m}}\ar[d]_{u^t} & B^t \ar[d]^{u^s}\\
    A \ar[r]_m & B}\]
  which we then paste together to get a 2-cell, and the forgetful-free
  composite gives us the monads.  Thus we should define a morphism of
  monads to be a morphism $m:A\to B$ with a 2-cell $\phi:sm\to mt$ such that 
  the diagrams
  $$\xymatrix{
ssm \ar[r]^{s\phi}\ar[d]_{\mu m} & smt \ar[r]^{\phi t} & mtt  \ar[d]^{m\mu} &&
m \ar[dr]^{m\eta} \ar[dl]_{\eta m} \\
sm\ar[rr]_\phi && mt & sm \ar[rr]_{\lambda} && mt}$$
  commute. We could {\em define} a morphism of monads to be a morphism 
$m:A\to B$ with a lifting $\bar{m}:A^t\to B^s$ as above,
except that the EM-objects need not exist in a general 2-category; but
when they do exist, the two descriptions are equivalent.
\item The 2-cells in $\mnd(\K)$ are 2-cells 
  \[\xymatrix{ A \rtwocell^m_n{\rho} & B}\]
  in \K with a compatibility condition, which you could express as saying
  that $\rho$ lifts to a $\overline{\rho}$ between EM-objects, or
  you could equivalently express as saying that
  \[\xymatrix{sm \ar[r]^\phi\ar[d]_{s\rho} & mt \ar[d]^{\rho t}\\
    sn\ar[r]^\psi & nt}\]
  commutes.
\end{itemize}

There's a full embedding $\id\maps \K\hookrightarrow \mnd(\K)$ sending $A$ to
the identity monad $(A,1)$ on $A$, and doing the obvious thing on 1-cells
and 2-cells.  This is particularly clear in the EM-objects
picture, since if $t=1$ then $A^t=A$, so obviously $m$ will lift uniquely to 
an $\bar{m}$, which is what fully faithfulness of \id\ says.

A trivial observation is that for any monad we can always choose to
forget the monad and be left with the object, and this is left adjoint
to \id.  The more interesting thing, however, is the existence of a 
right adjoint: this amounts exactly to a choice of an EM-object for 
each monad in \K.  Why?  Look at the universal property:
If $A\mapsto A^t$ is the right adjoint, this says that
\[\K(C,A^t) \cong \mnd(\K)((C,1), (A,t))\]
The key point is that the right hand side is equal to $\K(C,A)^{\K(C,t)}$, 
since an object (that is, a morphism $(C,1)\to(A,t)$) involves $a$ and
$\alpha$ as in 
\[\xymatrix{C \ar[r]^a\ar[d]_1 \drtwocell\omit{\alpha} & A \ar[d]^t\\
  C\ar[r]_a & A}\]
subject to exactly the conditions which make $(a,\alpha)$ into a generalized
algebra.

Now the really beautiful thing happens: we can start looking at duals
of \K\ and see what happens.  Consider first $\K\co$, where we reverse
the 2-cells but not the 1-cells.  A monad in $\K\co$ is then a
\emph{comonad} in \K.  And an EM-object in $\K\co$ is the obvious
analogue for comonads.  If $\K=\Cat$, we get ordinary comonads, and
the EM-object is the usual category of coalgebras for the comonad.

That's nice, but not incredibly surprising.  What's more interesting
is what happens in $\K\op$. When we consider monads not in \K but in
$\K\op$ we have to reverse the direction of the 1-cells, as in 
$$\xymatrix @R1pc {
& A \ar[dl]_{t} \\
A && A \lluppertwocell^{t}{^<1>\mu} \ar[ul]_{t} &&
A && A \lltwocell^{t}_{1}{^\eta}
}$$
but this is nothing but a monad in \K!

But what about the EM-object?  The arrows are reversed, so we get a
different universal property.  An algebra for this monad consists of
\[\xymatrix{ C & A \ar[l]_a \\
  & A \ar[u]_{t} \uluppertwocell<4>^a{^<1>\alpha} }\]
The wonderful thing is that in the case $\K=\Cat$ this is the same thing as 
a map $A_t \to C$ where $A_t$ is the Kleisli category. Recall that 
the Kleisli category can be defined as the full subcategory of the
Eilenberg-Moore category consisting of the {\em free} $T$-algebras (the
algebras of the form $(tx,\mu x)$ for some $x\in A$), or equivalently
as the category with the same objects as $A$, but with morphisms 
from $x$ to $y$ given by the morphisms in $A$ from $x$ to $ty$ (the 
monad structure is then used to make this into a category). The latter
description is more convenient here: given $(a:A\to C,\alpha:at\to a)$
as above, the induced functor $A_t\to C$ sends an object $x\in A$ to 
$ax$, and a morphism $\phi:x\to ty$ to the composite 
$$\xymatrix{ax \ar[r]^{a\phi} & aty \ar[r]^{\alpha y} & ay}$$
in $C$.

For a general 2-category \K, the 
Eilenberg-Moore objects in $\K\co$ are called Kleisli objects (in \K).
It's true in any 2-category that
the EM-object is the terminal adjunction giving rise to the monad, and
the Kleisli object is the initial one, but the universal property given
above is richer in that it refers to maps with arbitrary domains.

Using $\K\coop$, of course, gives you Kleisli objects for comonads.

\subsection{The monad structure of \mnd}
\label{sec:mnd-as-monad}

Now, where does the construction $\mnd(\K)$ really live?  Consider
the category \twocat of 2-categories and 2-functors.  
Completely banish from your mind all
concerns about size, which doesn't have any role here.  So far we've
constructed a 2-category $\mnd(\K)$ for any 2-category \K.
and this is clearly completely functorial, so we get a functor
\[\mnd\maps \twocat\to\twocat\]
and the inclusion \id\ is clearly natural in \K, so we get a natural 
transformation
\[\xymatrix{\twocat \rtwocell^{1}_{\mnd}{\id} & \twocat}\]
A certain sort of person is tempted to wonder whether this is part of 
the structure of a monad on 2-Cat!  We do have a composition map
$$\xymatrix @C3pc 
{{}\twocat \rrtwocell^{\mnd^2}_{\mnd}{\quad\comp} && {}\twocat}$$
and what it does is one of the most striking aspects of the formal
theory of monads. 

This composition map sends a monad in $\mnd(\K)$ to a monad in \K.
What is a monad in $\mnd(\K)$?  It consists of 
\begin{itemize}
\item a monad $(A,t)$ in \K (an object of $\mnd(\K)$)
\item an endo-1-cell, which consists of a morphism $A\too[s] A$ in \K\
  with a 2-cell $ts\too[\lambda] st$ (with conditions)
\item A multiplication $(s,\lambda)(s,\lambda) \to (s,\lambda)$,
  corresponding to $s^2\too[\nu] s$ (with conditions)
\item a unit $1\to (s,\lambda)$ corresponding to $1\to s$ (with
  conditions)
\end{itemize}
As well as the conditions for these to be 1-cells and 2-cells in
$\mnd(\K)$, we need the conditions for this to be a monad there.
These make $s$ itself into a monad on $A$ in \K.  The 2-cell $\lambda$
is now what's called a \emph{distributive law} between these two
monads, which is exactly what you need to ``compose'' these two monads and
get another monad.

Think about this as being like the tensor product of rings.  $R\ot S$ is
the tensor product of the underlying abelian groups, with
multiplication
$$\xymatrix @C3pc { R\ot S\ot R\ot S \ar[r]^-{1\ot\mathrm{tw}\ot1} &
R\ot R\ot S\ot S \ar[r]^-{m_R\ot m_S} & R\ot S. }$$ 
The point is we're trying to do something
very similar, but here we're in a world where the tensor product is
not commutative, so we don't have the twist.  So $\lambda$ plays the
role of the twist; it's a ``local'' commutativity that
only applies to these two objects.  The conditions put on it are
exactly what we need to make the composite $st$ into a monad.

For example, the multiplication on $st$ is then
\[stst \too[s\lambda t] sstt \too[ss\mu] sst \too[\mu t] st\] The
notion of distributive law, in the ordinary case of categories, is due
to Jon Beck, and he proved that we have a bijection between
distributive laws $ts\to st$ and ``compatible'' monad structures on
$st$, and also to liftings of $s$ to $A^t$ (whenever $A^t$ exists).
It's not as well known as it should be and is frequently
rediscovered.

You can also do this for $\K\co$ or $\K\op$ or $\K\coop$, of course. 
A distributive law in $\K\op$ is formally the same as a distributive
law in \K, but now rather than liftings of $s$ to $A^t$, it gives you
extensions of $t$ to $A_s$ along the left adjoint $f_s:A\to A_s$.

\begin{remark} Operads are monoids in a monoidal category, so there is
a corresponding notion of distributive law between operads. Furthermore,
the passage from the monoidal category of collections to the monoidal 
category of endofunctors is strong monoidal, so distributive laws
between operads induce distributive laws between the induced monads,
and this process is compatible with the formation of the composite
operad/monad. Just as not every monad arises from an operad, not 
every distributive law between monads arises from a distributive 
law between operads, even when the monads themselves do arise from
operads.
%\textbf{Peter: the notion of action of one operad on another, as
%  formulated, involves diagonals, so it only makes sense cartesian?
%  Algebraic distributivity $x(y+z) = xy + xz$ involves the diagonal on
%  $x$.}
\end{remark}

\begin{example}
Groups are particular monoids in \Set, so there is a corresponding
notion of distributive law. If a group $G$ acts on a group $H$, then
there is a distributive law $G\t H\to H\t G$ sending $(g,h)$ to $(g.h,g)$,
and the induced ``composite'' is the semidirect product $H\sd G$.
This generalizes to arbitrary monoids in a cartesian monoidal category.
\end{example}

\subsection{Eilenberg-Moore objects as limits}
\label{sec:making-it-into}

There are two ways to see Eilenberg-Moore objects as weighted 
limits.  Remember that way back in Section~\ref{sect:lax-functors}, 
we saw that monads $t$ in \K\ correspond to lax functors
$\tilde{t}\maps 1\to \K$.  Then the \emph{lax limit} of $\tilde{t}$ is exactly
the EM-object $A^t$.

I haven't explicitly discussed lax limits of lax functors, but
it's not hard to extend the definition of lax limit to cover this case.
Alternatively one can
replace the lax functor by the corresponding 2-functor out of the
``lax morphism classifier'', and then just take the lax limit of the
2-functor. Let's see how this would work. 

First recall how $\tilde{t}$ is defined. It sends $*$ to $A$, and $1_*$ to an
endomorphism $t\maps A\to A$, the unit is the lax unit comparison, and
the multiplication is the lax composition comparison.  To understand
the lax limit of these sorts of things, we should think about lax
cones.  A lax cone would involve a vertex $C$ of \K, with just one
component $C\too[a] A$, and a lax naturality 2-cell for every 1-cell
in $1$:
\[\xymatrix{C \ar[r]^a\ar@{=}[d] \drtwocell\omit{\alpha} & A \ar[d]^t\\
  C\ar[r]_a & A}\]
and some conditions.

The lax morphism classifier on 1 is a 2-category \mnd with a bijection
\[\frac{\text{2-functors } \mathbf{mnd} \too \K}{\text{lax functors }
  1\too \K}\]
but such lax functors are in turn the same as monads in \K. Thus \mnd is the 
universal 2-category containing a monad.
Remember that a monad in \K\ is the same as a monoid in a
hom-category, and we know the universal monoidal category containing a
monoid is the ``algebraic $\DD$'', the category \ordf of (possibly empty)
finite ordinals. This is not the $\DD$ of simplicial sets:
an extra object has been added.  Thus $\mathbf{mnd}$ has one object $*$ and
$\mathbf{mnd}(*,*) = \ordf$. 

Now we have a limit notion ($\tilde{t}\mapsto A^t$), and we want to 
know the corresponding weight $J:\mnd\to\Cat$, so that $\{J,\tilde{t}\}=A^t$.
We saw in Section~\ref{sect:colimits} that the recipe for calculating $J$ is 
to consider the Yoneda functor $\mnd\to[\mnd,\Cat]\op$ and form the limit of 
it, or equivalently the colimit of $\mnd\op\to[\mnd,\Cat]$. The colimit is the 
Kleisli object; since we are in a presheaf 2-category $[\mnd,\Cat]$ it is 
computed pointwise. The weight is called \alg; it's now a straightforward 
exercise to calculate it.

Of course, in general, \alg-weighted limits may or may not exist. 
Subject to the existence of the relevant limits, they can be built up
from other limits we already know: 
\begin{itemize}
\item First form the inserter of $\xymatrix{A \rtwocell^t_1{\omit} &
    A}$.  This is an $A_1\too[k] A$ equipped with a 2-cell $tk\too[\kappa] k$.
\item Then take the equifier of $k (\eta k)$ and $1$ to get an
  $A_2\too[k'] A_1$ such that the identity law holds.
\item Finally take the equifier of something else to get the
  associativity.
\end{itemize}
In particular, this shows that EM-objects are PIE-limits, in fact
\emph{finite} PIE-limits.

\subsection{Limits in \talgl and \talgc}
\label{sec:limits-talg_c-talg}

\talgl and \talgc are, recall, the 2-categories of strict $T$-algebras 
with lax and with colax morphisms.
Recall also that we had nice pseudo-limits in \talg; here it's much
harder.

In \talgl, you have oplax limits, and in \talgc you have lax
limits (it's a twisted world we live in!) These are much more 
restricted classes of limits, not including inserters, equifiers,
comma objects and many of our other favourite limits.

I described how to construct inserters and equifiers in \talg: form
the limit downstairs and show that the thing you get canonically
becomes an algebra. This involves morphisms $(f,\bar{f})$ and $(g,\bar{g})$
between $T$-algebras $(A,a)$ and $(B,b)$, and 2-cells $fk\to gk$ for some $k$.
If you look carefully at the construction, you'll see that $\bar{f}$ needs
to be invertible, but $\bar{g}$ can be arbitrary. So you can form
inserters and equifiers in \talgl provided that one of the 1-cells (the
one that is, or tries to be, the domain of the 2-cells) is actually pseudo. 

Dually, in \talgc, it's the other 1-cell which needs to be pseudo.
Now the Eilenberg-Moore object of a monad $(A,t)$ can be constructed
using the inserter $k:C\to A$ of $t$ and $1_A$, and then an equifier
(see Section~\ref{sec:making-it-into}).
Furthermore $1_A$ will always be strict, and it turns out that \talgc
does have Eilenberg-Moore objects for monads. The most important 
case is where $T$-algebras are monoidal categories and so \talgc
has opmonoidal functors. Then a monad in \talgs is an ordinary monad
for which the category is monoidal, the endofunctor opmonoidal, and
the natural transformations are opmonoidal natural transformations;
this is sometimes called a Hopf monad.

\subsection{The limit-completion approach}
\label{sec:ftm-2}

We can now see EM-objects as weighted limits in the strict sense, and
there's a well-developed theory of free completions under classes of
weighted limits.  So we can form the free completion $\EM(\K)$ of a
2-category \K\ under EM-objects; or we can form the corresponding
colimit completion $\KL(\K)$, which freely
adds Kleisli objects. These are related: $\EM(\K)=\KL(\K\op)\op$.

The colimit side is more familiar to construct.  To freely add all
colimits to an ordinary category, we take the presheaf category; to
add a restricted class, we take the closure in the presheaf category
under the colimits we want to add.  So here, to get $\KL(\K)$, we take
the closure of the representables in $[\K\op,\Cat]$ under Kleisli
objects.  It's part of a general theorem that this works, at least
when \K\ is small.  

Sometimes it can be tricky to calculate exactly which things
appear in this completion process.  You start with the representables
and throw in the relevant colimits of representables. There will now
be new diagrams, and we may have to add colimits for these. This
can continue transfinitely. The nice thing
about the case of Kleisli objects is that, as we shall see, it stops after one 
step.  

Colimits in the functor category are constructed pointwise, so we 
construct Kleisli objects as in \Cat.  The key facts are:
\begin{itemize}
\item A left adjoint in \Cat\ is of Kleisli type if and only if it is
bijective on objects;
\item These are closed under composition.
\end{itemize}
Now, given a monad $t$ on $A$ we throw in the Kleisli object $A_t$ in
$[\K\op,\Cat]$, which may have a new monad $s$ on it.  We then throw
in its Kleisli object for $s$ to get $(A_t)_s$, but then the composite
\[A \too A_t \too (A_t)_s\]
is also a bijective-on-objects left adjoint, hence $(A_t)_s$ is also a
Kleisli object for a monad on $A$.  Thus this is a 1-step process.

Therefore, we can identify (up to equivalence) the objects of $\KL(\K)$
with monads in \K, and then explicitly describe morphisms and 2-cells
between them in terms of \K\ itself.

In the dual case $\EM(\K) = (\KL(\K\op))\op$ we get
\begin{itemize}
\item The objects are the monads in \K,
\item The morphisms are the monad morphisms (same as in $\mnd(\K)$), and
\item The 2-cells $\xymatrix{(A,t)\rtwocell^{(m,\phi)}_{(n,\psi)}
    &(B,s)}$ are 2-cells $m\to sn$
  (which should look ``Kleisli-like'') with some compatibility with $t$.
\end{itemize}
Composition is like in the Kleisli category.  Think of $sn$ as the
``free $s$-algebra on $n$'', so using the universal property of free
algebras, can express this as something $sm\to sn$, and express
compatibility that way.

Why is this a good thing to do?
\begin{enumerate}[(i)]
\item We still have a fully faithful inclusion $\id:\K\to\EM(\K)$, 
  and by general nonsense for limit-completions, a right adjoint to
  $\id$ is just a choice of EM-objects in \K.
\item It comes up in examples.  If we start with \Span, we've seen
  that categories are just monads in \Span, and that functors can be
  seen as special morphisms between such monads; now we can also deal
  with natural transformations. There is a 2-functor
  \[\Cat \to \KL(\Span)\]
  which is bijective-on-objects and locally fully faithful, so that 
  $\KL(\Span)$ captures precisely the notion of natural transformation.
  This works equally well for \CatE, for \VCat, or for
  generalized multicategories.
\item Remember that a distributive law is a monad in
  $\mnd(\K)$.  The multiplication and unit are 2-cells in $\mnd{\K}$, 
  so if we change the 2-cells, the notion of monad changes.  A monad in
  $\EM(\K)$ is more general: we call it a \emph{wreath}, since the
  composition operation is a wreath product.
\end{enumerate}

A wreath still lives on a monad $(A,t)$ in \K.  We have an endomorphism
$s:A\to A$ as before, along with a 2-cell $\lambda\maps ts\to st$ with some
conditions as before, but $s$ is no longer a monad: the
multiplication is now something $\nu:ss \to st$, and the unit $\sigma:1\to
st$.  You can still make sense of associativity and unit using
$\lambda$, but everything ends up in $st$.  Ultimately this gives a
monad structure on $st$, which is called the \emph{wreath product} or
{\em composite} of $s$ and $t$.

For example, consider the monoidal category \Set under cartesian
product. This can be regarded as a one-object bicategory, and so,
after strictification, as a 2-category. Let $G$ be a group acting
on an abelian group $A$, and consider a normalized 2-cocycle
$G\times G \too[\rho] A$.  We consider $A$ as a monoid (in \Set).
$G$ happens also to
be a monoid (in fact a group), but the monoid structure isn't used 
directly. Rather we have the action
\[\lambda\maps G\times A \to A\times G\]
\[ (g,a) \mapsto ({}^ga,g)\]
and the ``$A\t G$-valued multiplication''
\[\nu\maps  G\times G \too A\times G\]
\[ (g,h) \mapsto (\rho(g,h),gh)\]
which gives a wreath, and so induces a monoid structure $A\sd G$ (which
is actually a group).  The multiplication is the usual one coming from
the cocycle.

There's a corresponding thing for Hopf algebras, giving a type of
``twisted smash product''.

\subsection{The module-theoretic approach}
\label{sec:another-point-view}

Here's another point of view. It's particularly suggestive if
we take \K to be the monoidal category (1-object bicategory) \Ab
of abelian groups. Then a monad (monoid) in \Ab is a ring $R$:
the objects of $\EM(\Ab)$ are the rings.

We defined a morphism $(f,\phi)\maps (A,t)\to (B,s)$ in $\EM(\K)$ 
to consist of
a 1-cell $f\maps A\to B$ and a 2-cell $\phi:sf\to ft$ subject to 
two equations. A 1-cell $R\to S$ in $\EM(\Ab)$
consists of an abelian group $M$ and a map $S\ot M\to M\ot R$. Think of 
this as being a bimodule structure on $M\ot R$; the left action is
\[S\ot M\ot R \too M\ot R\ot R \too M\ot R\] and the right action
is the free one, and the conditions on $\phi$ are equivalent to the
bimodule axioms.  Thus the 1-cells are the \emph{right-free bimodules}.  
The 2-cells are the bimodule homomorphisms.

Composition of 1-cells is the ordinary module composition, but because
of the freeness condition, don't need to use any coequalizers.  If we
were to look at $\KL(\K)$, we'd get the \emph{left-free} modules. 

One could also consider arbitrary modules. This is an important 
construction, but it requires the bicategory to have coequalizers
in the hom-categories in order to define composition; and these
coequalizers to be preserved by whiskering on either side in order
for this composition to be associative (up to isomorphism), and so
this has rather a different flavour.

\subsection{References to the literature}

The formal theory of monads goes back of course to \cite{ftm}; for
the account using limit-completions, and the notion of wreath see
the much later sequel \cite{ftm2}. Distributive laws (for ordinary
monads) are due to Beck \cite{Beck:dist-law}. The Eilenberg-Moore object was
described as a lax limit of a lax functor in ``Two constructions on 
lax functors'' \cite{Street-twoconstructions}: it was the first
construction, the Kleisli object was the second. This was done
using weighted limits in \cite{Street-Cat-limits}.
For limits in $\textrm{T-Alg}_l$, including Eilenberg-Moore objects 
for comonads, see \cite{talgl}; for Hopf monads see 
\cite{Moerdijk-Hopf-monads} and also \cite{McCrudden-opmonoidal}.

\section{Pseudomonads}
\label{sect:pseudomonads}

These are formally very similar to monoidal categories.  A pseudomonad 
involves a thing
$T$, which plays the role of a category, a multiplication $m:T^2\to T$, 
a unit $i:1\to T$, an associativity isomorphism 
\[\xymatrix{T^3 \ar[r]^{mT} \ar[d]_{Tm} \ar@{}[dr]|{\cong} & 
T^2 \ar[d]^m \\
  T^2\ar[r]_m & T}\]
unit isomorphisms $\lambda, \rho$, and so on, all looking very like a
monoidal category.
Just as monads can be defined in any 2-category or bicategory, pseudomonads
can be defined in any Gray-category or tricategory. 

The monoidal 2-category
\Cat (with cartesian structure) can be regarded as a one-object tricategory,
and a pseudomonad in this tricategory is precisely a monoidal category.
The associativity pentagon becomes a cube, relating ways to go from $T^4$
to $T$, involving a bunch of $\mu$'s and a pseudonaturality isomorphism.
In monoidal categories, one side of the cube corresponds to
\[\xymatrix{((A\ot B)\ot C)\ot D \ar[d]\\
  (A \ot (B\ot C))\ot D \ar[dr] &&
  A\ot (B\ot (C\ot D))\\
  & A \ot ((B\ot C) \ot D)\ar[ur]
}\]
and the other side corresponds to
\[\xymatrix{
  & (A\ot B) \ot (C\ot D) \ar@{=}[dr] & \\
  ((A\ot B)\ot C)\ot D \ar[ur] &&   (A\ot B) \ot (C\ot D) \ar[d]\\
  && A\ot (B\ot (C\ot D))
}\]
where in general, the equality will be replaced by an isomorphism,
saying that it doesn't matter whether we tensor $A$ and $B$ first,
then $C$ and $D$, or vice versa.

Our unit isomorphisms have the form
\[\xymatrix{T \ar[r]^{iT} \ar[d]_{Ti} \ar[dr]_(0.3){1} 
\drtwocell\omit{<-2>} \drtwocell\omit{^<2>} &
 T^2 \ar[d]^{m} \\
 T^2\ar[r]_{m} & T.}\]
If we were to consider lax monads, so that $\lambda$ and $\rho$ were
not necessarily invertible, their direction would change, but for
coherence problems it's useful to look at $\mu$, $\lambda$, and $\rho$
as rewriting rules, and then one wants them to go from the more 
complicated expression to the simpler.

It is convenient to work with Gray-categories rather than tricategories;
by the coherence result that every tricategory is triequivalent to a 
Gray-category there is no loss of generality. Note, however, that 
the one-object tricategory corresponding to \Cat is not a Gray-category,
although it is a very special sort of tricategory.

One reason for working with Gray-categories is that we can then make
use of the substantial machinery developed for enriched categories.

\subsection{Coherence}
\label{sec:coherence}

The coherence result describes the fact that \textbf{there's a
universal Gray-category with a pseudomonad in it}: there's a
Gray-category \Psm such that for any Gray-category \aA, to give a Gray-functor
$\Psm\to \aA$ is equivalent to giving a pseudomonad in \aA. 
Corresponding to the identity Gray-functor $\Psm\to\Psm$ there is a
pseudomonad in \Psm, and this is the universal pseudomonad.

\Psm is a sort of cofibrant replacement of \mnd. More precisely,
\Psm like \mnd has a single object $*$, and $\Psm(*,*)$  is a 
cofibrant replacement of $\mnd(*,*)$. It's not a pseudomorphism
classifier: that would be too large; we need a smaller cofibrant replacement.
Recall that $\mnd(*,*)=\ordf$, the category of finite ordinals,
or ``algebraists' simplicial category''.
The underlying category of $\Psm(*,*)$ (which is a 2-category, since \Psm is a
Gray-category) is freely generated by the face and degeneracy maps in
\ordf\ (forget the relations we expect to hold)
\[\xymatrix{\ar[r] &
  \ar@<-2mm>[r]  \ar@<2mm>[r] & \dots \ar[l]
}\]

Since this graph  $G$ generates \ordf, we have a map $FG\to \ordf$ which
is bijective on objects and surjective on objects, so we can factor it
as a bijective-on-objects-bijective-on-arrows 2-functor followed 
by an locally-fully-faithful one (throw in isomorphisms between the things 
that would become equal in \ordf), to get
\[\xymatrix{FG \ar[rr] \ar[dr]&& \ordf \\
  & \Psm(*,*) \ar[ur]}\]
To construct the pseudofunctor classifier $\ordf'$, we would forgot all the 
way down to the underlying graph of \ordf, rather than a \emph{generating} graph 
for it, and that would give a much larger cofibrant resolution, including,
for example, a generating operation $T^n\to T$ for any $n$

We also saw  something like this for the Gray tensor product, which
was obtained by factorizing the map from the funny tensor to the ordinary one.

You now have to define the composition in \Psm
\[\mathbf{Psm}(*,*)\ot \mathbf{Psm}(*,*) \too \mathbf{Psm}(*,*)\]
to make it a Gray-category.  You basically take the composition in
\ordf, use that to define it on the generators, then build it up to deal
with arbitrary 1-cells, but since the relations only hold up to
isomorphism, that's why the Gray-tensor appears.

Now you prove that this has the universal property that I said it does, so
it really does classify pseudo-monads in a Gray-category.  I'm
certainly not going to do that.  Roughly, how does it go?  Given a
pseudomonad, we have
\[\xymatrix{1 \ar[r]^i & T
  \ar@<-1pc>[r]^{Ti}  \ar@<1pc>[r]^{iT} & \ar[l]_{m} T^2 & \dots
}\]
and so on, which defines the putative Gray-functor $\Psm\to\aA$ on 
objects, 1-cells, and 2-cells. The fun starts when we come to the
3-cells: we have $\mu$, $\lambda$, and $\rho$, and we need to build up 
all the other required 3-cells. The idea is that for any 2-cell $f$ in 
\Psm (any 1-cell in the above picture, generated by $m$'s and $i$'s), 
there is a normal form $\fbar$ and a unique isomorphism $f\cong \fbar$ 
built up out of the 3-cells in \Psm that one might expect to call $\mu$,
$\lambda$, and $\rho$. Thus any 3-cell $f\cong g$ in \Psm can be 
written as a composite $f\cong\fbar=\gbar\cong g$, and this can be used
to define the Gray-functor $\Psm\to\aA$ on a 3-cells. The details of
the rewrite system that these normal forms come from are a bit technical.

\subsection{Algebras}
\label{sec:algebras}

The next step is to construct a particular weight $\Psa:\Psm\to\Gray$ 
such that for any Gray-functor $\TT:\Psm\to\aA$, 
the weighted limit $\{\Psa,\TT\}$ is the object of
pseudoalgebras, pseudomorphisms, and algebra 2-cells (all suitably
defined) for the pseudomonad corresponding to \TT.  Again, this is
sort of a ``cofibrant replacement'' for the corresponding one for
2-categories, although the domain has changed.

I won't do this, but I do want to make one point. It is the fact that we are 
working with Gray-categories rather than 3-categories which 
causes the pseudomorphisms to appear here. Recall that for ordinary
monads, we talked about the fact that to give something $C\to A^t$ is
the same as $a\maps C\to A$ with an action $\alpha:ta\to a$, where
$c\mapsto (\alpha c:tac\to ac)$, and $\gamma\maps c\to d$ is sent to
\[\xymatrix{tac \ar[r]^{\alpha c}\ar[d]_{ta\gamma} & ac \ar[d]^{a\gamma}\\
  tad\ar[r]_{\alpha d} & ad}\]
and the fact that $a\gamma$ is a homomorphism can be seen as the
naturality of $\alpha$. There's an analogous fact for operads
and Lawvere theories: the actions are natural with respect to
homomorphisms.

When we come up to the Gray situation, we are thinking of
pseudonatural transformations, hence the square commutes up to
isomorphism, so we get pseudo-morphisms, not strict ones.  That's the
``reason'' for making the formal theory of pseudo-monads live in the
Gray context.  Even if you wanted only to consider 3-categories \aA,
the fact of working over \textbf{Gray} gives you the pseudomorphisms.

\subsection{References to the literature}

The basic definitions involving pseudomonads in Gray-categories were given 
in \cite{Marmolejo-Dist1, Marmolejo-KZ}, or in the equivalent language of 
pseudomonoids in \cite{DS:Hopf-algebroid}; for the universal pseudomonad, and 
the Gray-limit approach to pseudoalgebras, see \cite{psm}. For further
work see also \cite{Cheng-Hyland-Power,Street-frobenius,Lauda-frobenius}.

\section{Nerves}
\label{sect:nerves}

In this section we use \DD for the ``topologists' delta'', the
category of {\em non-empty} finite ordinals. As usual, we write
$[n]$ for the ordinal $\{0<1<\ldots<n\}$. This section is particularly
light on details; see \cite{2-nerve} for more.

The nerve of an ordinary category \C is the simplicial set
$N\C$ in which
\begin{itemize}
\item a 0-simplex is an object
\item a 1-simplex is a morphism
\item a 2-simplex is a composable pair and its composite
\end{itemize}
and so on. This process gives a fully faithful embedding 
\[\Cat\hookrightarrow[\DD\op,\Set]=\SSet\]
into the category of simplicial sets.

The nerve of a {\em bicategory} \B is the simplicial set $N\B$ in which
\begin{itemize}
\item a 0-simplex is an object
\item a 1-simplex is a morphism
\item a 2-simplex is a 2-cell living in a triangle
  
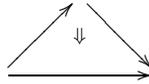
\[\xymatrix{ & \ar[dr] \ar@{}[d]|{\Downarrow}  \\ 
\ar[rr] \ar[ur] && }\] 
\end{itemize}
and so on.
These 2-simplices are being overworked; they have to express both 
(something about) composition of 1-cells, as well as what
the 2-cells are. Actually they don't ever encode what the composite of
two 1-cells is, only what the morphisms out of the composite are.
(This is like defining the tensor product of modules only in terms
of bilinear maps, never actually specifying a choice of universal
bilinear morphism.) Now this has its advantages, but it does make it 
hard to say when such composites are being preserved. In fact we get 
a fully faithful embedding
\[\bicat_{\mathrm{nlax}} \hookrightarrow [\DD\op,\Set]\]
where ``nlax'' indicates that we are taking the  {\em normal lax functors}
as morphisms: these are the lax functors which strictly preserve identities.

A lot of the time you want to talk about homomorphisms or strict
homomorphisms rather than lax
ones. If you want to get your hands on those there are various
possibilities.  One is to have a bit more structure than a simplicial
set: specify as extra data a class of simplices called the {\em thin}
simplices. Every degenerate simplex is thin, but there can also be 
non-degenerate
thin simplices. The resulting structure is called a {\em stratified simplicial
set}. The {\em stratified nerve} of a 2-category is the usual nerve,
made equipped with a suitable stratification, in which a 2-simplex
is thin when the 2-cell it contains is in fact an identity. (There is
also a different stratification, often used for nerves of bicategories,
in which a 2-simplex is then when the 2-cell is invertible.) One can
now characterize the stratified simplicial sets which are stratified
nerves of 2-categories (or bicategories); and indeed similarly for 
strict or weak $\omega$-categories. The stratified simplicial sets
which arise as nerves of $\omega$-categories are called {\em complicial
sets}.

A different way to specify this extra structure is to use simplicial
objects not in \Set but in \Cat. For a bicategory \B, the 2-nerve
$N_2\B$ of \B (or just $N\B$ from now on) is a functor 
$N\B:\DD\op\to\Cat$ 
\begin{itemize}
\item the category $(N\B)_0$ of 0-simplices is the discrete category 
consisting of the objects of \B.
\item the category $(N\B)_1$ of 1-simplices has morphisms of \B as objects
and 2-cells of \B as morphisms. So far this looks like some kind of 
enriched nerve.
\item the category $(N\B)_2$ of 2-simplices doesn't need to include the 
2-cells, since we already have them; we can therefore take 2-simplices
to be {\em invertible} 2-cells
  
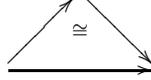
\[\xymatrix{ & \ar[dr] \ar@{}[d]|{\cong}  \\ 
\ar[rr] \ar[ur] && }\] 
and morphisms of 2-simplices to consist of three 2-cells satisfying
the evident compatibility conditions. (The domain and codomain 2-simplices
will need to have the same three objects.)
\item the category $(N\B)_3$ of 3-simplices is the category of 
tetrahedra all of whose faces are isomorphisms.
\end{itemize}
and so on.
We'd like a functorial description of the 2-nerve. Consider \NHom, the
2-category of bicategories, normal homomorphisms, and icons (recall that these
are oplax natural transformations all of whose 1-cell
components are identities).  Now $\Cat\hookrightarrow \NHom$, where \Cat\ 
is the locally discrete 2-category consisting of categories, functors, and
only identity natural transformations, embedding as a full
sub-2-category consisting of the locally discrete bicategories.  (An
icon between functors can only be an identity.)  And of course
we have $\DD\hookrightarrow \Cat$, so the composite fully faithful 
$J\maps\DD\hookrightarrow\NHom$ induces
\[\NHom(J,1)\maps \NHom \too{} [\DD\op,\Cat]\]
sending \B\ to $\NHom(J-,\B)$.

For instance, $J$ sends $[0]\in\DD$ goes to the terminal bicategory, and a 
normal homomorphism from that into \B\ is just an object of \B, with no
room for icons, so $\NHom(J[0],\B)$ is indeed the discrete category 
$(N\B)_0$ of objects of \B. Similarly $J$ sends $[1]\in \DD$ goes to the 
arrow category $\two$, and a normal homomorphism from this into \B\ is an 
arrow in \B, and $\NHom(J[1],\B)$ is just $(N\B)_1$, and so on.

\begin{theorem}
  $\NHom(J,1) = N$ is a fully faithful 2-functor (in a completely
  strict sense) and has a left biadjoint.
\end{theorem}

The fact that $N=\NHom(J,1)$ is a straightforward direct calculation.
The existence of the left biadjoint can be proved using techniques
of 2-dimensional universal algebra.

How can we characterize the image of $N$?  $X\in [\DD\op,\Cat]$ is isomorphic
to some $N\B$ if and only if
\begin{enumerate}[(a)]
\item $X_0$ is discrete;
\item $X$ is \emph{3-coskeletal}; that is, isomorphic to the right Kan
  extension of its 3-truncation --- the idea is that 4-simplices and
  higher are uniquely determined by their boundary;
\item $X_2 \to \cosk_1(X)_2$ is a \emph{discrete isofibration}.  A
  functor $p\maps A\to B$ is a discrete isofibration if given $e\in E$
  and $\beta\maps b\cong pe$, there exists a \emph{unique} $\epsilon\maps
  e'\cong e$ with $p\epsilon = \beta$.  This implies that if
  \[\xymatrix{X\rtwocell{\epsilon}& E}\]
  and $p\epsilon = \id$, then $\epsilon = \id$;
\item $X_3 \to \cosk_1(X)_3$ (could also use the 2-coskeleton) is also
  a discrete isofibration
\item The Segal maps are equivalences.
\end{enumerate}

A {\em Tamsamani weak 2-category}, or just {\em Tamsamani 2-category},
since no strict notion is considered, is a functor $X:\DD\op\to\Cat$ 
satisfying (a) and (d); thus the 2-nerve of a bicategory is a Tamsamani
2-category. Tamsamani suggests a way of getting from a bicategory 
to a Tamsamani 2-category, but it is not the 2-nerve construction
given here.

The inclusion of \NHom into Tamsamani 2-categories looks like it
should be a biequivalence, but it's not quite. It would be if you 
broadened the definition of morphism of Tamsamani 2-category to include
what might be called the ``normal pseudonatural transformations''.

Finally a warning:
what you might guess for the nerve of a bicategory is to have
\begin{itemize}
\item $N\B_0$ the objects
\item $N\B_1 = \sum_{x,y} \B(x,y)$
\item $N\B_2 = \sum_{x,y,z} \B(y,z)\times \B(x,y)$
\item $N\B_3 = \sum_{w,x,y,z} \B(y,z)\times \B(x,y)\t\B(w,x)$
\end{itemize}
in analogy with the case of nerves of ordinary categories.
If you try to do this, the
simplicial identities fail, due to the failure of associativity.
Actually, what we do is to take the pseudo-limit of
the composition functor 
\[ \sum_{x,y,z} \B(y,z)\times \B(x,y) \too\sum_{x,y} \B(x,y)\]
as $(N\B)_2$. And this continues: for composable triples, we have
\[\xymatrix{\B^3 \ar[r]\ar[d] \ar@{}[dr]|{\cong} & \B^2 \ar[d]\\
  \B^2\ar[r] & \B}\]
and $N\B_3$ is the pseudo-limit of this whole diagram.  Going on,
we can construct each $(N\B)_n$ as the pseudolimit of some 
higher cube.

\subsection{References to the literature}

The notion of nerve of a bicategory is due to Street. For nerves of
$\omega$-categories as stratified simplicial sets, see 
\cite{Verity-complicial}, and the references therein. The notion of 
2-nerve of a bicategory is described in \cite{2-nerve}.
Tamsamani's definition of weak n-category is in \cite{Tamsamani}.

\subsection*{Acknowledgements}

It is a pleasure to acknowledge support and encouragement from a number
of sources. I am grateful to the Institute for Mathematics and its
Applications, Minneapolis for hosting and supporting the workshop on
higher categories in 2004, and to John Baez and Peter May who organized the 
workshop and who encouraged me to publish these notes. The material here
was based on lectures I gave at the University of Chicago in 2006, at 
the invitation of Peter May and Eugenia Cheng. I'm grateful to them
for their hospitality, and the interest that they and the
topology/categories group at Chicago took in these lectures. I'm 
particularly grateful to Mike Shulman, whose excellent TeXed notes of the
lectures were the basis for the companion. 

%I also acknowledge with gratitude the financial support of the 
%Australian Research Council and DETYA.

\bibliographystyle{plain}
%\bibliography{my}

\end{document}